\documentclass{l4dc2027}
\usepackage{times}
\newtheorem{mainthm}{Theorem}
\newtheorem{assumption}{Assumption}
\usepackage{booktabs}
\usepackage{tikz-cd}
\usepackage{caption}
\newcommand{\R}{\mathbb{R}}
\newcommand{\g}{\mathfrak{g}}
\newcommand{\so}{\mathfrak{so}}
\newcommand{\se}{\mathfrak{se}}
\newcommand{\SO}{\mathsf{SO}}
\newcommand{\SE}{\mathsf{SE}}
\newcommand{\Exp}{\operatorname{Exp}}
\newcommand{\Log}{\operatorname{Log}}
\newcommand{\dist}{\operatorname{dist}}
\newcommand{\vecop}{\operatorname{vec}}
\newcommand{\tr}{\operatorname{tr}}
\newcommand{\Id}{\mathrm{I}}
\newcommand{\amax}{\alpha_{\max}}
\newcommand{\spr}{\operatorname{spr}}
\newcommand{\Vm}{V_{\mathrm{m}}}
\newcommand{\Wm}{W_{\mathrm{m}}}
\newcommand{\Pm}{P_{\mathrm{m}}}
\newcommand{\calA}{\mathcal{A}}
\newcommand{\calB}{\mathcal{B}}
\newcommand{\calC}{\mathcal{C}}
\newcommand{\calX}{\mathcal{X}}
\newcommand{\calK}{\mathcal{K}}
\newcommand{\calN}{\mathcal{N}}
\newcommand{\calZ}{\mathcal{Z}}
\newcommand{\calL}{\mathcal{L}}
\newcommand{\calR}{\mathcal{R}}
\newcommand{\rot}{\vartheta}
\newcommand{\rG}{r_{\mathsf{G}}}
\newcommand{\Dch}{\mathcal{D}}
\newcommand{\Uch}{\mathcal{U}}
\newcommand{\Cl}{C_{\ell}}
\newcommand{\Lip}{\operatorname{Lip}}
\newcommand{\sat}{\operatorname{sat}}
\newcommand{\proj}{\Pi_{\amax}}

\title[Maximal Lyapunov Functions for RL Controllers on Lie Groups]{Learning Maximal Lyapunov Functions for Deep Reinforcement Learning-based Controllers on Matrix Lie Groups}

\author{\Name{Hamza Tariq} \Email{ht299@njit.edu}\\ \addr New Jersey Institute of Technology, Newark, NJ, USA.\AND \Name{Matthieu Barreau} \Email{barreau@kth.se}\\ \addr Digital Futures and KTH Royal Institute of Technology, Stockholm, Sweden.\AND \Name{Adeel Akhtar} \Email{adeel.akhtar@njit.edu}\\ \addr New Jersey Institute of Technology, Newark, NJ, USA.}

\begin{document}

\maketitle

\vspace{-0.36in}

\begin{abstract}

\looseness -1 Deep reinforcement learning (RL) is effective at generating policies for orientation and pose control on matrix Lie groups, but gives no estimate of their region of attraction (ROA). We estimate this ROA by a sublevel set of a maximal Lyapunov function, the undiscounted sum of squared distances to the goal, which a neural network learns from its one-step equation. In practice, an actuator may scale its command or act along rotated axes, which can invalidate this estimate. Once this error is identified, we recompute the region for the same policy, without retraining. The regions are verified by interval arithmetic on $\SO(2)$ and checked by sampling and adversarial search on $\SO(3)$, $\SE(2)$ and $\SE(3)$. On 39 policies trained with four RL algorithms, they match or exceed the regions of squared-norm and linearization-based Lyapunov functions. With actuator errors, they exceed both on all 55 closed loops and those of a neural Lyapunov baseline on 49. In MuJoCo simulations of a CubeSat and an underwater vehicle, no trajectory that starts in an estimated region leaves it, including the regions recomputed after actuator errors. Code and animations are available at \url{https://gradslab.github.io/maxlyap_lie_rl}.
\end{abstract}


\section{Introduction}
\looseness -1 Rotations and rigid motions, which describe the orientation and pose of robots, can be represented as elements of matrix Lie groups. Exploiting the symmetry of these groups accelerates the training of reinforcement learning (RL) policies and reduces their tracking error~\citep{welde2025leveraging,jaquier2021geometry}. Likewise, on the rotation group $\SO(3)$, bounded body-frame tangent vectors as actions give faster and more reliable training than quaternions, Euler angles, or rotation matrices~\citep{schuck2025lie,schuck2026primer}. However, training does not establish the set of initial states from which the policy reaches the goal, that is, the region of attraction (ROA) of the closed loop. Moreover, an actuator on the hardware may deliver a scaled version of its command, called a gain error~\citep{xiao2013attitude}, or act along axes rotated away from the commanded ones, called a misalignment. Such actuator errors increase the gap between simulation and hardware~\citep{zhao2020sim}, so a ROA estimated in simulation may not hold on the hardware. Indeed, once such an error is identified, the same policy has a different closed loop and possibly a different ROA. Hence, there is a need to estimate the ROA of a trained policy, and to recompute it after such an error to decide whether the policy can be kept.

\looseness -1 Lyapunov methods compute a ROA as a sublevel set of a Lyapunov function. For complex dynamics, a Lyapunov function in closed form is difficult if not impossible to derive. Therefore, neural Lyapunov functions are learned and formally verified jointly with the controller~\citep{chang2019neural,dai2021lyapunov,wu2023neural,yang2024lyapunov,shi2026certified}. However, the controller is then obtained together with its certificate, so this approach does not apply to a policy that has already been trained. In contrast, a certificate can be learned for the closed loop of a fixed controller~\citep{richards2018lyapunov,long2025certifying,liu2025physics}. This separates certification from training, so a policy trained by any RL algorithm can be analyzed, and the analysis can be repeated when the dynamics change. We follow this second approach. For an RL policy, a natural candidate for such a certificate is its value function, but a discounted value function is a Lyapunov function only under extra conditions~\citep{postoyan2017stability,granzotto2021finite}. \citet{long2025certifying} therefore add a learned correction to the value function and require only that the sum decrease on average over several steps. However, for pretrained policies, they do not estimate a ROA. Moreover, their systems evolve in $\R^{n}$, and matrix Lie groups are not addressed.

Another research direction uses maximal Lyapunov functions, which are finite exactly on the ROA~\citep{vannelli1985maximal} and are characterized by a Zubov-type equation~\citep{zubov1964methods}. For discrete-time systems on $\R^{n}$, the sum of a running cost along trajectories is a maximal Lyapunov function on the whole ROA~\citep{giesl2007determination}. \citet{serry2025safe} characterize safe domains of attraction, subsets of the ROA, by this sum and a discrete-time Zubov equation, with a verified neural approximation. On Lie groups, \citet{akhtar2026learning} learn neural maximal Lyapunov functions for known continuous-time dynamics on $\SO(n)\times\R^{n(n-1)/2}$ from a Zubov-type equation. In continuous time, neural candidates built on local Taylor terms or homogeneous Lyapunov functions have also been proposed~\citep{barreau2024unsupervised,liu2026universal}. However, none of these works learns a maximal Lyapunov function for the discrete-time closed loop of a trained RL policy on a matrix Lie group, or repeats the analysis when the applied actuator commands change. This is difficult because a trained policy does not make the goal an exact fixed point, and the goal need not be stable. In addition, near the goal, the decrease of the function over one step is too small to learn.

In this work, our key idea is to learn the trained policy's own accumulated cost, the undiscounted sum of squared distances to the goal. Provided the goal is a locally exponentially stable fixed point, this accumulated cost is a maximal Lyapunov function of the closed loop. It also satisfies an undiscounted Bellman equation over one step: its value is the cost of this step plus its value at the next state. Hence, a network can learn it, and its Lyapunov conditions can be checked, from single transitions rather than from long simulations or a partial differential equation. Near the goal, where the decrease is too small to learn, we use the exact quadratic term from the Lyapunov equation for the linearized system. Moreover, unlike in \citet{schuck2026primer}, our policy observes only the state-to-goal error, so one region serves every goal. Finally, when an actuator error is identified, we recompute the region for the same policy on the changed closed loop. To reach this objective, we make the following contributions:
\begin{itemize}\setlength{\itemsep}{0pt}
\item We consider RL policies on matrix Lie groups that steer the state to a goal and observe only the state-to-goal error, and show that a single closed loop then covers every goal~(Proposition~\ref{thm:descent}).
\item We approximate the maximal Lyapunov function by a neural network and give conditions under which a sublevel set of the learned function lies inside the ROA (Theorem~\ref{thm:main}). We also prove that a network satisfying these conditions exists, and we bound the approximation error by the one-step residual (Corollary~\ref{cor:exists} in Appendix~\ref{app:proofs} and Proposition~\ref{prop:residual}).
\item We recompute the region after an identified actuator gain error or misalignment, without retraining the policy. The regions match or exceed those of baseline Lyapunov functions on almost all evaluated closed loops, with and without actuator errors.
\item We test the regions in MuJoCo simulations of a CubeSat and an underwater vehicle, in which no trajectory leaves its estimated region, also after actuator errors on the underwater vehicle.
\end{itemize}

Sections~\ref{sec:background} to~\ref{sec:sims} give the RL problem, the estimation of the ROA, the guarantees and the experiments. The appendices give the proofs and further details.

\vspace{-0.075in}

\paragraph{Notation.} We write $\|\cdot\|_F$ for the Frobenius norm, $\spr(M)$ for the spectral radius of $M$, $(x)_{+} = \max\{x, 0\}$ for the positive part of a real number $x$, $\Phi^{k}$ for $\Phi$ applied $k$ times, and $B_r$ and $\bar B_r$ for the open and closed Euclidean balls of radius $r$ about the origin $0$. We write $\|\cdot\|$ for the Euclidean norm, and $O(\|z\|^{k})$ for a function bounded by a constant times $\|z\|^{k}$ near $0$.

\vspace{-0.1in}

\section{Reinforcement learning on matrix Lie groups}
\label{sec:background}
\vspace{-0.05in}

\subsection{The system and exponential coordinates}
\looseness -1 The orientation of a rigid body is a rotation matrix, and its pose is a rotation together with a position. We focus here on the following matrix Lie groups
\begin{equation}
  \SO(n) = \{R \in \R^{n\times n} : R^{\top}R = \Id,\ \det R = 1\}, \quad
  \SE(n) = \bigl\{\bigl(\begin{smallmatrix} R & p \\ 0 & 1 \end{smallmatrix}\bigr) : R \in \SO(n),\ p \in \R^{n}\bigr\},
\end{equation}
with $n = 2$ or $3$, since the states of most robotic systems lie on them~\citep{akhtar2022hybrid,akhtar2020controller}. More generally, let $\mathsf{G}$ be a matrix Lie group, that is, a closed subgroup of the group of invertible $d \times d$ matrices~\citep[Definition 1.4]{hall2015lie}, and let $\g$ be its Lie algebra, of dimension $m$. We fix a basis $\hat e_1, \dots, \hat e_m$ of $\g$ (Appendix~\ref{app:lie}), and the hat map $\xi \mapsto \hat\xi = \sum_i\xi_i\hat e_i$ assigns to every vector $\xi$ in $\R^{m}$ a matrix $\hat\xi$ in $\g$. We consider the left-invariant kinematic system $\dot g = g\hat\xi$ with state $g$ in $\mathsf{G}$ and input $\xi$ in $\R^{m}$~\citep{murray1994mathematical}. The matrix $\hat\xi = g^{-1}\dot g$ is the body velocity, an element of $\g$, and $\xi$ is its coordinate vector: the body angular velocity on $\SO(n)$, and the body angular and linear velocities on $\SE(n)$. Under a zero-order hold of the input over a sampling period $\Delta t$, the state at the sampling instants satisfies exactly
\begin{equation}
  g_{k+1} = g_k\Exp(\tau_k), \qquad \tau_k = \Delta t\,\xi_k, \qquad \Exp(\tau) = \exp(\hat\tau),
  \label{eq:sys}
\end{equation}
where $\tau_k$ in $\R^{m}$ is the increment, $\exp$ is the matrix exponential, and $\Exp : \R^{m} \to \mathsf{G}$ is the exponential map of $\mathsf{G}$ in the coordinates of the basis. Since $\Exp(\tau_k)$ lies in $\mathsf{G}$, the state stays on the group. We bound the length of every increment, as \citet{schuck2026primer} do on $\SO(3)$, so that one step moves the state by a bounded amount. To this end, every increment is projected onto the closed ball $\bar B_{\amax} \subset \R^{m}$ of radius $\amax > 0$, the clip length, by $\proj(\tau) = \tau$ if $\|\tau\| \le \amax$ and $\proj(\tau) = \amax\,\tau/\|\tau\|$ otherwise. Near the identity, the matrix logarithm $\Log$ inverts $\Exp$ and gives the exponential coordinates $z = \Log g$ in $\R^{m}$~\citep{akhtar2026learning}. More precisely, on every matrix Lie group there are a chart radius $\rG > 0$ and an open set $\Dch \subseteq \R^{m}$ containing $B_{\rG}$ on which $\Exp$ is a diffeomorphism onto its image, with inverse $\Log : \Exp(\Dch) \to \Dch$ (Proposition~\ref{prop:general} in Appendix~\ref{app:lie}). On $\SO(n)$ and $\SE(n)$, $\rG = \pi$ (Proposition~\ref{prop:assumption}). Moreover, the clip length satisfies $0 < \amax < \rG$.

\vspace{-0.1in}

\subsection{The error state, the learning task, and the closed loop}
\label{sec:loopdef}
\looseness -1 The control objective is to steer the state $g \in \mathsf G$ of~\eqref{eq:sys} to a goal $\bar g \in \mathsf G$. The error is $\tilde g = \bar g^{-1}g$, which equals the identity $\Id$ exactly at the goal. Multiplying \eqref{eq:sys} on the left by $\bar g^{-1}$ gives us
\begin{equation}
  \tilde g_{k+1} = \tilde g_k\Exp(\tau_k), \qquad \tilde g_k = \bar g^{-1}g_k .
  \label{eq:errdyn}
\end{equation}
The goal does not appear in \eqref{eq:errdyn}, so the learning task is posed on the error, with the goal fixed at $\Id$. The policy observes $\vecop\tilde g$, where $\vecop : \mathsf{G} \to \R^{n_o}$ stacks the matrix entries that are not constant on $\mathsf{G}$: all $n^{2}$ entries on $\SO(n)$, and the $n^{2} + n$ entries of $R$ and $p$ on $\SE(n)$. As is common in RL~\citep{haarnoja2018sac}, $\sat$ saturates each entry of an action $a$ in $\R^{m}$ to $[-1, 1]$, and $a$ is applied as the increment $\tau = \proj(\amax\sat(a))$. With the cost $\phi(\tilde g) = \|\Log\tilde g\|$ and a scale $\bar\phi > 0$, the reward is $-\phi(\tilde g_{k+1})/\bar\phi$. The initial error is drawn from a Borel probability measure $\mu_0$ on $\mathsf{G}$ with compact support, and returns are discounted by $\beta$ (Appendix~\ref{app:policies}). An RL algorithm trains a neural network $\mu : \R^{n_o} \to \R^{m}$, which is the mean of a stochastic policy or the deterministic policy itself. The trained policy is $u : \mathsf{G} \to \R^{m}$ with $u(\tilde g) = \amax\sat(\mu(\vecop\tilde g))$, and its closed loop on the errors is
\begin{equation}
  \tilde g_{k+1} = \Phi_u(\tilde g_k), \qquad \Phi_u : \mathsf{G} \to \mathsf{G}, \quad \Phi_u(\tilde g) = \tilde g\,\Exp\bigl(\proj(u(\tilde g))\bigr).
  \label{eq:closedloop}
\end{equation}
\begin{proposition}
\label{thm:descent}
Let $\mathsf{G}$ be a matrix Lie group and $u : \mathsf{G} \to \R^{m}$ a map. For every goal $\bar g \in \mathsf{G}$ and initial state $g_0 \in \mathsf{G}$, the states $g_{k+1} = g_k\Exp(\proj(u(\bar g^{-1}g_k)))$ satisfy $g_k = \bar g\tilde g_k$, with $\tilde g_k$ given by \eqref{eq:closedloop} from $\tilde g_0 = \bar g^{-1}g_0$, and $g_k \to \bar g$ if and only if $\tilde g_k \to \Id$.
\end{proposition}
\noindent The proof is in Appendix~\ref{app:proofs}. By Proposition~\ref{thm:descent}, the closed loop for every goal $\bar g$ is conjugate to \eqref{eq:closedloop} by the left translation $g \mapsto \bar g g$, so it suffices to study \eqref{eq:closedloop}. Let $E \in \R^{m\times m}$ be the matrix of an identified actuator error, diagonal for a gain error and orthogonal for a misalignment. The policy is unchanged, but its increment is multiplied by $E$ and projected again, which gives the closed loop
\begin{equation}
  \Phi_{u_E}(\tilde g) = \tilde g\,\Exp\bigl(\proj(u_E(\tilde g))\bigr), \qquad u_E(\tilde g) = E\,\proj(u(\tilde g)).
  \label{eq:closedloopE}
\end{equation}
This is \eqref{eq:closedloop} with $u_E$ in place of $u$, so the results below apply to it whenever their hypotheses hold for $u_E$. We study \eqref{eq:closedloop} in exponential coordinates on a closed ball $\bar B_\gamma \subset \R^{m}$ with $\gamma > 0$ and $\gamma + \amax < \rG$. One step from this ball cannot leave $\Exp(\Dch)$ (Proposition~\ref{prop:domain} in Appendix~\ref{app:proofs}), so the closed loop in coordinates is the continuous map
\begin{equation}
  z_{k+1} = F(z_k), \qquad F = \Log{}\circ\Phi_u\circ\Exp : \bar B_\gamma \to \Dch .
  \label{eq:F}
\end{equation}
The region of attraction of the goal~\citep[Remark 6.22]{bullo2005geometric}, and its part in coordinates whose trajectories $z_{k+1} = F(z_k)$ stay in the ball, are
\begin{equation}
  \calA = \bigl\{\tilde g \in \mathsf{G} : \lim_{k\to\infty}\Phi_u^{k}(\tilde g) = \Id\bigr\}, \quad
  \calA_\gamma = \{z_0 \in \bar B_\gamma : z_k \in \bar B_\gamma \text{ for every } k \ge 0 \text{ and } z_k \to 0\}.
\end{equation}

\paragraph{The ROA estimation task.} We estimate $\calA_\gamma$ by a subset $\Omega \subseteq \calA_\gamma$, called an inner estimate of $\calA_\gamma$. By Proposition~\ref{thm:descent}, the policy reaches every goal $\bar g$ from every state in $\bar g\Exp\Omega$. We take $\Omega$ to be a sublevel set of a function that we learn, that is, $\Omega = \{z \in \bar B_\gamma : V(z) \le \rho\}$ with $V : \Dch \to [0, \infty)$ continuous and a level $\rho > 0$. Such a set is always compact. We seek $V$ and $\rho$ for which $\Omega$ is an inner estimate with two further properties. First, it is positively invariant, that is, $F(\Omega) \subseteq \Omega$, so trajectories that start in $\Omega$ stay in it. Second, it contains a neighborhood of $0$. We call $\Omega$, and its image $\Exp\Omega$ on the group, the region. Among the regions with these properties, we want one that is as large as possible. Recall that the learning task draws each initial error at random, from the probability measure $\mu_0$ on $\mathsf{G}$. We measure the size of a region by the probability $\mu_0(\Exp\Omega)$ that such an initial error lies in it, that is, by the fraction of these initial errors that the region covers. Here $u$ is the trained policy after a correction that makes the goal a fixed point (Section~\ref{sec:loop}). After an identified actuator error, we seek a region with the same properties for the closed loop \eqref{eq:closedloopE}.

\vspace{-0.12in}

\section{Estimating the region of attraction}
\vspace{-0.05in}

\subsection{Fixed point and stability test}
\label{sec:loop}
\looseness -1 Our construction requires that $\Phi_u$ keeps the goal fixed and that the goal is locally exponentially stable, and training for return guarantees neither. The closed loop keeps the goal fixed only if $\mu(\vecop\Id) = 0$, which training does not enforce. We therefore subtract the constant offset $\mu(\vecop\Id)$, as \citet{dai2021lyapunov} and \citet{yang2024lyapunov} do for neural controllers in $\R^{n}$, and replace $u$ by $u(\tilde g) = \amax\,\sat(\mu(\vecop\tilde g) - \mu(\vecop\Id))$, so that $u(\Id) = 0$ and $\Phi_u(\Id) = \Id$. From now on $u$ denotes the corrected policy. It is smooth near $\Id$, because the network has smooth activations (Appendix~\ref{app:policies}) and the saturation and projection are inactive near the goal, so $F$ is smooth near $0$ with $F(0) = 0$.

\looseness -1 The goal is locally exponentially stable if there are $r_0 > 0$, $c_0 \ge 1$ and $q \in [0, 1)$ such that every trajectory of \eqref{eq:F} with $z_0 \in B_{r_0}$ stays in $\bar B_\gamma$ and satisfies $\|z_k\| \le c_0q^{k}\|z_0\|$ for every $k \ge 0$. Since $F$ is smooth near $0$, this holds if and only if $\spr(A) < 1$, where $A = DF(0)$ is the Jacobian of $F$ at $0$ (Lemma~\ref{lem:spr} in Appendix~\ref{app:proofs}). The stability test estimates $A$ by central differences (Appendix~\ref{app:training}) and keeps a closed loop only if the spectral radius of the estimate is at most $0.99$. This numerical test does not prove that $\spr(A) < 1$. Among the evaluated closed loops, the largest estimate is $0.873$ (Appendix~\ref{app:results}). The results below therefore assume that the goal is locally exponentially stable. Under this assumption, the closed loop on $\SO(n)$ has a second fixed point by the Lefschetz fixed-point theorem (Proposition~\ref{prop:lefschetz} in Appendix~\ref{app:proofs}), so the ROA is not the whole group.

\vspace{-0.1in}

\subsection{Maximal Lyapunov function and neural candidate}
We approximate a maximal Lyapunov function of $\Phi_u$ by a learned candidate. In the discrete-time version of the definition of \citet{vannelli1985maximal}, such a function $W : \mathsf{G} \to [0, \infty]$ is finite exactly on $\calA$ and continuous there, zero at $\Id$, positive with $W(\Phi_u(\tilde g)) < W(\tilde g)$ on $\calA \setminus \{\Id\}$, and tends to infinity along every sequence in $\calA$ that converges to a point outside $\calA$.

\vspace{-0.075in}

\paragraph{The accumulated cost.}\label{sec:local} Let $\ell(\tilde g) = \|\tilde g - \Id\|_{F}$ be the distance between the matrices $\tilde g$ and $\Id$. Each step is charged the running cost $\ell(\tilde g)^{2}$, which is smooth on every matrix Lie group, needs no logarithm and vanishes only at $\Id$. The accumulated cost of the policy is $\Vm(\tilde g) = \sum_{k \ge 0} \ell(\Phi_u^{k}(\tilde g))^{2}$, with values in $[0, \infty]$. It is the undiscounted cost-to-go for $\ell^{2}$, not the value function of the training objective, which uses the cost $\phi$ and the discount $\beta$. Its terms are nonnegative, so separating the term $k = 0$ gives the one-step equation
\begin{equation}
  \Vm(\tilde g) = \ell(\tilde g)^{2} + \Vm(\Phi_u(\tilde g)) \qquad \text{for every } \tilde g \text{ in } \mathsf{G} .
  \label{eq:recursion}
\end{equation}

Since the goal is locally exponentially stable, $\Vm$ is finite exactly on $\calA$ and is a maximal Lyapunov function of $\Phi_u$ (Proposition~\ref{prop:maximal} in Appendix~\ref{app:proofs}), as \citet[Theorems 2, 4 and 5]{serry2025safe} show for maps on $\R^{n}$. Unlike \citet{serry2025safe}, who train a network on the discrete-time Zubov equation, we train on \eqref{eq:recursion} directly. The second-order term of the accumulated cost needs no simulation, unlike the cost itself. Near the goal the running cost satisfies $\ell(\Exp z)^{2} = z^{\top}Qz + O(\|z\|^{3})$, with the positive definite matrix $Q_{ij} = \tr(\hat e_i^{\top}\hat e_j)$, which is $2\Id$ on $\SO(n)$ and $\operatorname{diag}(2\Id, \Id)$ on $\SE(n)$ (Appendix~\ref{app:lie}). Since $A = DF(0)$ has $\spr(A) < 1$, the discrete-time Lyapunov equation
\begin{equation}
  A^{\top}\Pm A - \Pm = -Q
  \label{eq:stein}
\end{equation}
has the unique solution $\Pm = \sum_{k \ge 0}(A^{\top})^{k}QA^{k}$ \citep{lasalle1986stability}, symmetric and positive definite because $\Pm \succeq Q \succ 0$. In exponential coordinates the accumulated cost is $\Wm = \Vm\circ\Exp$. Since $u$ is smooth near $\Id$, $\Wm(z) = z^{\top}\Pm z + O(\|z\|^{3})$ as $z \to 0$~(Lemma~\ref{lem:local} in Appendix~\ref{app:proofs}).

\vspace{-0.075in}

\paragraph{The neural candidate.}\label{sec:candidate} \looseness -1 With a network $\psi : \R^{m} \to \R^{n_\psi}$, a coefficient $\kappa > 0$ and a radius $r_a > 0$, the candidate is
\begin{equation}
  V(z) = \bigl(1 - s(z)\bigr)\,z^{\top}\Pm z + s(z)\,V_0(z), \qquad s(z) = \frac{\|z\|^{2}}{\|z\|^{2} + r_a^{2}}, \qquad z \in \Dch ,
  \label{eq:cert}
\end{equation}
with the network part $V_0(z) = \|\psi(z) - \psi(0)\|^{2} + \kappa^{2}\,\ell(\Exp z)^{2}$, the architecture of \citet{akhtar2026learning}, which adapts the Lyapunov-Net of \citet{gaby2022lyapunovnet} to Lie groups. The candidate is continuous, vanishes at $0$ and is positive elsewhere. For any weights, $V(z) = z^{\top}\Pm z + O(\|z\|^{4})$, so $V$ has the same quadratic term as $\Wm$. It therefore decreases strictly at every $z \ne 0$ near the goal (Lemma~\ref{lem:localdec} in Appendix~\ref{app:proofs}), in a neighborhood whose radius is not computed. The quadratic term is needed because near the goal the one-step decrease $\ell^{2}$ goes to zero, while the fitting error of a network trained on the whole ball does not. \citet{giesl2008construction} multiplies the approximation by a Taylor polynomial. \citet{barreau2024unsupervised} and \citet{liu2026universal} fix the leading term of a continuous-time candidate from the linearization or from a homogeneous approximation, and \citet{bechihi2025region} correct it so that the second-order term solves a Lyapunov equation. Here the quadratic term is the second-order term of the accumulated cost itself, from \eqref{eq:stein}.

\vspace{-0.1in}

\subsection{Training and the level rule}
\label{sec:learning}
We fix a threshold $r_b > 0$ and a level $\rho > 0$. Let $\Omega_\rho = \{z \in \bar B_\gamma : V(z) \le \rho\}$ be the sublevel set, $N_b = \{z \in \Omega_\rho : \ell(\Exp z) < r_b\}$ its part near the goal, and $\delta V(z) = V(F(z)) - V(z)$ the one-step change of $V$. Near the goal $V$ and $\delta V$ tend to zero, so a sampled test of decrease has no absolute margin there, and $N_b$ is treated separately. We use three Lyapunov conditions: decrease, $\delta V < 0$ on $\Omega_\rho \setminus N_b$; entry, $F(N_b) \subseteq \Omega_\rho$; and containment, $F(\Omega_\rho) \subseteq \bar B_\gamma$. Conditions restricted to a sublevel set inside a box are those of \citet{yang2024lyapunov} in $\R^{n}$, and the separate ball around the fixed point follows~\citet[Corollary 3.13]{giesl2007determination}. If they hold, then $\Omega_\rho$ is positively invariant and every trajectory from it reaches $N_b$ in finitely many steps (Proposition~\ref{prop:cut} in Appendix~\ref{app:proofs}). If, in addition, every trajectory from $N_b$ converges to $0$, then $\Exp(\Omega_\rho) \subseteq \calA$ (Corollary~\ref{cor:region}). The network is trained toward these conditions with two kinds of loss. The residual of \eqref{eq:recursion} pulls $V$ toward $\Wm$. The other terms push $V$ to satisfy the conditions that the level rule will test.

A counterexample is a point at which a condition fails. With a margin $\zeta$ in $(0, r_b^{2}]$, the loss $\calL_i$ has five terms: the mean squared residual of \eqref{eq:recursion}; a hinge $(\delta V + \zeta)_{+}$ outside $N_b$, with a weight $\lambda$ raised by dual ascent, and with a fixed weight on counterexamples; a shell term and an exit term, which raise $V$ just outside the ball and at points whose next point leaves it; and a penalty on the spectral norms of the weights of $\psi$. Algorithm~\ref{alg:train} summarizes training and refinement, and Algorithm~\ref{alg:train-full} in Appendix~\ref{app:training} gives the full procedure, with constants that are the same for every RL algorithm and, except for the number of iterations, on every group.

\begin{algorithm2e}[!t]
\caption{Training and refinement of the candidate \eqref{eq:cert}, abridged from Algorithm~\ref{alg:train-full}.}
\label{alg:train}
\small\DontPrintSemicolon
\KwIn{the closed loop $F$ on $\bar B_\gamma$, the threshold $r_b$, the margin $\zeta$, the iteration counts $T$ and $T'$}
\KwOut{the network $\psi$ and the coefficient $\kappa$}
estimate $A = DF(0)$ by central differences and compute $\Pm$ from \eqref{eq:stein}\;
initialize $\psi$, $\kappa$ and the hinge weight $\lambda$, set $\calC, \calX \gets \emptyset$, and fill a pool $\mathcal{P}$ with points of $B_\gamma$\;
\For{$i = 0, \dots, T - 1$}{
  every 10 iterations, search $B_\gamma$ for counterexamples, adding the points with $\delta V > -\zeta$ to $\calC$ and the points with $\ell(\Exp z) \ge r_b$ and $\|F(z)\| > \gamma$ to $\calX$\;
  draw batches $Z_i$ from $\mathcal{P}$, $C_i$ from $\calC$, $X_i$ from $\calX$ and a shell sample $Y_i$, take one Adam step on the loss $\calL_i$, and update $\lambda$ by dual ascent on the hinge\;
}
\For{$j = 0, \dots, 6$}{
  compute a quick estimate of the region of the current candidate on fresh samples, with one round of the search of the level rule, and stop if no sampled point fails decrease and the search does not lower the level, or if $j = 6$\;
  pin the counterexamples found, and run the loop above again from the current $\psi$ and $\kappa$ for $T'$ iterations, with a new hinge weight, pool, buffers and optimizer state, drawing half of each $C_i$ and $X_i$ from the pinned points\;
}
\Return{the candidate with the largest quick estimate}
\end{algorithm2e}

Let $\calB = \{z \in \bar B_\gamma : \ell(\Exp z) \ge r_b,\ \delta V(z) \ge 0\} \cup \{z \in \bar B_\gamma : \|F(z)\| > \gamma\}$ be the set of points at which decrease or containment fails. The level rule places the level just below the smallest value of $V$ at the points of $\calB$ that it finds. It first evaluates two independent samples of the ball. It then searches $\Omega_\rho \setminus N_b$ for points of $\calB$ by gradient ascent on $\delta V$ and on $\|F(z)\| - \gamma$, lowers $\rho$ after each round that finds some, and stops after a round that finds none, or after 50 rounds. This search is run twice, from different starting points, and entry is checked on a separate sample of $N_b$ (Appendix~\ref{app:training}). The reported region always comes from the level rule, on samples that the refinement never used.

\vspace{-0.075in}

\paragraph{Baselines and oracle.}\looseness -1 Two standard Lyapunov functions, the quadratic baselines, go through the same level rule, with the same $r_b$. The squared-norm baseline is $V_d(z) = \|z\|^{2}$. The linearization-based baseline is $V_P(z) = z^{\top}Pz$, with $P$ the solution of $A^{\top}PA - P = -\Id$, the Lyapunov function of the linearized closed loop. A third baseline, the neural Lyapunov baseline $V_h$, is the candidate \eqref{eq:cert} trained with the residual removed and without the refinement, with the same architecture, iterations and constants. Its only decrease loss is then a hinge on sampled decrease, which with the counterexamples is the usual loss of neural Lyapunov functions~\citep{chang2019neural,richards2018lyapunov,dai2021lyapunov}. It goes through the same level rule, with the same two runs of the search. The oracle is $\Wm$ itself, from 3000 simulated steps at each sample point, with infinite cost if at the end the rotation error, the angle of the rotation part of $\tilde g$, or the translation error, the length of its translation part, is still above $0.05$. Its level comes from the two samples alone, without the search, so its region is a sampled estimate that can be optimistic (Appendix~\ref{app:results}). That region, the largest sublevel set of $\Wm$ in the ball on which no sampled point fails containment, is the reference for the fraction a function recovers. It estimates the region that the candidate approximates (Section~\ref{sec:guarantee}), and it can be smaller than the part of the ball that converges, since converging trajectories may leave the ball.

\vspace{-0.1in}

\section{Guarantees}
\label{sec:guarantee}
\vspace{-0.05in}

\looseness -1 Throughout this section, $\mathsf{G}$ is a matrix Lie group and three standing hypotheses hold: (H1) the corrected policy $u : \mathsf{G} \to \R^{m}$ is continuous, and smooth near $\Id$; (H2) $\gamma + \amax < \rG$; (H3) the goal is locally exponentially stable. For $c \ge 0$ let $K_c = \{z \in \bar B_\gamma : \Wm(z) \le c\}$, and let $\bar c = \sup\{c \ge 0 : F(K_c) \subseteq \bar B_\gamma\}$. Then $\bar c \in (0, \infty]$, and for every $c \in (0, \bar c)$, the set $K_c$ is a compact, positively invariant inner estimate that contains a neighborhood of $0$ (Lemma~\ref{lem:sublevel} in Appendix~\ref{app:proofs}).

\vspace{-0.1in}

\begin{mainthm}
\label{thm:main}
Let (H1)--(H3) hold. Let $c \in (0, \bar c)$ and $\varepsilon > 0$ with $4\varepsilon < c$ and $2\varepsilon < r_b^{2}$, and suppose that every $z$ in $\bar B_\gamma$ with $\ell(\Exp z) < r_b$ lies in $K_{c - 4\varepsilon}$. Let $V$ be a candidate \eqref{eq:cert} such that\\
\textup{(A)} $|V(z) - \Wm(z)| \le \varepsilon$ for every $z$ in $K_c$, and \textup{(B)} $V(z) > c - \varepsilon$ for every $z$ in $\bar B_\gamma \setminus K_c$.\\
Then, with $\rho_* = c - 2\varepsilon$:\\
\textup{(a)} $K_{c - 3\varepsilon} \subseteq \Omega_{\rho_*} \subseteq K_{c - \varepsilon}$;\\
\textup{(b)} for every level $\rho \in (0, c - \varepsilon]$, the set $\Omega_\rho$ lies in $K_c$, containment holds on it, decrease holds with $\delta V(z) \le -(r_b^{2} - 2\varepsilon)$ for $z$ in $\Omega_\rho \setminus N_b$, and $\Exp(\Omega_\rho) \subseteq \calA$;\\
\textup{(c)} for every $\rho \in [c - 3\varepsilon, c - \varepsilon]$, entry also holds, and $\Omega_\rho$ is a compact, positively invariant inner estimate of $\calA_\gamma$ that contains a neighborhood of $0$.
\end{mainthm}
\vspace{-0.05in}
\looseness -1 The proof is in Appendix~\ref{app:proofs}. Theorem~\ref{thm:main} is conditional: (A) and (B) are hypotheses on the candidate, and the results of this section do not show that training finds a candidate that satisfies them. In training, the residual and the quadratic term serve (A), the exit and shell terms raise $V$ near the edge of the ball, which serves (B), and the hinge and the counterexamples enforce decrease where the fit is poor. Such candidates exist: for every $c$ in $(0, \bar c)$, a small enough threshold $r_b$ and every small enough $\varepsilon > 0$, a ReLU network with one hidden layer gives a candidate that satisfies all hypotheses of Theorem~\ref{thm:main} (Corollary~\ref{cor:exists} in Appendix~\ref{app:proofs}). Hypothesis (A) concerns the distance to $\Wm$, which training does not evaluate. Training evaluates the residual of \eqref{eq:recursion}, and the next proposition bounds the distance by it.

\vspace{-0.075in}

\begin{proposition}
\label{prop:residual}
Let (H1)--(H3) hold, and let $c \in (0, \bar c)$. For every $r > 0$ there is an integer $k_r$, which depends on the policy, $c$ and $r$ but not on $V$, and which Appendix~\ref{app:proofs} gives explicitly, with the following property. If $V : \bar B_\gamma \to \R$ and $\eta \ge 0$ satisfy $|V(z) - \ell(\Exp z)^{2} - V(F(z))| \le \eta$ for every $z$ in $K_c$, then $\sup_{z \in K_c}|V(z) - \Wm(z)| \le k_r\,\eta + \sup_{z \in K_c \cap B_r}|V(z) - \Wm(z)|$.
\end{proposition}
\vspace{-0.075in}
\noindent In the bound of Proposition~\ref{prop:residual}, the error near the goal is of order $r^{2}$ for the candidate~\eqref{eq:cert}~(Appendix~\ref{app:proofs}), so it is small when $r$ is small. Training minimizes an average of the residual over samples, while Proposition~\ref{prop:residual} needs its largest value on $K_c$. With $r_b = 0.1$ the theorem needs $|V - \Wm| < r_b^{2}/2 = 0.005$ on all of $K_c$. There the accumulated cost reaches the reported levels of up to $265.5$, while the fit of the candidates trained on their closed loop is a median relative error between $0.05$ and $6.7$ percent (Appendix~\ref{app:results}). So no trained candidate is known to satisfy (A).

Only on $\SO(2)$ are the regions formally verified, by interval arithmetic in double precision on all 11 closed loops, at a level lower by a relative $1.2\times10^{-6}$ on one of them (Appendix~\ref{app:so2proof}). On the other three groups they are checked on samples and by search, and a search round that finds no counterexample does not show that none exists. The neighborhood of the goal on which $V$ decreases strictly (Lemma~\ref{lem:localdec}) may not contain $N_b$, so on these three groups the region near the goal rests on the sampled entry check and on simulation. On 93 of the 94 closed loops, none of $10^{6}$ held-out points drawn in the ball lies in the region and fails decrease or containment. This bounds the measure of the failing points but does not exclude them (Appendix~\ref{app:ablations}).

\vspace{-0.1in}

\section{Experiments}
\label{sec:sims}
\looseness -1 Three seeds of each of PPO~\citep{schulman2017ppo}, TRPO~\citep{schulman2015trpo}, SAC~\citep{haarnoja2018sac} and TD3~\citep{fujimoto2018td3} are trained on $\SO(2)$, $\SO(3)$, $\SE(2)$ and $\SE(3)$, with $\gamma = \pi - \pi/8$ and the clip length $\amax = 0.314$ on the rotation groups and $0.351$ on the rigid-motion groups (Appendix~\ref{app:policies}). Thirty-nine of the 48 policies pass the stability test, and the nine that fail include all six SAC and TD3 policies on $\SE(3)$ (Table~\ref{tab:fail} in Appendix~\ref{app:results}). These six settle near the goal without reaching it (Table~\ref{tab:policies} in Appendix~\ref{app:policies}), and after the correction the goal is an unstable fixed point, so the method gives no region for them. Together with the 55 closed loops with an actuator error described below, 94 closed loops are evaluated. The size of a region is estimated by the fraction of $2\times10^{5}$ initial errors drawn from $\mu_0$ that lie in it, and a trajectory converges if its rotation and translation errors are below $0.05$ after 600 simulated steps. The ball $\bar B_\gamma$ holds $87.6$, $75.4$, $82.0$ and $65.1$ percent of the samples of $\mu_0$ on the four groups, and every sampled initial error in the ball converges on every policy. The ball, not the region of attraction, therefore limits the regions, except on TRPO seed 2 on $\SO(2)$, which has a second fixed point inside the ball (Appendix~\ref{app:results}).

\vspace{-0.075in}

\paragraph{Policies without an actuator error.} Most of the 39 policies do not separate the methods. On 31 of them the candidate, both baselines and the oracle give the whole ball, to within $0.01$ percentage points (Tables~\ref{tab:main} and \ref{tab:all} in Appendix~\ref{app:results}). On the other 8 the candidate gives a larger region than both baselines, by more than one percentage point on 6. Five of these are on $\SE(2)$ and $\SE(3)$, where the candidate's region exceeds that of $V_d$ by $6.0$ to $24.3$ percentage points and that of $V_P$ by $6.4$ to $24.3$. The sixth is TD3 seed 1 on $\SO(3)$, where both baselines and the quadratic term of the candidate alone give below $0.01$ percent, and the entry condition fails for $V_P$ near the goal (Appendix~\ref{app:results}), while the candidate gives the whole ball. On no policy is the candidate smaller than a baseline by more than two of the $2\times10^{5}$ sample points. The candidate recovers at least 99 percent of the oracle region on 38 of the 39 policies, and 81 percent on PPO seed 1 on $\SE(3)$ (Appendix~\ref{app:results}).

\begin{figure}[!t]
\centering
\includegraphics[width=0.75\linewidth]{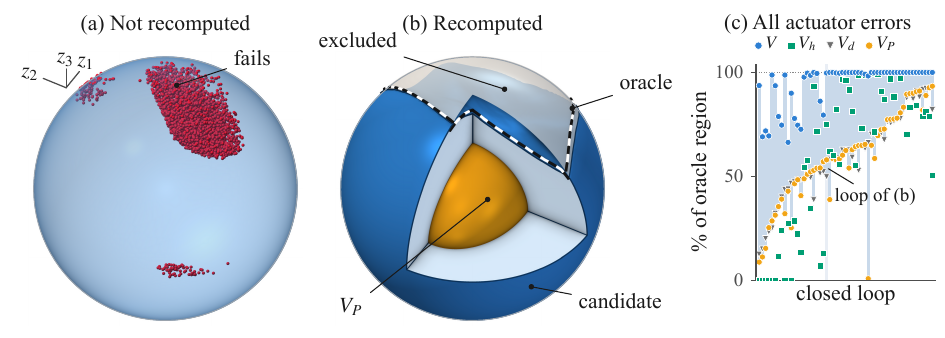}

\vspace{-0.05in}

\caption{PPO seed 0 on $\SE(2)$ with the actuators misaligned by $40^\circ$, in exponential coordinates inside the ball $\bar B_\gamma$. (a) The region computed without the actuator error. Under the misalignment, decrease or containment fails at the red points. (b) The region recomputed with the actuator error (blue, one octant cut away). It leaves out the gray part of the ball, and the dashed line is the boundary of the oracle region. The region of $V_P$, the larger quadratic baseline region here, is orange. (c) The regions of $V$, $V_h$, $V_d$ and $V_P$ as percentages of the oracle region on the 55 closed loops with an actuator error, ordered by the larger of the two quadratic baseline regions. The shaded column is the closed loop of (b).}
\label{fig:actuator}
\vspace{-0.1in}
\end{figure}

\vspace{-0.075in}

\paragraph{Actuator errors.} A gain error multiplies the increment by a constant, here $0.5$, $0.75$, $1.25$ or $1.5$, or each of its coordinates by its own gain drawn from $[0.7, 1.3]$, the per-axis gains. A misalignment rotates the increment by $5^\circ$, $10^\circ$, $15^\circ$, $30^\circ$ or $40^\circ$. Both give a matrix $E$ in \eqref{eq:closedloopE}. Robust analysis certifies one region for a set of input uncertainties in $\R^{n}$~\citep{yin2022stability}. Here the error is identified, and the region is recomputed. To find where the baselines fall short, we screened 300 closed loops with these errors, all passing the stability test, by comparing the convergent part of the ball with the larger quadratic baseline region (Appendix~\ref{app:screen}). The two differ by less than $1.7$ percentage points on $\SO(2)$ and $\SO(3)$, and by at least 5 points on 55 of the 198 closed loops on $\SE(2)$ and $\SE(3)$. Only these 55 are evaluated~(Table~\ref{tab:pert}, and Table~\ref{tab:pert-all} in Appendix~\ref{app:results}). For 29 of them, the candidate starts from that of the same policy without the actuator error. It is refitted on 15 and kept with a new level on 14, where its region near the goal rests on the sampled entry check (Appendix~\ref{app:training}).

\vspace{-0.075in}

\paragraph{Results with an actuator error.} An actuator error can break the conditions on the region computed without it. On 21 of the 32 closed loops with a misalignment, including all 16 at $30^\circ$ and $40^\circ$, that region fails decrease or containment at sampled points (Figure~\ref{fig:actuator}(a)), and the search finds more. Every sampled initial error in it still converges (Appendix~\ref{app:ablations}). Without the conditions, however, convergence is known only at the simulated points, so the region has to be recomputed. The recomputed region is larger than both quadratic baseline regions on all 55 closed loops (Figure~\ref{fig:actuator}(b, c)), by a median of $24.0$ percentage points over $V_d$ and $24.4$ over $V_P$, and by at least $6.2$ and $5.4$ points. It is also larger than the neural Lyapunov baseline on 49 of the 55, equal to it within $0.01$ on 5, and $0.8$ points smaller on one, TD3 seed 1 on $\SE(2)$ at $40^\circ$, with a median difference of $18.1$ points over all 55. The neural Lyapunov baseline exceeds the larger quadratic baseline on 29 closed loops and falls below it on 26; on $\SE(3)$, it yields no region because the entry condition fails at the levels the search leaves (Appendix~\ref{app:results}). The recomputed region recovers a median of $99.9$ percent of the oracle region, against $73.0$ for $V_h$, $59.1$ for $V_d$ and $58.5$ for $V_P$. It recovers less than 99 percent on 19 closed loops, 17 of them PPO on $\SE(3)$.

\begin{table}[!t]
\centering
\caption{The 55 closed loops with an actuator error selected by the screens of Appendix~\ref{app:screen}, by group and actuator error. Entries other than the number of closed loops are in percent.}
\label{tab:pert}
\scriptsize\setlength{\tabcolsep}{1.3pt}
\begin{tabular}{@{}llrrrrrrrrrr@{}}
\toprule
group & actuator error & loops & oracle & $V$ & $V_d$ & $V_P$ & $V_h$ & $V$/oracle & $V_d$/oracle & $V_P$/oracle & $V_h$/oracle \\
\midrule
$\SE(2)$ & gain 0.5--1.5 & 6 & 82.0 & 82.0 & 57.5--75.7 & 53.2--72.9 & 64.7--81.5 & 100.0 & 70.2--92.3 & 64.9--88.9 & 78.9--99.4 \\
& per-axis gains & 3 & 82.0 & 81.9--82.0 & 55.3--75.8 & 59.4--75.9 & 64.7--74.5 & 99.9--100.0 & 67.5--92.5 & 72.4--92.6 & 78.9--90.8 \\
& misalignment 5--15$^\circ$ & 8 & 82.0 & 82.0 & 62.7--74.6 & 63.6--74.6 & 57.6--82.0 & 100.0 & 76.5--91.0 & 77.6--91.1 & 70.2--100.0 \\
& misalignment 30--40$^\circ$ & 7 & 81.5--82.0 & 80.6--82.0 & 31.7--67.2 & 0.6--76.5 & 41.4--82.0 & 98.3--100.0 & 38.6--82.0 & 0.7--93.4 & 50.6--100.0 \\
\midrule
$\SE(3)$ & gain 0.5--1.5 & 8 & 65.1 & 47.3--65.1 & 30.3--41.4 & 16.4--42.6 & 0.0--64.1 & 72.6--100.0 & 46.6--63.5 & 25.2--65.4 & 0.0--98.5 \\
& per-axis gains & 6 & 65.1 & 48.6--65.1 & 22.5--42.3 & 23.1--42.3 & 7.5--62.7 & 74.6--99.9 & 34.5--64.9 & 35.4--64.9 & 11.5--96.2 \\
& misalignment 5--15$^\circ$ & 8 & 65.1 & 43.2--65.1 & 26.0--38.2 & 25.4--38.1 & 15.5--46.7 & 66.3--100.0 & 40.0--58.7 & 39.0--58.5 & 23.9--71.8 \\
& misalignment 30--40$^\circ$ & 9 & 63.5--65.1 & 44.8--65.1 & 7.9--45.9 & 5.5--38.2 & 0.0--46.5 & 69.0--100.0 & 12.4--70.5 & 8.7--58.6 & 0.0--71.4 \\
\bottomrule
\end{tabular}
\vspace{-4pt}
\end{table}

\vspace{-0.075in}

\paragraph{Cost.} On two CPU threads, training and refinement take $2.5$ to $6.7$ hours per trained or refitted closed loop with an actuator error, against 18 to 48 minutes for the oracle (Appendices~\ref{app:training} and \ref{app:results}). The candidate's conditions involve a single step, so the cost of checking them does not depend on how long the trajectories take to converge. Outside $\SO(2)$, the result is still a sampled estimate.

\vspace{-0.075in}

\paragraph{Ablations.} On five closed loops, one part of the method is changed at a time (Table~\ref{tab:ablation} in Appendix~\ref{app:ablations}). The hinge-only variant, without the residual, and the variant that replaces the residual with a regression on $\Wm$ each lose more than 50 percentage points on at least one of them. The variant without refinement loses up to $15.0$ points. The variant without the quadratic term loses up to $9.3$, and gains $1.1$ on one closed loop. The region away from the goal comes from the network: on 14 of the 16 closed loops formed by the 8 policies on which the candidate exceeds both baselines and 8 closed loops with an actuator error, the quadratic term $z^{\top}\Pm z$ alone gives $27.9$ to $75.4$ points less than the candidate. Appendix~\ref{app:ablations} reports the sensitivity to $r_a$ and $r_b$.

\vspace{-0.075in}

\paragraph{Deployment on two vehicles.}\looseness -1 The regions assume that each increment is carried out exactly, so two corrected policies are run on simulated vehicles. In MuJoCo \citep{todorov2012mujoco}, PPO seed 0 on $\SO(3)$ steers the attitude of a 3U CubeSat and PPO seed 0 on $\SE(3)$ steers the pose of a BlueROV2 underwater vehicle, through proportional-derivative controllers. In 200 episodes of 40 policy steps on each vehicle, from initial errors drawn uniformly in $\Omega_\rho$, none leaves $\Omega_\rho$, on the closed loop or on the vehicle (Table~\ref{tab:deploy}, Figure~\ref{fig:deploy}, Appendix~\ref{app:deploy}). Nor does any leave it with an initial body rate or a halved sampling period. On all 28 closed loops of PPO on $\SE(3)$ with an actuator error, where the BlueROV2 tracks the distorted increments, none of the 200 episodes leaves the recomputed region, while one leaves the region computed without the error at $40^\circ$ (Appendix~\ref{app:deploy}). Without disturbances or model mismatch, the runs test how closely the controller tracks the increments.

\vspace{-0.15in}

\section{Conclusion}
\vspace{-0.05in}

In this work, we train RL policies on matrix Lie groups and estimate their ROA. We recompute the ROA estimate after an identified actuator gain error or misalignment, without retraining the policy. The regions are checked on samples and by adversarial search, except on $\SO(2)$, where interval arithmetic verifies them. On the 55 closed loops with an actuator error on $\SE(2)$ and $\SE(3)$ where the baselines fall short, the recomputed regions exceed those of both quadratic baselines on all and that of the neural Lyapunov baseline on 49. The region computed without the error fails the sampled conditions on 21 of the 32 closed loops with a misalignment. However, every sampled initial error in the ball still converges on all 55, so the region needs to be recomputed; there is no need to retrain the policy. Future work includes formal verification beyond $\SO(2)$, and higher-dimensional systems.

\noindent \textbf{Generative AI Use:} Generative AI tools were used improve the syntax and grammar of some paragraphs in the manuscript. The authors reviewed and verified all AI-assisted content.

\bibliography{ref}

\appendix

\section{Matrix Lie groups and the chart assumption}
\label{app:lie}
Throughout, $\mathsf{G}$ is a matrix Lie group with Lie algebra $\g$ of dimension $m$, a basis $\hat e_1, \dots, \hat e_m$ of $\g$ identifies $\g$ with $\R^{m}$, and $\|\cdot\|$ is the Euclidean norm of the coordinates. The curve $t \mapsto \Exp(t\xi)$ is the one-parameter subgroup with generator $\hat\xi$ \citep[Definition 2.13 and Theorem 2.14]{hall2015lie}, so $\Exp(t_1\xi)\Exp(t_2\xi) = \Exp((t_1 + t_2)\xi)$. Left translation carries the inner product of $\g$ to every tangent space, which gives a left-invariant Riemannian metric on $\mathsf{G}$ with distance $\dist$ \citep[Theorem 5.38]{bullo2005geometric}.

\paragraph{The step bound.} The curve $t \mapsto g\Exp(t\tau)$, $0 \le t \le 1$, joins $g$ to $g\Exp(\tau)$. Its velocity is $g\Exp(t\tau)\hat\tau$ \citep[Proposition 5.5]{bullo2005geometric}, the left translate of $\hat\tau$, so in the left-invariant metric its speed is the constant $\|\tau\|$, and $\dist(g, g\Exp\tau) \le \|\tau\|$.

\paragraph{The chart assumption.} The results of the paper use three properties of the exponential coordinates, collected in the following assumption. Proposition~\ref{prop:general} shows that every matrix Lie group of positive dimension satisfies it for some chart radius, and Proposition~\ref{prop:assumption} gives the chart radius $\pi$ for the rotation and rigid-motion groups of the experiments.

\begin{assumption}
\label{ass:chart}
There are a continuous function $\rot : \mathsf{G} \to [0, \infty)$ and a number $\rG > 0$, the chart radius, with the following three properties. (i) The inequalities $\rot(gh) \le \rot(g) + \rot(h)$ and $\rot(\Exp\xi) \le \|\xi\|$ hold for all $g$, $h$ in $\mathsf{G}$ and $\xi$ in $\R^{m}$. (ii) The exponential map is injective on an open set $\Dch \subseteq \R^{m}$ with $B_{\rG} \subseteq \Dch$, and its image is the open set $\Uch = \{g \in \mathsf{G} : \rot(g) < \rG\}$. (iii) The inverse $\Log : \Uch \to \Dch$ is smooth.
\end{assumption}

The function $\rot$ measures how far an element is from $\Id$, and $\Dch$ may be larger than the ball.

\begin{proposition}
\label{prop:general}
Every matrix Lie group of positive dimension satisfies Assumption~\ref{ass:chart} for some chart radius $\rG > 0$, with $\rot(g) = \min\{\dist(\Id, g), \rG\}$ and $\dist$ the distance of the left-invariant metric, which is infinite between elements that no curve in $\mathsf{G}$ joins.
\end{proposition}

\begin{proof}
The group $\mathsf{G}$ need not be connected. Its identity component $\mathsf{G}_0$, the set of elements joined to $\Id$ by a continuous path in $\mathsf{G}$, is an open subgroup, because $\mathsf{G}$ is locally path connected \citep[remark after Corollary 3.45]{hall2015lie}. It is itself a matrix Lie group with the same Lie algebra \citep[Corollary 3.52]{hall2015lie}. The distance $\dist(\Id, g)$ is finite exactly on $\mathsf{G}_0$ and continuous there in the manifold topology. Off $\mathsf{G}_0$ it is infinite, so with values in $[0, \infty]$ it is continuous on $\mathsf{G}$, since $\mathsf{G}_0$ is open and closed.

Part (i). The distance satisfies the triangle inequality, and left translation by $g$ preserves it, so $\dist(\Id, gh) \le \dist(\Id, g) + \dist(g, gh) = \dist(\Id, g) + \dist(\Id, h)$, in $[0, \infty]$. Cutting off at any $r > 0$ keeps the inequality: $\min\{a + b, r\} \le \min\{a, r\} + \min\{b, r\}$ for $a, b$ in $[0, \infty]$, because the right side is at least $r$ when $a \ge r$ or $b \ge r$ and equals $a + b$ otherwise. The step bound with $g = \Id$ gives $\dist(\Id, \Exp\xi) \le \|\xi\|$, hence $\rot(\Exp\xi) \le \|\xi\|$. The function $\rot$ is continuous as the minimum of two continuous functions.

Parts (ii) and (iii). Let $U_{\varepsilon_0}$ be the ball of radius $\varepsilon_0$ about $0$ in the Hilbert--Schmidt norm of the matrices, which is the Frobenius norm. By \citet[Theorem 3.42]{hall2015lie} there is $\varepsilon_0$ in $(0, \log 2)$ such that every $M$ in $\exp(U_{\varepsilon_0}) \cap \mathsf{G}$ satisfies $\log M \in \g$, and then $\exp$ maps $U_{\varepsilon_0} \cap \g$ onto $\exp(U_{\varepsilon_0}) \cap \mathsf{G}$ with the matrix logarithm as inverse \citep[Corollary 3.44]{hall2015lie}. The matrix logarithm is given by a convergent power series on $\{\|M - \Id\| < 1\}$ \citep[Theorem 2.8]{hall2015lie}, so it is smooth there, and $\varepsilon_0 < \log 2$ places $\exp(U_{\varepsilon_0})$ inside this set. Two norms on the finite-dimensional space $\g$ bound each other, so there is $r_* > 0$ with $B_{r_*} \subseteq U_{\varepsilon_0} \cap \g$ in coordinates. So $\Exp$ is injective on $B_{r_*}$, its image $\Exp(B_{r_*})$ is open in $\mathsf{G}$, and $\Log = \log$ is smooth on this image, as the restriction of a real-analytic map on the matrices to the embedded submanifold $\mathsf{G}$. The image $\Exp(B_{r_*})$ is an open neighborhood of $\Id$ in $\mathsf{G}$, hence in $\mathsf{G}_0$, and the metric balls $\{\dist(\Id, \cdot) < \rho_{\mathsf{G}}\}$, which lie in $\mathsf{G}_0$, form a basis of neighborhoods of $\Id$ in $\mathsf{G}_0$, since the topology of a Riemannian distance is the manifold topology. Hence $\{\dist(\Id, \cdot) < \rho_{\mathsf{G}}\} \subseteq \Exp(B_{r_*})$ for some $\rho_{\mathsf{G}} > 0$. Put $\rG = \min\{\rho_{\mathsf{G}}, r_*\}$, $\rot = \min\{\dist(\Id, \cdot), \rG\}$, $\Uch = \{\rot < \rG\} = \{\dist(\Id, \cdot) < \rG\}$ and $\Dch = \{\xi \in B_{r_*} : \Exp\xi \in \Uch\}$. Then $\Uch$ is open because $\rot$ is continuous, $\Dch$ is open because $\Exp$ is continuous, $\Exp$ is injective on $\Dch \subseteq B_{r_*}$ with image $\Uch$, and $\Log : \Uch \to \Dch$ is smooth. Last, $B_{\rG} \subseteq \Dch$, since for $\|\xi\| < \rG$ part (i) gives $\rot(\Exp\xi) \le \|\xi\| < \rG$, so $\Exp\xi$ is in $\Uch$ and $\xi$ is in $B_{r_*}$. On a connected group the cutoff changes nothing near $\Id$, and on $\SO(n)$ and $\SE(n)$ Proposition~\ref{prop:assumption} uses the distance on the rotation group, defined below, without a cutoff.
\end{proof}

The proof gives a chart radius but not its value. The value is computed below for the rotation and rigid-motion groups, by the same route: a function $\rot$ with the triangle inequality, an explicit set on which $\Exp$ is injective, and a smooth inverse. Every compact group has an inner product on $\g$ that is invariant under conjugation \citep[Proposition 21.6]{gallier2020differential}. With a basis orthonormal for it, the left-invariant metric of this appendix is bi-invariant \citep[Proposition 21.3]{gallier2020differential}, and Proposition~\ref{prop:general} supplies $\rot$ and $\rG$ for that metric, as Proposition~\ref{prop:lefschetz} needs.

\paragraph{The rotation and rigid-motion groups.} The special orthogonal group is $\SO(n) = \{R \in \R^{n\times n} : R^{\top}R = \Id,\ \det R = 1\}$, with Lie algebra $\so(n)$ of skew-symmetric matrices. The special Euclidean group $\SE(n)$ is the set of pairs $(R, p)$ with $R$ in $\SO(n)$ and $p$ in $\R^{n}$, realized as the block matrices $\bigl(\begin{smallmatrix} R & p \\ 0 & 1 \end{smallmatrix}\bigr)$, with Lie algebra $\se(n) = \{\bigl(\begin{smallmatrix} X & v \\ 0 & 0 \end{smallmatrix}\bigr) : X \in \so(n),\ v \in \R^{n}\}$. Both are closed subgroups of a general linear group, hence matrix Lie groups \citep[Definition 1.4]{hall2015lie}. Both are connected, since $\SO(n)$ is the image of $\so(n)$ under $\exp$ \citep[Theorem 2.6]{gallier2020differential} and $\SE(n)$ is diffeomorphic to $\SO(n)\times\R^{n}$. The first is compact and the second is not. Their dimensions are $m = n(n-1)/2$ and $m = n(n+1)/2$. The inner product $\tfrac12\tr(X^{\top}Y)$ on $\so(n)$ does not change under conjugation by $\SO(n)$, so it extends to a bi-invariant metric on $\SO(n)$ \citep[Proposition 21.3]{gallier2020differential}. Let $E_{ij}$ be the matrix with a single unit entry in row $i$ and column $j$. For a basis of $\so(n)$ formed by the matrices $\pm(E_{ij} - E_{ji})$ with $i < j$, in some order, this inner product is the Euclidean one on $\R^{m}$ for every $n$. The paper uses such a basis. On $\SO(3)$ it is the basis of the cross-product matrices of the unit vectors, so the coordinates $\xi$ of a body velocity $\hat\xi$ are its angular velocity vector. On $\SE(n)$ the basis consists of the matrices $\bigl(\begin{smallmatrix} X & 0 \\ 0 & 0 \end{smallmatrix}\bigr)$ with $X$ in this basis of $\so(n)$, followed by the matrices $\bigl(\begin{smallmatrix} 0 & e_i \\ 0 & 0 \end{smallmatrix}\bigr)$ with $e_i$ the unit vectors of $\R^{n}$, so the coordinates are $z = (\omega, v)$. The upper right block of the matrix exponential of $\bigl(\begin{smallmatrix} \hat\omega & v \\ 0 & 0 \end{smallmatrix}\bigr)$ is $\sum_{k \ge 1}\hat\omega^{k-1}v/k! = \int_0^{1}\Exp(t\omega)\,dt\ v = J(\omega)v$, which gives the exponential map of $\SE(n)$. For $n = 3$ it is the closed form of \citet[eq.~(2.36)]{murray1994mathematical}.

\paragraph{The distance on the rotation group.} Every rotation $R$ is orthogonally similar to a direct sum of planar rotations by angles $\theta_i(R)$ in $(-\pi, \pi]$ and an identity block. Since $\SO(n)$ is compact, it is complete, and every pair of points is joined by a minimal geodesic \citep[Theorems 16.16 and 16.17 and the remark after the latter]{gallier2020differential}. For a bi-invariant metric, the geodesics through $\Id$ are the one-parameter subgroups \citep[Proposition 21.21]{gallier2020differential}. A minimal geodesic from $\Id$ to $R$ is therefore a curve $t \mapsto \Exp(t\xi)$ with $\Exp\xi = R$ and length $\|\xi\|$, and every other logarithm of $R$ gives a curve at least as long. Hence $\dist(\Id, R) = \min\{\|\xi\| : \Exp\xi = R\}$. A logarithm $\hat\xi$ of $R$ is orthogonally similar to a direct sum of planar generators with angles $\varphi_j$ and a zero block, with $\|\xi\|^{2} = \sum_j\varphi_j^{2}$. Its exponential has the eigenvalues $e^{\pm i\varphi_j}$ and $1$, which coincide with the eigenvalues $e^{\pm i\theta_i(R)}$ and $1$ of $R$, with multiplicity. After padding both lists with zero angles and reordering, $\varphi_j \equiv \pm\theta_j(R)$ modulo $2\pi$, and $\theta_j(R)$ is the representative of least modulus, so $|\varphi_j| \ge |\theta_j(R)|$. Hence $\|\xi\|^{2} \ge \sum_j\theta_j(R)^{2}$, with equality for the logarithm built from the angles $\theta_j(R)$ and the same similarity, and $\dist(\Id, R)^{2} = \sum_i\theta_i(R)^{2}$. We call this distance the rotation angle $\rot(R)$. It is the function $\rot$ of Assumption~\ref{ass:chart} for $\SO(n)$.

\begin{proposition}
\label{prop:assumption}
Let $\mathsf{G}$ be $\SO(n)$ or $\SE(n)$, let $\rot(g) = \dist_{\SO(n)}(\Id, R)$ for the rotation part $R$ of $g$, and let $\Dch = B_\pi$ on $\SO(n)$ and $\Dch = \{(\omega, v) : \|\omega\| < \pi\}$ on $\SE(n)$. Then $\mathsf{G}$ satisfies Assumption~\ref{ass:chart} with chart radius $\rG = \pi$.
\end{proposition}

\begin{proof}
The rotation part of a product is the product of the rotation parts, and $\dist(\Id, RS) \le \dist(\Id, R) + \dist(R, RS) = \dist(\Id, R) + \dist(\Id, S)$ by left invariance, which is part (i) for products. By the step bound, $\rot(\Exp\omega) \le \|\omega\|$, and on $\SE(n)$ the rotation part of $\Exp(\omega, v)$ is $\Exp\omega$, so $\rot(\Exp(\omega, v)) \le \|\omega\| \le \|(\omega, v)\|$. The ball $B_\pi$ lies in $\Dch$ on both groups.

On $\SO(n)$, a skew-symmetric matrix $\hat\xi$ is orthogonally similar to a direct sum of planar generators with angles $\theta_i$ and a zero block, and $\|\xi\|^{2} = \sum_i\theta_i^{2}$, so $\|\xi\| < \pi$ forces $|\theta_i| < \pi$ for every $i$. Its eigenvalues $\pm i\theta_i$, and possibly $0$, lie in the strip $\{\lambda \in \mathbb{C} : |\operatorname{Im}\lambda| < \pi\}$, on which the matrix exponential is injective \citep[Theorem 10.15]{gallier2020differential}, so $\Exp$ is injective on $B_\pi$. The image of $B_\pi$ lies in $\{\rot < \pi\}$ by the bound above. Conversely, a rotation with $\rot(R) < \pi$ is the exponential of the logarithm built from its angles, whose norm is $\rot(R)$. The image is therefore $\{\rot < \pi\}$, which is open because $\rot$ is continuous. The eigenvalues of $\operatorname{ad}_{\hat\xi}$ are $\pm i(\theta_j \pm \theta_k)$ and $\pm i\theta_j$, of modulus at most $\sqrt2\,\|\xi\| < 2\pi$ on $B_\pi$. So the differential of $\Exp$ is invertible there \citep[Proposition 20.2]{gallier2020differential}, and an injective local diffeomorphism is a diffeomorphism onto its open image.

On $\SE(n)$, the equality $\Exp(\omega, v) = \Exp(\omega', v')$ forces $\Exp\omega = \Exp\omega'$, hence $\omega = \omega'$, and then $J(\omega)v = J(\omega)v'$. The eigenvalues of $J(\omega)$ are $(e^{i\theta} - 1)/(i\theta)$ and $1$, nonzero for $|\theta| < 2\pi$, so $J(\omega)$ is invertible and $v = v'$. The image is the set of pairs whose rotation part has $\rot < \pi$, and the inverse $(R, p) \mapsto (\Log R, J(\Log R)^{-1}p)$ is smooth.
\end{proof}

For $n \le 3$ a rotation has at most one nonzero angle, and $\{\rot < \pi\}$ is the set of rotations without eigenvalue $-1$. On $\SO(3)$ it is $\{\tr R \ne -1\}$, with the closed-form logarithm of \citet[Proposition 5.7]{bullo2005geometric}, and it is the regular set of \citet{akhtar2026learning}. For $n \ge 4$ the two sets differ: angles $(0.9\pi, 0.9\pi)$ give no eigenvalue $-1$ and $\rot = 0.9\sqrt2\,\pi > \pi$. For $\SO(2)$ and $\SO(3)$ the boundary sphere of $B_\pi$ is the tangential cut locus of $\Id$ \citep[Definition 16.14 and the sentence after it]{gallier2020differential}.

\paragraph{The running cost, the cost bound and the matrix $Q$.} On every matrix Lie group, $\Exp z - \Id = \sum_{k \ge 1}\hat z^{k}/k!$, so $\ell(\Exp z) \le \sum_{k \ge 1}\|\hat z\|_F^{k}/k! = e^{\|\hat z\|_F} - 1 \le \|\hat z\|_F\,e^{\|\hat z\|_F}$, and $\|\hat z\|_F^{2} = z^{\top}Qz \le \lambda_{\max}(Q)\|z\|^{2}$ with $Q_{ij} = \tr(\hat e_i^{\top}\hat e_j)$ as in Section~\ref{sec:local}. Hence $\ell(\Exp z)^{2} \le \Cl\,\|z\|^{2}$ for every $z$ in $\bar B_{\rG}$, the cost bound, with $\Cl = \lambda_{\max}(Q)\exp(2\lambda_{\max}(Q)^{1/2}\rG)$. The same expansion gives $\ell(\Exp z)^{2} = \|\hat z\|_F^{2} + O(\|z\|^{3}) = z^{\top}Qz + O(\|z\|^{3})$, the expansion of Section~\ref{sec:local}, and $Q$ is positive definite because the $\hat e_i$ are linearly independent.

On $\SO(2)$ and $\SO(3)$ the running cost $\ell^{2}$ equals $2\tr(\Id - \tilde g)$, the trace error function of geometric attitude control~\citep{bullo1999tracking}, and depends on the rotation angle alone. A region defined on the group through it would again be a ball in exponential coordinates, and the ball $\bar B_\gamma$ excludes only the rotations by more than $\gamma$.

On $\SO(n)$ and $\SE(n)$ the cost bound holds with $\Cl = 2$ for every $z$. A planar rotation by $\theta$ differs from the identity by $8\sin^{2}(\theta/2)$ in squared Frobenius norm. So on $\SO(n)$, $\ell(R)^{2} = \|R - \Id\|_F^{2} = 2\tr(\Id - R) = \sum_i 8\sin^{2}(\theta_i(R)/2)$. For $R = \Exp z$ the angles of $R$ are the angles $\varphi_i$ of $\hat z$ reduced modulo $2\pi$, and $8\sin^{2}(\theta/2)$ has period $2\pi$, so $\ell(\Exp z)^{2} = \sum_i 8\sin^{2}(\varphi_i/2)$ with $\sum_i\varphi_i^{2} = \|z\|^{2}$. The bound $8\sin^{2}(\theta/2) \le 2\theta^{2}$ gives $\ell(\Exp z)^{2} \le 2\|z\|^{2}$, and the expansion $8\sin^{2}(\theta/2) = 2\theta^{2} - \theta^{4}/6 + \dots$ gives $\ell(\Exp z)^{2} = 2\|z\|^{2} + O(\|z\|^{4})$, so $Q = 2\Id$, in agreement with $\|\hat z\|_F^{2} = 2\|z\|^{2}$ for the inner product $\tfrac12\tr(X^{\top}Y)$. On $\SE(n)$ the translation of $\Exp(\omega, v)$ is $J(\omega)v$, and $J(\omega)$ is an average of orthogonal matrices, so $\|J(\omega)v\| \le \|v\|$ and the same bound holds. Since $J(\omega) = \Id + \tfrac12\hat\omega + O(\|\omega\|^{2})$ and $v^{\top}\hat\omega v = 0$, $\|J(\omega)v\|^{2} = \|v\|^{2} + O(\|z\|^{4})$, so $Q = \operatorname{diag}(2\Id, \Id)$. On $\SO(3)$ the running cost is twice the attitude error function $\tfrac12\tr(G_P(\Id - R))$ of \citet[eq.~(11.18)]{bullo2005geometric} with gain matrix $\Id$, for which $G_P = \tr(\Id)\Id - \Id = 2\Id$.

\section{Proofs}
\label{app:proofs}
As $z \to 0$, $O(\|z\|^{k})$ denotes a function bounded by a constant times $\|z\|^{k}$ near $0$, and $o(\|z\|^{k})$ one whose quotient by $\|z\|^{k}$ tends to $0$.

\subsection{Proof of Proposition~\ref{thm:descent}}
\begin{proof}
Let $\tilde g_k = \bar g^{-1}g_k$. Multiplying the step on the left by $\bar g^{-1}$, as in \eqref{eq:errdyn}, gives $\tilde g_{k+1} = \Phi_u(\tilde g_k)$. So $\tilde g_k$ is the trajectory of \eqref{eq:closedloop} from $\tilde g_0 = \bar g^{-1}g_0$, and $g_k = \bar g\tilde g_k$. Equivalently, with the left translation $L_{\bar g}(g) = \bar g g$, the closed-loop map on the states, $\Phi_{u,\bar g}(g) = g\,\Exp(\proj(u(\bar g^{-1}g)))$, satisfies $\Phi_{u,\bar g} = L_{\bar g}\circ\Phi_u\circ L_{\bar g}^{-1}$. Left translation by $\bar g$ is a homeomorphism of $\mathsf{G}$ that maps $\Id$ to $\bar g$, so $g_k \to \bar g$ if and only if $\tilde g_k \to \Id$.
\end{proof}

\paragraph{Scope of Proposition~\ref{thm:descent}.} The proposition uses only that the step is right multiplication by an element that depends on the error $\tilde g$. It fails for the left error $g\bar g^{-1}$, whose step multiplies it on the right by $\bar g\Exp(\tau)\bar g^{-1}$, which depends on the goal, for the additive error $g - \bar g$, which is not an element of $\mathsf{G}$, and for dynamics that are not left-invariant, such as a gravity torque or a disturbance in the spatial frame. The proposition still holds under an actuator error that changes only the increment applied in the body frame, such as the gain errors and misalignments of Section~\ref{sec:sims}. Such an error multiplies the increment by the matrix $E$ of \eqref{eq:closedloopE} before the second projection onto $\bar B_{\amax}$, and both projections are inactive near the goal, so the closed loop with the error keeps $\Id$ fixed, is smooth near it and takes increments in $\bar B_{\amax}$. Hypotheses (H1) and (H2) of Section~\ref{sec:guarantee} then hold for it, and (H3) holds whenever its Jacobian at $0$ has spectral radius below $1$ (Lemma~\ref{lem:spr}), which the stability test checks numerically.

\subsection{One step stays in the chart}
\begin{figure}[!t]
\centering
\begin{tikzcd}[row sep=large, column sep=huge, ampersand replacement=\&]
\bar g\,\Exp(\bar B_\gamma) \arrow[r, "\Phi_{u,\bar g}"] \arrow[d, "L_{\bar g}^{-1}"'] \& \bar g\,\Uch \arrow[d, "L_{\bar g}^{-1}"] \\
\Exp(\bar B_\gamma) \arrow[r, "\Phi_u"] \arrow[d, "\Log"'] \& \Uch \arrow[d, "\Log"] \\
\bar B_\gamma \arrow[r, "F"] \& \Dch
\end{tikzcd}

\vspace{8pt}
\begin{tikzpicture}[>=stealth, block/.style={draw, rectangle, minimum height=7mm, inner sep=4pt}, font=\small]
\node (in) at (0,0) {$\tilde g_k$};
\node[block] (vec) at (1.45,0) {$\vecop$};
\node[block] (mu) at (3.9,0) {$\mu(\cdot) - \mu(\vecop\Id)$};
\node[block] (sat) at (6.6,0) {$\amax\sat$};
\node[block] (proj) at (9.0,0) {$\proj$};
\node[block] (exp) at (10.9,0) {$\Exp$};
\node[draw, circle, inner sep=1.5pt] (mul) at (12.4,0) {$\times$};
\node (out) at (13.9,0) {$\tilde g_{k+1}$};
\draw[->] (in) -- (vec);
\draw[->] (vec) -- (mu);
\draw[->] (mu) -- (sat);
\draw[->] (sat) -- node[above] {$u(\tilde g_k)$} (proj);
\draw[->] (proj) -- node[above] {$\tau_k$} (exp);
\draw[->] (exp) -- (mul);
\draw[->] (mul) -- (out);
\draw[->] (0.62,0) |- (12.4,-0.95) -- node[right, pos=0.4] {$\tilde g_k$} (mul);
\fill (0.62,0) circle (1.2pt);
\end{tikzpicture}
\vspace{2pt}
\caption{One step of the closed loop. Top: for a goal $\bar g$, with the left translation $L_{\bar g}(g) = \bar g g$ and the closed-loop map $\Phi_{u,\bar g}(g) = g\,\Exp(\proj(u(\bar g^{-1}g)))$ on the states, the diagram commutes, that is, $\Phi_{u,\bar g} = L_{\bar g}\circ\Phi_u\circ L_{\bar g}^{-1}$ on $\bar g\,\Exp(\bar B_\gamma)$ and $F = \Log\circ\Phi_u\circ\Exp$ on $\bar B_\gamma$. All vertical maps are bijections. Bottom: block diagram of the corrected policy of Section~\ref{sec:loop}. The error $\tilde g_k$ gives the increment $\tau_k = \proj(u(\tilde g_k))$, and the next error is the product $\tilde g_{k+1} = \tilde g_k\Exp(\tau_k)$.}
\label{fig:diagram}
\end{figure}
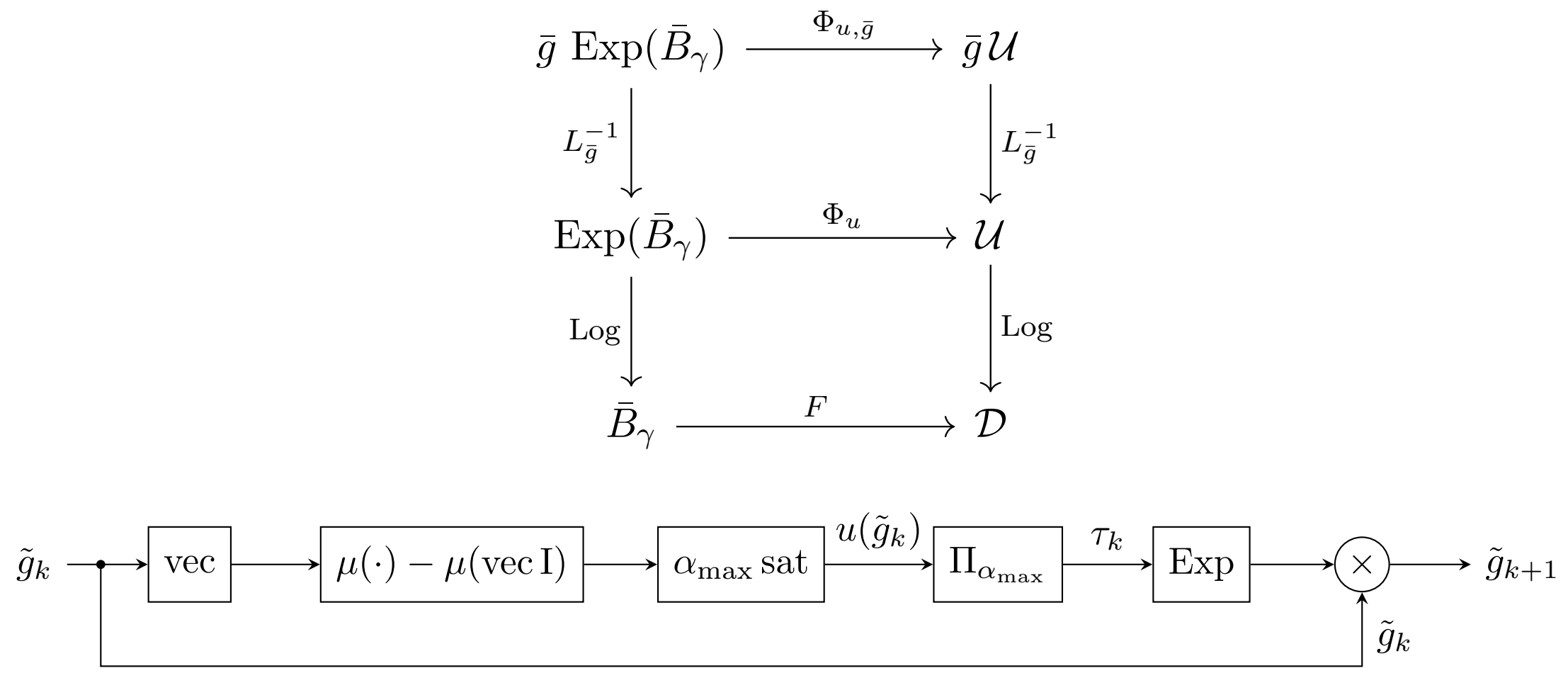

\begin{proposition}
\label{prop:domain}
Suppose that Assumption~\ref{ass:chart} holds with chart radius $\rG$, that $u : \mathsf{G} \to \R^{m}$ is continuous and that $\gamma + \amax < \rG$. Then $\Phi_u(\Exp z)$ is in $\Uch$ for every $z$ in $\bar B_\gamma$, and the map $F = \Log{}\circ\Phi_u\circ\Exp$ from $\bar B_\gamma$ to $\Dch$ is continuous.
\end{proposition}
\begin{proof}
Let $z$ be in $\bar B_\gamma$ and put $\tau = \proj(u(\Exp z))$, so $\|\tau\| \le \amax$. Part (i) of Assumption~\ref{ass:chart} gives
\[
  \rot\bigl(\Phi_u(\Exp z)\bigr) = \rot\bigl(\Exp(z)\Exp(\tau)\bigr) \le \rot(\Exp z) + \rot(\Exp\tau) \le \|z\| + \|\tau\| \le \gamma + \amax < \rG ,
\]
so $\Phi_u(\Exp z)$ is in $\Uch$. The map $F$ is the composite of the continuous maps $\Exp$, $\Phi_u$ and $\Log$ (Assumption~\ref{ass:chart}(iii)), so it is continuous. Figure~\ref{fig:diagram} shows one step of the closed loop in states, errors and coordinates.
\end{proof}

\subsection{A second fixed point on compact groups}
A region of attraction cannot contain a second fixed point. On a compact connected group such as $\SO(n)$, with coordinates whose inner product is invariant under conjugation, every corrected policy whose goal is locally exponentially stable has one (Proposition~\ref{prop:lefschetz}), by the Lefschetz fixed-point theorem, as \citet{kalabic2017mpc} argue for model predictive control. So the region of attraction of such a closed loop is never the whole group. The argument uses the Lefschetz number $L(f)$ of a smooth map $f$ from a compact oriented manifold to itself \citep[Chapter 3, Section 4]{guillemin1974differential}. A Lie group is orientable, because left translation of a basis of $\g$ gives a global frame of its tangent bundle. Three facts are needed. The number $L(f)$ equals $L(f')$ when $f$ and $f'$ are smoothly homotopic. The number $L$ of the identity map is the Euler characteristic. If every fixed point $x$ of $f$ satisfies $\det(Df(x) - \Id) \ne 0$, then $L(f)$ is the sum of the signs of these determinants. \citet[Theorem 4]{kalabic2017mpc} use the theorem in the same way, for a map on a compact manifold that is homotopic to the identity and has finitely many fixed points. The nondegeneracy hypothesis below follows from local exponential stability, since $\spr(DF(0)) < 1$ (Lemma~\ref{lem:spr}) excludes $1$ from the spectrum. The proposition needs the left-invariant metric of this appendix to be bi-invariant, which holds exactly when the inner product of the coordinates is invariant under conjugation by $\mathsf{G}$ \citep[Proposition 21.3]{gallier2020differential}. It also needs Assumption~\ref{ass:chart} with $\rot$ the distance of that metric cut off at $\rG$, which Proposition~\ref{prop:general} supplies for every compact group. On $\SO(n)$ the inner product $\tfrac12\tr(X^{\top}Y)$ is invariant under conjugation, its distance from $\Id$ is the rotation angle, and Assumption~\ref{ass:chart} holds with $\rG = \pi$ and the rotation angle cut off at $\pi$, since the cutoff changes neither $\Uch$ nor part (i). The proposition says nothing about a group that is not compact, such as $\SE(n)$.

\begin{proposition}
\label{prop:lefschetz}
Let $\mathsf{G}$ be a compact connected matrix Lie group of positive dimension whose coordinate inner product is invariant under conjugation, so the left-invariant metric of this appendix is bi-invariant, and suppose that Assumption~\ref{ass:chart} holds with $\rot(g) = \min\{\dist(\Id, g), \rG\}$. Let $u : \mathsf{G} \to \R^{m}$ be continuous, smooth near $\Id$, with $u(\Id) = 0$ and $\det(DF(0) - \Id) \ne 0$. Then $\Phi_u$ has a fixed point other than $\Id$.
\end{proposition}

\begin{proof}
Let $\dist$ be the bi-invariant distance, so $\Uch = \{g : \dist(\Id, g) < \rG\}$, and let $\tau(\tilde g) = \proj(u(\tilde g))$. Two elements $h_1$, $h_2$ at distance less than $\rG$ are joined by a unique minimal geodesic, which depends continuously on its end points. Indeed, $\mathsf{G}$ is compact, hence complete, so a minimal geodesic from $\Id$ to an element $h$ exists \citep[Theorems 16.16 and 16.17]{gallier2020differential}. For a bi-invariant metric it is a one-parameter subgroup $t \mapsto \Exp(t\xi_h)$ \citep[Proposition 21.21]{gallier2020differential}, and its length is $\|\xi_h\| = \dist(\Id, h)$, because the speed of this curve in the left-invariant metric is the constant $\|\xi_h\|$. If $\dist(\Id, h) < \rG$, then $\xi_h$ lies in $B_{\rG} \subseteq \Dch$ and $h$ lies in $\Uch$, so $\xi_h = \Log h$ by the injectivity of $\Exp$ on $\Dch$. Left translation by $h_1$ is an isometry, so the unique minimal geodesic from $h_1$ to $h_2$ is $t \mapsto h_1\Exp(t\Log(h_1^{-1}h_2))$, which is continuous in $(h_1, h_2, t)$ because $\Log$ is smooth.

\emph{Step 1: fixed points.} Since $\tau(\tilde g)$ lies in $\bar B_{\amax} \subseteq \Dch$ and $\Exp$ is injective on $\Dch$ with $\Exp(0) = \Id$, the equation $\Phi_u(\tilde g) = \tilde g$ holds exactly when $\tau(\tilde g) = 0$, that is, when $u(\tilde g) = 0$.

\emph{Step 2: a smooth replacement.} Suppose, for a contradiction, that $\Id$ is the only fixed point of $\Phi_u$. Choose a closed ball $\calN$ about $\Id$ on which $\Phi_u$ is smooth. The function $\tilde g \mapsto \dist(\Phi_u(\tilde g), \tilde g)$ is continuous and positive on the compact set $\mathsf{G} \setminus \operatorname{int}\calN$, so it has a positive minimum $\delta_0$ there. The step bound of Appendix~\ref{app:lie} gives $\delta_0 \le \amax < \rG$. Regard $\Phi_u$ as a map into $\R^{d\times d}$. By the Whitney approximation theorem \citep[Theorem 6.15]{lee2000smooth}, there is a smooth map into $\R^{d\times d}$ that equals $\Phi_u$ on $\calN$ and is arbitrarily close to it. Composing with the smooth retraction of a tubular neighborhood of $\mathsf{G}$ onto $\mathsf{G}$ \citep[Theorem 6.17 and Proposition 6.18]{lee2000smooth} gives a smooth map $\Phi_{\mathrm s} : \mathsf{G} \to \mathsf{G}$ that equals $\Phi_u$ on $\calN$. If the approximation is close enough, then $\dist(\Phi_{\mathrm s}(\tilde g), \Phi_u(\tilde g)) < \delta_0/2$ for every $\tilde g$, because the retraction is uniformly continuous on a compact neighborhood of $\mathsf{G}$ and the Riemannian distance is continuous in the matrix topology.

\emph{Step 3: the Lefschetz number of $\Phi_{\mathrm s}$.} For every $\tilde g$ outside the interior of $\calN$, the triangle inequality gives $\dist(\Phi_{\mathrm s}(\tilde g), \tilde g) > \delta_0 - \delta_0/2 > 0$, so $\Phi_{\mathrm s}$ has no fixed point outside $\calN$. On $\calN$ it equals $\Phi_u$, whose only fixed point is $\Id$. In exponential coordinates the derivative of $\Phi_{\mathrm s}$ at $\Id$ is $DF(0)$, and $\det(DF(0) - \Id) \ne 0$ by hypothesis. Hence $L(\Phi_{\mathrm s}) = \operatorname{sign}\det(DF(0) - \Id) = \pm1$ \citep[Chapter 3, Section 4, p.~121]{guillemin1974differential}.

\emph{Step 4: a homotopy to the identity.} Moving each point $\Phi_{\mathrm s}(\tilde g)$ to $\Phi_u(\tilde g)$ along the minimal geodesic, which is unique since $\delta_0/2 < \rG$, is a homotopy from $\Phi_{\mathrm s}$ to $\Phi_u$. The map $(\tilde g, t) \mapsto \tilde g\Exp(t\tau(\tilde g))$, $0 \le t \le 1$, is a homotopy from the identity to $\Phi_u$. So $\Phi_{\mathrm s}$ is homotopic to the identity, and since both maps are smooth, they are smoothly homotopic \citep[Proposition 6.20]{lee2000smooth}. Hence $L(\Phi_{\mathrm s})$ equals the Euler characteristic $\chi(\mathsf{G})$ \citep[Chapter 3, Section 4, p.~120]{guillemin1974differential}.

\emph{Step 5: the Euler characteristic.} Fix $\xi$ with $0 < \|\xi\| < \rG$, which exists since $m \ge 1$. The map $h \mapsto h\Exp\xi$ is smoothly homotopic to the identity through $h \mapsto h\Exp(t\xi)$, and it has no fixed point because $\Exp\xi \ne \Id$ by the injectivity of $\Exp$ on $\Dch$. A map without fixed points has Lefschetz number zero, so $\chi(\mathsf{G}) = 0$ \citep[Chapter 3, Section 4, p.~120]{guillemin1974differential}. This contradicts $L(\Phi_{\mathrm s}) = \pm1$.
\end{proof}

The hypothesis $u(\Id) = 0$ cannot be dropped. For example, if $\dist(\Id, h_0) < \amax$ and $h_0 \ne \Id$, then the continuous policy $u(\tilde g) = \Log(\tilde g^{-1}h_0\tilde g)$ gives $\Phi_u(\tilde g) = h_0\tilde g$, which has no fixed point. For a fixed network, the other fixed points are the zeros of $\mu(\vecop\tilde g) - \mu(\vecop\Id)$, so they do not depend on $\amax$ or on $\gamma$. The continuous-time obstruction of \citet[Theorem 1]{bhat2000topological} is related, and \citet{chaturvedi2011rigid} survey the topology of attitude control. On a group that is not compact the Lefschetz number is not available. \citet[Remark 6.22]{bullo2005geometric} record that the region of attraction of a continuous vector field is contractible, which rules out a global region of attraction in continuous time on every group that is not contractible, such as $\SE(n)$, but gives no fixed point.

\subsection{Exponential stability and the spectral radius}
The proof of the next lemma follows the linearization argument of \citet[Theorem 7.1]{lasalle1986stability} and adds the exponential estimate.

The goal is locally exponentially stable if there exist $r_0 > 0$, $c_0 \ge 1$ and $q$ in $[0, 1)$ such that, for every $z$ in $B_{r_0}$ and every $k \ge 0$, the iterate $F^{k}(z)$ lies in $\bar B_\gamma$ and satisfies $\|F^{k}(z)\| \le c_0q^{k}\|z\|$.

\begin{lemma}
\label{lem:spr}
Suppose that $u$ is smooth near $\Id$ and that $u(\Id) = 0$. Then the goal is locally exponentially stable if and only if $\spr(DF(0)) < 1$.
\end{lemma}

\begin{proof}
Let $A = DF(0)$. Suppose first that $\spr(A) < 1$. Then $A^{k} \to 0$, and the series $P = \sum_{k \ge 0}(A^{\top})^{k}A^{k}$ converges to the symmetric positive definite solution of $A^{\top}PA - P = -\Id$ \citep[Theorem 6.3]{lasalle1986stability}. It satisfies $P \succeq \Id$, because its first term is $\Id$. Let $b \ge 1$ be the largest eigenvalue of $P$, and let $V_P(z) = z^{\top}Pz$ be the linearization-based baseline, so $\|z\|^{2} \le V_P(z) \le b\|z\|^{2}$. Since $F(z) = Az + o(\|z\|)$, the change $V_P(F(z)) - V_P(z)$ equals $-\|z\|^{2} + o(\|z\|^{2})$. Hence there is $r_1$ in $(0, \gamma]$ with $V_P(F(z)) - V_P(z) \le -\tfrac12\|z\|^{2} \le -V_P(z)/(2b)$ for $\|z\| < r_1$. Put $\varpi = 1 - 1/(2b)$, which lies in $[\tfrac12, 1)$. The set $\{V_P < r_1^{2}\}$ lies in $B_{r_1}$ and is mapped into itself by $F$, and on it $V_P(F^{k}(z)) \le \varpi^{k}V_P(z)$. Every $z$ in $B_{r_1/\sqrt b}$ lies in this set and therefore satisfies $\|F^{k}(z)\|^{2} \le V_P(F^{k}(z)) \le \varpi^{k}b\|z\|^{2}$. The definition of local exponential stability holds with $r_0 = r_1/\sqrt b$, $c_0 = \sqrt b$ and $q = \sqrt\varpi$.

\looseness -1 Conversely, suppose that the goal is locally exponentially stable. Since $F(0) = 0$, the chain rule gives $D(F^{k})(0) = A^{k}$. For a unit vector $y$, the bound of the definition gives $\|A^{k}y\| = \lim_{t \to 0^{+}}\|F^{k}(ty)\|/t \le c_0q^{k}$, so $\|A^{k}\| \le c_0q^{k}$. Hence $\spr(A) = \lim_{k \to \infty}\|A^{k}\|^{1/k} \le q < 1$. \end{proof}

\subsection{The accumulated cost}
\begin{proposition}
\label{prop:maximal}
Suppose that $u$ is continuous and that the goal is locally exponentially stable. Then $\calA = \{\tilde g \in \mathsf{G} : \Vm(\tilde g) < \infty\}$, and $\Vm$ is a maximal Lyapunov function of $\Phi_u$.
\end{proposition}
\begin{proof}
Let $r_0$, $c_0$ and $q$ be as in the definition of local exponential stability. The set $\Exp(B_{r_0})$ is a neighborhood of $\Id$, because $\Exp$ maps the open set $\Dch$ onto the open set $\Uch$ with a continuous inverse (Assumption~\ref{ass:chart}). Let $\tilde g$ be in $\calA$. Since $\Phi_u^{k}(\tilde g) \to \Id$, there is $k_0$ with $\Phi_u^{k_0}(\tilde g) = \Exp(z_0)$ for some $z_0$ in $B_{r_0}$. The iterates $F^{j}(z_0)$ stay in $\bar B_\gamma$, so the identity $\Exp{}\circ F = \Phi_u{}\circ\Exp$ applies to each of them and gives $\Phi_u^{k}(\tilde g) = \Exp(F^{k - k_0}(z_0))$ for $k \ge k_0$. The cost bound then gives
\[
  \ell\bigl(\Phi_u^{k}(\tilde g)\bigr)^{2} \le \Cl\bigl\|F^{k - k_0}(z_0)\bigr\|^{2} \le \Cl c_0^{2}q^{2(k - k_0)}\|z_0\|^{2}.
\]
The tail of the accumulated cost is dominated by a convergent geometric series, so $\Vm(\tilde g) < \infty$. Conversely, if $\Vm(\tilde g) < \infty$, then $\|\Phi_u^{k}(\tilde g) - \Id\|_F = \ell(\Phi_u^{k}(\tilde g)) \to 0$, so $\tilde g$ is in $\calA$.

It remains to check the other properties of a maximal Lyapunov function. Continuity of $\Vm$ on $\calA$, and $\Vm(\tilde g_j) \to \infty$ along every sequence $(\tilde g_j)$ in $\calA$ that converges to a point outside $\calA$, are Proposition~\ref{prop:zubov}(ii), whose proof uses only the first paragraph of this proof. The map $\Phi_u$ fixes $\Id$ and $\ell(\Id) = 0$, so $\Vm(\Id) = 0$. For $\tilde g$ in $\calA \setminus \{\Id\}$, the first term of the accumulated cost gives $\Vm(\tilde g) \ge \ell(\tilde g)^{2} > 0$, and \eqref{eq:recursion} gives $\Vm(\Phi_u(\tilde g)) = \Vm(\tilde g) - \ell(\tilde g)^{2} < \Vm(\tilde g)$.
\end{proof}

Exponential stability cannot be replaced by convergence alone: on $\SO(2)$ a policy with $F(z) = z/(1 + z^{2})$ near $0$ drives every small $z_0 \ne 0$ to $0$, but with $F^{k}(z_0)^{2} \sim 1/(2k)$, so $\Vm(\Exp z_0) = \infty$ although $\Exp z_0$ is in $\calA$.

\begin{corollary}
\label{cor:unique}
Let $\mathcal{J} \subseteq \calA$ satisfy $\Phi_u(\mathcal{J}) \subseteq \mathcal{J}$, with $\Id$ in $\mathcal{J}$, and let $\tilde W : \mathcal{J} \to \R$ be continuous at $\Id$, with $\tilde W(\Id) = 0$ and $\tilde W(\tilde g) = \ell(\tilde g)^{2} + \tilde W(\Phi_u(\tilde g))$ for every $\tilde g$ in $\mathcal{J}$. Then $\tilde W = \Vm$ on $\mathcal{J}$.
\end{corollary}
\begin{proof}
Let $\tilde g$ be in $\mathcal{J}$. Applying the equation $k$ times along the trajectory, which stays in $\mathcal{J}$, gives $\tilde W(\tilde g) = \sum_{j < k}\ell(\Phi_u^{j}(\tilde g))^{2} + \tilde W(\Phi_u^{k}(\tilde g))$. As $k \to \infty$, $\Phi_u^{k}(\tilde g) \to \Id$ because $\tilde g$ is in $\calA$, so the last term tends to $\tilde W(\Id) = 0$, and the partial sums converge to $\Vm(\tilde g)$.
\end{proof}

For a map on $\R^{n}$ with an exponentially stable fixed point, \citet[Theorem 2.8]{giesl2007determination} defines the same sum with a running cost such as the squared distance to the fixed point, and shows that it solves the one-step equation on the basin and has compact sublevel sets. With a safe set added, Proposition~\ref{prop:maximal} and Corollary~\ref{cor:unique} correspond to \citet[Theorems 2, 4, 5, 6 and 10]{serry2025safe}. The two results above state the construction on a matrix Lie group, with the running cost $\ell^{2}$.

\begin{lemma}
\label{lem:local}
Suppose that the goal is locally exponentially stable and that $u$ is twice continuously differentiable near $\Id$. Then $\Wm(z) = z^{\top}\Pm z + O(\|z\|^{3})$ as $z \to 0$.
\end{lemma}
\begin{proof}
Let $r_0$, $c_0$ and $q$ be as in the definition of local exponential stability. Replacing $q$ by $\max\{q, \tfrac12\}$ keeps the definition true, so assume $q > 0$. Since $F(0) = 0$, the chain rule gives $D(F^{k})(0) = A^{k}$, and for a unit vector $y$ the bound of the definition gives $\|A^{k}y\| = \lim_{t \to 0^{+}}\|F^{k}(ty)\|/t \le c_0q^{k}$, so $\|A^{k}\| \le c_0q^{k}$ for every $k$. The map $F$ is twice continuously differentiable near $0$, because $\Exp$ and $\Log$ are smooth, $u$ is twice continuously differentiable near $\Id$, and the projection and the saturation are the identity near the goal. Likewise $z \mapsto \ell(\Exp z)^{2}$ is smooth, with value $0$ and gradient $0$ at $z = 0$. Taylor's theorem therefore gives a radius $r$ in $(0, r_0]$ and constants $C_1$, $C_2$ with $\|F(w) - Aw\| \le C_1\|w\|^{2}$ and $|\ell(\Exp w)^{2} - w^{\top}Qw| \le C_2\|w\|^{3}$ for $\|w\| \le r$. Fix $z$ with $\|z\| < r/c_0$, so every iterate satisfies $\|F^{k}(z)\| \le c_0q^{k}\|z\| < r$. Put $d_k = F^{k}(z) - A^{k}z$, so $d_0 = 0$. Applying $F$ to $F^{k}(z) = A^{k}z + d_k$ and subtracting $A^{k+1}z$ gives $d_{k+1} = (F(F^{k}(z)) - AF^{k}(z)) + Ad_k$, and unrolling this recursion gives $d_k = \sum_{j < k}A^{k-1-j}(F(F^{j}(z)) - AF^{j}(z))$. The bounds $\|A^{k-1-j}\| \le c_0q^{k-1-j}$ and $\|F(F^{j}(z)) - AF^{j}(z)\| \le C_1\|F^{j}(z)\|^{2} \le C_1c_0^{2}q^{2j}\|z\|^{2}$ and the geometric series give
\[
  \|d_k\| \le C_1c_0^{3}\|z\|^{2}\,q^{k-1}\sum_{j < k}q^{j} \le C_3\,q^{k}\|z\|^{2}, \qquad C_3 = \frac{C_1c_0^{3}}{q(1 - q)} .
\]
Next, substituting $F^{k}(z) = A^{k}z + d_k$ into $w^{\top}Qw$ and collecting terms gives
\[
  \ell\bigl(\Exp F^{k}(z)\bigr)^{2} = (A^{k}z)^{\top}Q(A^{k}z) + 2(A^{k}z)^{\top}Qd_k + d_k^{\top}Qd_k + \bigl(\ell(\Exp F^{k}(z))^{2} - F^{k}(z)^{\top}QF^{k}(z)\bigr) .
\]
The last three terms are bounded by $2\|Q\|c_0C_3q^{2k}\|z\|^{3}$, by $\|Q\|C_3^{2}q^{2k}\|z\|^{4}$ and by $C_2c_0^{3}q^{3k}\|z\|^{3}$. For $\|z\| \le 1$ their sum is at most $C_4q^{2k}\|z\|^{3}$ with $C_4 = 2\|Q\|c_0C_3 + \|Q\|C_3^{2} + C_2c_0^{3}$. Last, the terms are summed over $k$. Since the iterates stay in $\bar B_\gamma$, the identity $\Exp{}\circ F = \Phi_u{}\circ\Exp$ applies to each of them, and $\Wm(z) = \sum_{k \ge 0}\ell(\Exp F^{k}(z))^{2}$. The first terms sum to $z^{\top}\bigl(\sum_{k \ge 0}(A^{\top})^{k}QA^{k}\bigr)z = z^{\top}\Pm z$, and the remainders sum to at most $C_4\|z\|^{3}/(1 - q^{2})$.
\end{proof}

If $\Wm(z) = z^{\top}\Pm z + T_3(z) + o(\|z\|^{3})$ with a nonzero cubic form $T_3$, the candidate cannot match $T_3$, since $V(z) = z^{\top}\Pm z + O(\|z\|^{4})$ for every choice of weights (Section~\ref{sec:candidate}). Near the goal $V - \Wm = -T_3(z) + o(\|z\|^{3})$, which is of order $\|z\|^{3}$ along every direction in which $T_3$ does not vanish, and the network can reduce the difference only away from the goal, where $s$ is not small.

\subsection{Approximation by the candidate}
The candidate \eqref{eq:cert} vanishes at $0$, since $V_0(0) = 0$ for every choice of weights, and it is positive elsewhere: $V(z) \ge (1 - s(z))\lambda_{\min}(\Pm)\|z\|^{2} + s(z)\kappa^{2}\ell(\Exp z)^{2}$, and the first term on the right is positive for $z \ne 0$, since $s(z) < 1$.
\begin{proposition}
\label{prop:ua}
Let $\calK \subseteq \Dch$ be compact with $0$ in $\calK$, let $\kappa > 0$, and let $V^{*} : \calK \to \R$ be continuous with $V^{*}(0) = 0$ and $V^{*}(z) \ge \kappa^{2}\ell(\Exp z)^{2}$ for every $z$ in $\calK$. Then for every $\varepsilon > 0$ there is a ReLU network $\psi : \R^{m} \to \R$ with one hidden layer whose network part satisfies $|V_0(z) - V^{*}(z)| < \varepsilon$ for every $z$ in $\calK$.
\end{proposition}

\begin{proof}
The construction follows \citet[Proposition 2]{akhtar2026learning}. The function $\psi_1 = (V^{*} - \kappa^{2}\ell(\Exp\cdot)^{2})^{1/2}$ is continuous on $\calK$, with $\psi_1(0) = 0$, and $\psi_1^{2} + \kappa^{2}\ell(\Exp\cdot)^{2} = V^{*}$. Let $M_1 = \max_{\calK}\psi_1$ and choose $\eta_1$ in $(0, 1]$ with $4\eta_1(M_1 + 1) < \varepsilon$. A ReLU network $\psi$ with one hidden layer and $|\psi - \psi_1| < \eta_1$ on $\calK$ exists, because the rectifier is not a polynomial, by \citet[Theorem 1]{leshno1993multilayer} and \citet[Theorem 3.1]{pinkus1999approximation}. With $a(z) = \psi(z) - \psi(0)$, the network part is $V_0 = a^{2} + \kappa^{2}\ell(\Exp\cdot)^{2}$, and the bound $|a(z) - \psi_1(z)| \le |\psi(z) - \psi_1(z)| + |\psi(0) - \psi_1(0)| < 2\eta_1$ holds, so
\[
  |V_0(z) - V^{*}(z)| = |a(z) - \psi_1(z)|\,|a(z) + \psi_1(z)| < 2\eta_1\,(2M_1 + 2\eta_1) \le 4\eta_1(M_1 + 1) < \varepsilon .
\]
\end{proof}

Propositions~\ref{prop:ua} and \ref{prop:anchored} and Corollary~\ref{cor:exists} are stated for a network with one output. A network $\psi : \R^{m} \to \R^{n_\psi}$ whose further $n_\psi - 1$ outputs are zero has the same network part, so the three results hold for every $n_\psi$.

\begin{proposition}
\label{prop:anchored}
Let $\calK \subseteq \Dch$ be compact with $0$ in $\calK$, and let $V^{*} : \calK \to \R$ be continuous with $V^{*}(z) = z^{\top}\Pm z + o(\|z\|^{2})$ at $0$ and $V^{*} \ge \ell(\Exp\cdot)^{2}$ on $\calK$. Then for every $\varepsilon > 0$ there are $r_a > 0$, $\kappa > 0$ and a ReLU network $\psi : \R^{m} \to \R$ with one hidden layer whose candidate \eqref{eq:cert} satisfies $|V(z) - V^{*}(z)| < \varepsilon$ for every $z$ in $\calK$.
\end{proposition}
\begin{proof}
Let $V_{\mathrm q}(z) = z^{\top}\Pm z$ be the quadratic term of the candidate, so $(1 - s(z))\,V_{\mathrm q}(z) \le \lambda_{\max}(\Pm) r_a^{2}\|z\|^{2}/(\|z\|^{2} + r_a^{2}) \le \lambda_{\max}(\Pm) r_a^{2}$ for every $z$ and every $r_a$. Fix $\kappa$ in $(0, \tfrac12]$ and $r_a > 0$ with $\lambda_{\max}(\Pm) r_a^{2} < \varepsilon/2$.

The target of the network part comes first. Put $V_{\mathrm r} = V^{*} - (1 - s)V_{\mathrm q}$ and $V_1 = V_{\mathrm r}/s$ on $\calK \setminus \{0\}$. Since $V_{\mathrm r}/s = V_{\mathrm q} + (V^{*} - V_{\mathrm q})(\|z\|^{2} + r_a^{2})/\|z\|^{2}$ and $V^{*} - V_{\mathrm q} = o(\|z\|^{2})$, the function $V_1$ tends to $0$ at $0$. With $V_1(0) = 0$, the function $V_1$ is continuous on $\calK$, because $s$ is continuous and positive away from $0$. Put $V_2 = \max\{V_1, \kappa^{2}\ell(\Exp\cdot)^{2}\}$. This function is continuous on $\calK$, vanishes at $0$, and is at least $\kappa^{2}\ell(\Exp\cdot)^{2}$. Proposition~\ref{prop:ua} gives a network $\psi$ with $|V_0 - V_2| < \varepsilon/2$ on $\calK$.

Next, substituting $V = (1 - s)V_{\mathrm q} + sV_0$ and $V^{*} = (1 - s)V_{\mathrm q} + sV_1$ gives $V - V^{*} = s(V_0 - V_2) + s(V_2 - V_1)$. The first term is at most $\varepsilon/2$ in absolute value, since $s \le 1$. The second term is nonnegative and equals $(s\kappa^{2}\ell(\Exp\cdot)^{2} - V_{\mathrm r})_{+} = (s\kappa^{2}\ell(\Exp\cdot)^{2} + (1 - s)V_{\mathrm q} - V^{*})_{+}$. Since $s \le 1$, $\kappa^{2} \le \tfrac14$ and $\ell(\Exp\cdot)^{2} \le V^{*}$, the first two summands are at most $V^{*}/4 + \lambda_{\max}(\Pm) r_a^{2}$, so the positive part is at most $\lambda_{\max}(\Pm) r_a^{2} < \varepsilon/2$. Hence $|V - V^{*}| < \varepsilon$ on $\calK$. \end{proof}

The accumulated cost satisfies the hypotheses of Proposition~\ref{prop:anchored} on $V^{*}$: for a compact set $\mathcal{M} \subseteq \calA \cap \Uch$ that contains $\Id$, the set $\calK = \Log(\mathcal{M})$ is compact, and $V^{*} = \Wm$ is continuous on it (Proposition~\ref{prop:zubov}), satisfies $V^{*} \ge \ell(\Exp\cdot)^{2}$ by the first term of the accumulated cost, and equals $z^{\top}\Pm z + o(\|z\|^{2})$ by Lemma~\ref{lem:local}.

The proof shows that, for every $\kappa$ in $(0, \tfrac12]$ and every $r_a > 0$, a network exists with $|V - V^{*}| < \varepsilon/2 + \lambda_{\max}(\Pm) r_a^{2}$ on $\calK$, so it needs a small radius $r_a$ for a given accuracy.

\subsection{Conditions on a sublevel set}
\begin{proposition}
\label{prop:cut}
Suppose that $u$ is continuous, that $V : \Dch \to [0, \infty)$ is continuous, and that decrease, entry and containment hold. Then $F(\Omega_\rho) \subseteq \Omega_\rho$, and for every $z$ in $\Omega_\rho$ there is $k \ge 0$ with $F^{k}(z)$ in $N_b$.
\end{proposition}

\begin{proof}
Let $z$ be in $\Omega_\rho$. If $z$ is in $N_b$, then $F(z)$ is in $\Omega_\rho$ by entry. Otherwise, containment gives $F(z)$ in $\bar B_\gamma$ and decrease gives $V(F(z)) < V(z) \le \rho$, so again $F(z)$ is in $\Omega_\rho$. This proves $F(\Omega_\rho) \subseteq \Omega_\rho$. The map $F$ is continuous (Proposition~\ref{prop:domain}) and $V$ is continuous, so $\Omega_\rho$ is closed in the compact ball $\bar B_\gamma$ and hence compact. The set $N_b$ is open in $\Omega_\rho$, so $\Omega_\rho \setminus N_b$ is compact. If it is empty, then the claim holds with $k = 0$. Otherwise, the continuous function $-\delta V$ is positive on it and so has a positive minimum $\delta_{\min}$ there. Suppose that $F^{k}(z)$ is not in $N_b$ for $k = 0, \dots, K$. Then $0 \le V(F^{K+1}(z)) \le V(z) - (K + 1)\delta_{\min} \le \rho - (K + 1)\delta_{\min}$, which is impossible for $K + 1 > \rho/\delta_{\min}$.
\end{proof}

\begin{corollary}
\label{cor:region}
Suppose, in addition, that $F^{k}(z) \to 0$ for every $z$ in $N_b$. Then $\Exp(\Omega_\rho) \subseteq \calA$.
\end{corollary}

\begin{proof}
Let $z$ be in $\Omega_\rho$. By the proposition, the iterates $F^{k}(z)$ stay in $\Omega_\rho \subseteq \bar B_\gamma$ and reach $N_b$ at some step $k_1$, after which they converge to $0$. Since $\Phi_u^{k}(\Exp z) = \Exp(F^{k}(z))$ and $\Exp$ is continuous, $\Phi_u^{k}(\Exp z) \to \Id$. \end{proof}

The separate treatment of a neighborhood of the goal follows \citet[Corollary 3.13]{giesl2007determination}, where the decrease of the approximation is shown outside a sublevel set of a quadratic Lyapunov function of the linearization, and that sublevel set is handled by the quadratic function itself. Here the set $N_b$ is handled by the entry condition and the convergence of its points. Under the hypotheses of Theorem~\ref{thm:main} both hold for every level in $[c - 3\varepsilon, c - \varepsilon]$, by parts (b) and (c), and the experiments check both on samples. Lemma~\ref{lem:localdec} below gives both without Theorem~\ref{thm:main} when $N_b$ lies in the set $D_{\mathrm{loc}}$ of the lemma and $\rho_{\mathrm{loc}} \le \rho$, but its radius is not computed.

\begin{lemma}
\label{lem:localdec}
Let (H1)--(H3) hold and let $V$ be the candidate \eqref{eq:cert} with any weights. Then $\delta V(z) = -z^{\top}Qz + O(\|z\|^{3})$ as $z \to 0$. Hence there is $r_{\mathrm{loc}} > 0$ with $F(B_{r_{\mathrm{loc}}}) \subseteq \bar B_\gamma$ and $\delta V(z) < 0$ for $0 < \|z\| \le r_{\mathrm{loc}}$. For $\rho_{\mathrm{loc}} = \lambda_{\min}(\Pm)\,r_a^{2}r_{\mathrm{loc}}^{2}/(r_a^{2} + r_{\mathrm{loc}}^{2})$, the set $D_{\mathrm{loc}} = \{z \in \bar B_\gamma : V(z) < \rho_{\mathrm{loc}}\}$ lies in $B_{r_{\mathrm{loc}}}$, is mapped into itself by $F$, and every trajectory from it converges to $0$. If $N_b \subseteq D_{\mathrm{loc}}$ and $\rho_{\mathrm{loc}} \le \rho$, then entry holds and the hypothesis of Corollary~\ref{cor:region} is satisfied.
\end{lemma}

\begin{proof}
\looseness -1 Write $V = V_{\mathrm q} + s(V_0 - V_{\mathrm q})$ with $V_{\mathrm q}(z) = z^{\top}\Pm z$. Since $\psi$ is Lipschitz with constant $\Lip(\psi)$, $|s(V_0 - V_{\mathrm q})| \le (\Lip(\psi)^{2} + \Cl\kappa^{2} + \lambda_{\max}(\Pm))\|z\|^{4}/r_a^{2}$. Since $F$ is smooth near $0$ with $F(0) = 0$, $F(z) = Az + O(\|z\|^{2})$, so $V_{\mathrm q}(F(z)) = z^{\top}A^{\top}\Pm Az + O(\|z\|^{3})$, and \eqref{eq:stein} gives the expansion. The bound $V(z) \ge (1 - s(z))\lambda_{\min}(\Pm)\|z\|^{2}$ from the positivity argument above is increasing in $\|z\|$ and equals $\rho_{\mathrm{loc}}$ at $\|z\| = r_{\mathrm{loc}}$, so $D_{\mathrm{loc}} \subseteq B_{r_{\mathrm{loc}}}$. For $z$ in $D_{\mathrm{loc}} \setminus \{0\}$, $V(F(z)) < V(z) < \rho_{\mathrm{loc}}$, so $F(z)$ is in $D_{\mathrm{loc}}$. The sequence $V(F^{k}(z))$ decreases to a limit, and every limit point $\bar z$ of the orbit lies in $\bar B_{r_{\mathrm{loc}}}$ with $\delta V(\bar z) = 0$ by continuity, so $\bar z = 0$. The last claim follows from $F(N_b) \subseteq D_{\mathrm{loc}} \subseteq \Omega_\rho$.
\end{proof}

The radius $r_{\mathrm{loc}}$ depends on the weights through the constants above and is not computed here. In continuous time, \citet[Proposition 3]{barreau2024unsupervised} and \citet[Proposition 9]{liu2026universal} prove the same local property for their candidates, and \citet[Corollary 3.13]{giesl2007determination} handles the ball around the fixed point by the quadratic term alone.

\subsection{Sublevel sets of the accumulated cost}
\begin{lemma}
\label{lem:sublevel}
Under (H1)--(H3), $\bar c > 0$, and for every $c$ in $(0, \bar c)$: (i) $K_c$ is compact and contains a neighborhood of $0$; (ii) $F(K_c) \subseteq K_c$, and $\Wm(F(z)) = \Wm(z) - \ell(\Exp z)^{2}$ for every $z$ in $K_c$; (iii) $\Exp(K_c) \subseteq \calA$ and $K_c \subseteq \calA_\gamma$; (iv) $\Wm(z) \ge \bar c$ for every $z$ in $\bar B_\gamma$ with $\|F(z)\| > \gamma$.
\end{lemma}

\begin{proof}
The function $\ell$ is continuous and vanishes only at $\Id$. Its sublevel sets $\{\tilde g \in \mathsf{G} : \ell(\tilde g) \le r\}$ with $r < 1$ are compact. Such a set is bounded, and every matrix $M$ in its closure in $\R^{d\times d}$ satisfies $\|M - \Id\|_F \le r < 1$ and is therefore invertible. Since $\mathsf{G}$ is closed in the group of invertible $d \times d$ matrices, the closure lies in $\mathsf{G}$. Decreasing compact sets with intersection $\{\Id\}$ eventually lie in every open neighborhood of $\Id$. The set $\Exp(B_\gamma)$ is such a neighborhood, because $\Exp$ maps $\Dch$ onto the open set $\Uch$ with a continuous inverse (Assumption~\ref{ass:chart}). So $\{\ell^{2} \le c_1\} \subseteq \Exp(B_\gamma)$ for some $c_1 > 0$. Let $z$ be in $K_{c_1}$. Since $\Vm \ge \ell^{2}$, \eqref{eq:recursion} gives $\ell(\Phi_u(\Exp z))^{2} \le \Vm(\Phi_u(\Exp z)) \le \Wm(z) \le c_1$, so $\Exp(F(z)) = \Phi_u(\Exp z)$ lies in $\Exp(B_\gamma)$. Both $F(z)$ and $B_\gamma$ lie in $\Dch$, on which $\Exp$ is injective, so $F(z)$ is in $B_\gamma$. Hence $F(K_{c_1}) \subseteq \bar B_\gamma$ and $\bar c \ge c_1 > 0$.

Let $c$ be in $(0, \bar c)$. The sets $K_c$ increase with $c$, so $F(K_c) \subseteq \bar B_\gamma$. (i) The function $\Vm$ is lower semicontinuous (Proposition~\ref{prop:zubov}(ii)), so $K_c$ is closed in the compact ball $\bar B_\gamma$. The function $\Vm$ is continuous on the open set $\calA$ with $\Vm(\Id) = 0$, so $\{\Vm < c\} \cap \Exp(B_\gamma)$ is an open neighborhood of $\Id$, and its image under $\Log$ is a neighborhood of $0$ inside $K_c$. (ii) Let $z$ be in $K_c$. Then $F(z)$ is in $\bar B_\gamma$, and since $\Exp(F(z)) = \Phi_u(\Exp z)$, \eqref{eq:recursion} gives $\Wm(F(z)) = \Vm(\Phi_u(\Exp z)) = \Wm(z) - \ell(\Exp z)^{2} \le c$. (iii) Points of $K_c$ have finite $\Wm$, so their images under $\Exp$ lie in $\calA$ by Proposition~\ref{prop:maximal}. For $z$ in $K_c$, the iterates $F^{k}(z)$ stay in $K_c \subseteq \bar B_\gamma$ by (ii), and $\Exp(F^{k}(z)) = \Phi_u^{k}(\Exp z) \to \Id$, so $F^{k}(z) \to 0$ because $\Log$ is continuous on $\Uch$. Hence $K_c \subseteq \calA_\gamma$. (iv) Let $z$ be in $\bar B_\gamma$ with $\|F(z)\| > \gamma$. If $\Wm(z) < \bar c$, then $z$ lies in $K_c$ for some $c$ in $(0, \bar c)$, and $F(K_c) \subseteq \bar B_\gamma$ contradicts $\|F(z)\| > \gamma$. So $\Wm(z) \ge \bar c$. \end{proof}

\subsection{The level rule and Theorem~\ref{thm:main}}
The proof of Theorem~\ref{thm:main} follows.
\begin{proof}
(a) Let $z$ be in $\Omega_{\rho_*}$. By (B), $z$ is in $K_c$, since otherwise $V(z) > c - \varepsilon > \rho_*$, and (A) gives $\Wm(z) \le V(z) + \varepsilon \le c - \varepsilon$. Conversely, if $z$ is in $K_{c - 3\varepsilon} \subseteq K_c$, then (A) gives $V(z) \le \Wm(z) + \varepsilon \le \rho_*$. (b) Let $\rho$ be in $(0, c - \varepsilon]$ and let $z$ be in $\Omega_\rho$. As in (a), (B) forces $z$ into $K_c$. Lemma~\ref{lem:sublevel}(ii) gives $F(z)$ in $K_c \subseteq \bar B_\gamma$, which is containment, and $\Wm(F(z)) = \Wm(z) - \ell(\Exp z)^{2}$. Applying (A) at $z$ and at $F(z)$ gives $\delta V(z) \le (\Wm(F(z)) + \varepsilon) - (\Wm(z) - \varepsilon) = 2\varepsilon - \ell(\Exp z)^{2}$. If $z$ is not in $N_b$, then $\ell(\Exp z) \ge r_b$, so $\delta V(z) \le -(r_b^{2} - 2\varepsilon) < 0$. Last, $\Exp(\Omega_\rho) \subseteq \Exp(K_c) \subseteq \calA$ by Lemma~\ref{lem:sublevel}(iii). (c) Let $\rho$ be in $[c - 3\varepsilon, c - \varepsilon]$ and $z$ in $N_b$, so $z$ is in $K_{c - 4\varepsilon}$ by hypothesis and $V(F(z)) \le \Wm(F(z)) + \varepsilon \le \Wm(z) + \varepsilon \le c - 3\varepsilon \le \rho$. With $F(z)$ in $\bar B_\gamma$ by (b), $F(z)$ is in $\Omega_\rho$. So decrease, entry and containment hold on $\Omega_\rho$, which is therefore positively invariant (Proposition~\ref{prop:cut}), and it is compact because $V$ is continuous. By (b), $\Omega_\rho$ lies in $K_c \subseteq \calA_\gamma$ (Lemma~\ref{lem:sublevel}(iii)). Last, $V$ is continuous with $V(0) = 0 < \rho$, so $\Omega_\rho$ contains a neighborhood of $0$.
\end{proof}

The level rule of Section~\ref{sec:learning} works in three stages. The first two stages use independent samples $\calZ_1$, drawn like the training points, and $\calZ_2$, drawn uniformly in the ball, and set the sample level $\rho_0$ to the smaller of $\max_{\calZ_1} V$ and $(1 - 10^{-9})$ times the smallest value of $V$ at a sampled point of $\calB$. The third stage, from $\rho = \rho_0$, is an adversarial search for points of $\calB$ in $\Omega_\rho \setminus N_b$, by gradient ascent on $\delta V$ and separately on $\|F(z)\| - \gamma$. Whenever it finds some, $\rho$ becomes $(1 - 10^{-9})$ times the smallest value of $V$ at them, and the search is repeated until a search round finds none, for at most 50 search rounds. The third stage is then run a second time in the same way, from the final level and with new starting points drawn from the same samples.

\begin{proposition}
\label{prop:levelrule}
Under the hypotheses of Theorem~\ref{thm:main}, $V(z) > \rho_*$ for every $z$ in $\calB$, so the level $\hat\rho$ chosen in Section~\ref{sec:learning} satisfies $\hat\rho \ge \min\{\max_{\calZ_1}V, (1 - 10^{-9})\rho_*\}$, whatever the samples and the points visited by the search.
\end{proposition}
\begin{proof}
A point $z$ of $\bar B_\gamma$ with $V(z) \le \rho_*$ lies in $\Omega_{\rho_*}$, on which decrease holds outside $N_b$ and containment holds, by Theorem~\ref{thm:main}(b). Such a point is not in $\calB$. Every sample in $\calB$ therefore has $V > \rho_*$, so the sample level satisfies $\rho_0 \ge \min\{\max_{\calZ_1}V, (1 - 10^{-9})\rho_*\}$. Each round of the third stage, in either run, sets the level to $(1 - 10^{-9})$ times the least value of $V$ at points of $\calB$ that it finds, and this value exceeds $\rho_*$. A round that finds no such point leaves the level unchanged.
\end{proof}

\subsection{When the chosen level gives the region sought}
By Theorem~\ref{thm:main}, decrease, entry and containment hold on every sublevel set of $V$ with level in $[c - 3\varepsilon, c - \varepsilon]$. By Proposition~\ref{prop:levelrule}, the chosen level lies in this interval if three conditions hold together: it does not exceed $c - \varepsilon$, $\varepsilon \ge 10^{-9}\rho_*$, and $\max_{\calZ_1}V \ge c - 3\varepsilon$. The reported region then has every property asked of $\Omega$ in Section~\ref{sec:loopdef}. The set $\Omega_\rho$ is compact because $V$ is continuous, and positively invariant by Proposition~\ref{prop:cut}. It is contained in $\calA_\gamma$ because its iterates stay in $\Omega_\rho \subseteq \bar B_\gamma$, $\Exp(\Omega_\rho) \subseteq \calA$ by Theorem~\ref{thm:main}(b), and $\Log$ is continuous on $\Uch$. It is a neighborhood of $0$ because (A) gives $K_{c - 4\varepsilon} \subseteq \Omega_{c - 3\varepsilon} \subseteq \Omega_\rho$, and $K_{c - 4\varepsilon}$ contains a neighborhood of $0$ by Lemma~\ref{lem:sublevel}(i), since $4\varepsilon < c$. The theorem gives no upper bound on the chosen level. If the rule stops above $c - \varepsilon$, then the reported region is supported only by the checks on samples and by search.

\subsection{Existence of a candidate}
\begin{corollary}
\label{cor:exists}
Let (H1)--(H3) hold, and let $c$ be in $(0, \bar c)$. For every $\varepsilon > 0$ there are $r_a > 0$, $\kappa > 0$ and a ReLU network $\psi : \R^{m} \to \R$ with one hidden layer whose candidate \eqref{eq:cert} satisfies (A) and (B). In addition, there is $r_b > 0$ such that the hypotheses of Theorem~\ref{thm:main} hold for every $\varepsilon < \min\{3c/16, r_b^{2}/2\}$.
\end{corollary}
\begin{proof}
Let $M_\ell = \max_{\bar B_\gamma}\ell(\Exp\cdot)^{2}$, put $C_* = \max\{c + 1, M_\ell\}$ and $V^{*} = \min\{\Wm, C_*\}$ on $\bar B_\gamma$. By Proposition~\ref{prop:zubov}(ii), $\Vm$ is continuous on the open set $\calA$ and tends to infinity at every point outside $\calA$. So $\min\{\Vm, C_*\}$ is continuous on $\mathsf{G}$, and $V^{*}$ is continuous on $\bar B_\gamma$. On $\{\Wm \le C_*\}$, $V^{*} = \Wm \ge \ell(\Exp\cdot)^{2}$, and elsewhere $V^{*} = C_* \ge M_\ell \ge \ell(\Exp\cdot)^{2}$. Near $0$, $V^{*} = \Wm$, so $V^{*}(z) = z^{\top}\Pm z + o(\|z\|^{2})$ by Lemma~\ref{lem:local}. Also $\bar B_\gamma \subseteq B_{\rG} \subseteq \Dch$, since $\gamma < \rG$. Proposition~\ref{prop:anchored} with $\calK = \bar B_\gamma$ gives a candidate with $|V - V^{*}| < \varepsilon$ on $\bar B_\gamma$. On $K_c$, $V^{*} = \Wm$, which gives (A). On $\bar B_\gamma \setminus K_c$, $V^{*} > c$, so $V > c - \varepsilon$, which gives (B). For the choice of $r_b$, the set $\{\Vm < c/4\}$ is an open neighborhood of $\Id$, and the sublevel sets $\{\ell \le r\}$ are compact and shrink to $\{\Id\}$ as $r \to 0$ (proof of Lemma~\ref{lem:sublevel}), so there is $r_b > 0$ with $\{\tilde g : \ell(\tilde g) \le r_b\} \subseteq \{\Vm < c/4\}$. If $\varepsilon < 3c/16$, then $c - 4\varepsilon > c/4$, so every $z$ in $\bar B_\gamma$ with $\ell(\Exp z) < r_b$ lies in $K_{c - 4\varepsilon}$, and $4\varepsilon < c$. If also $\varepsilon < r_b^{2}/2$, then all hypotheses of Theorem~\ref{thm:main} hold.
\end{proof}

\subsection{Proof of Proposition~\ref{prop:residual}}
\begin{proof}
For $z$ in $K_c$, Lemma~\ref{lem:sublevel} gives $F^{k}(z)$ in $K_c$ for every $k$ and $\Exp(F^{k}(z)) = \Phi_u^{k}(\Exp z) \to \Id$, so $F^{k}(z) \to 0$. If $K_c \setminus B_r$ is empty, put $k_r = 0$. Otherwise put $m_r = \min\{\ell(\Exp z)^{2} : z \in K_c \setminus B_r\}$, which is positive because the set is compact and $\ell(\Exp z)$ vanishes only at $z = 0$, and $k_r = \lfloor c/m_r \rfloor + 1$. If $F^{k}(z)$ lies outside $B_r$ for $k = 0, \dots, K$, then Lemma~\ref{lem:sublevel}(ii) gives $0 \le \Wm(F^{K+1}(z)) \le c - (K + 1)m_r$, so $K + 1 \le c/m_r$. Hence every $z$ in $K_c$ has some $k \le k_r$ with $F^{k}(z)$ in $B_r$. Put $e = V - \Wm$ on $K_c$. By Lemma~\ref{lem:sublevel}(ii), $\Wm(z) = \ell(\Exp z)^{2} + \Wm(F(z))$, so the hypothesis gives $|e(z) - e(F(z))| \le \eta$ for $z$ in $K_c$. For $z$ in $K_c$ and its $k \le k_r$, this bound along $z, F(z), \dots, F^{k}(z)$, which all lie in $K_c$, gives $|e(z)| \le k\eta + |e(F^{k}(z))| \le k_r\eta + \sup_{K_c \cap B_r}|e|$.
\end{proof}

A continuous-time counterpart, for approximate solutions of the Zubov equation with a ball around the origin removed, is \citet[Proposition 4]{liu2025physics}.

\subsection{The bound near the goal after Proposition~\ref{prop:residual}}
Let $r_0$, $c_0$ and $q$ be as in the definition of local exponential stability, let $r \le r_0$, and suppose that $\psi$ is Lipschitz with constant $\Lip(\psi)$. For $\|z\| < r$, $s(z) \le r^{2}/r_a^{2}$ and $\ell(\Exp z)^{2} \le \Cl\|z\|^{2}$ by the cost bound, so the candidate satisfies
\begin{align*}
  V(z) &\le \lambda_{\max}(\Pm)\|z\|^{2} + s(z)\bigl(\Lip(\psi)^{2} + \Cl\kappa^{2}\bigr)\|z\|^{2} \\
  &\le \bigl(\lambda_{\max}(\Pm) + (\Lip(\psi)^{2} + \Cl\kappa^{2})\,r^{2}/r_a^{2}\bigr)\|z\|^{2},
\end{align*}
and the proof of Proposition~\ref{prop:maximal} gives $\Wm(z) \le \Cl c_0^{2}\|z\|^{2}/(1 - q^{2})$. Both are nonnegative, so $|V - \Wm| \le \max\{V, \Wm\} \le \Gamma r^{2}$ on $B_r$, with
\[
  \Gamma = \max\bigl\{\lambda_{\max}(\Pm) + (\Lip(\psi)^{2} + \Cl\kappa^{2})\,r^{2}/r_a^{2},\ \Cl c_0^{2}/(1 - q^{2})\bigr\},
\]
which grows with $\Lip(\psi)^{2}/r_a^{2}$.

\subsection{What is not proved}
Section~\ref{sec:guarantee} states the main gap: no trained candidate is known to satisfy (A). The remaining gaps are the following. The spectral term only encourages a small Lipschitz constant, and no bound is derived from it. Gradient descent may stop at a poor local minimum. The refinement guarantees nothing, and its estimate of the region can be optimistic (Appendix~\ref{app:training}). Lemma~\ref{lem:localdec} gives decrease near the goal for every candidate, but its radius $r_{\mathrm{loc}}$ is not computed, so whether $N_b \subseteq D_{\mathrm{loc}}$ is not known. On $\SO(3)$, $\SE(2)$ and $\SE(3)$ the entry check and the convergence of $N_b$ therefore rest on samples and simulation. Corollary~\ref{cor:exists} is an existence result for wide networks with one hidden layer, not for the architecture of Table~\ref{tab:consts}, and its proof uses a small $r_a$, while the experiments fix $r_a = 0.3$. The sensitivity to $\gamma$ is not studied, and that to $r_a$ and $r_b$ only on a few closed loops (Appendix~\ref{app:ablations}). The numbers $\bar c$ and $\varepsilon$ are never computed, so the reported level may exceed $c - \varepsilon$. The theorem therefore motivates the loss terms but does not cover the reported regions.

\subsection{The Zubov transform}
Put $W_{\mathrm{z}} = 1 - e^{-\Vm}$, with $W_{\mathrm{z}} = 1$ at the points with $\Vm = \infty$. The function $W_{\mathrm{z}}$ takes values in $[0, 1]$, and $W_{\mathrm{z}}(\tilde g) = 0$ exactly when $\tilde g = \Id$. It satisfies the discrete-time Zubov equation
\begin{equation}
  1 - W_{\mathrm{z}}(\Phi_u(\tilde g)) = e^{\ell(\tilde g)^{2}}\bigl(1 - W_{\mathrm{z}}(\tilde g)\bigr) \qquad \text{for every } \tilde g \text{ in } \mathsf{G} .
  \label{eq:zuboveq}
\end{equation}
At a point with $\Vm(\tilde g) < \infty$, the equation is the one-step equation \eqref{eq:recursion} with $\Vm = -\log(1 - W_{\mathrm{z}})$, multiplied by $-1$ and exponentiated. At a point with $\Vm(\tilde g) = \infty$, the one-step equation gives $\Vm(\Phi_u(\tilde g)) = \infty$, and both sides are zero.

\begin{proposition}
\label{prop:zubov}
Suppose that $u$ is continuous and that the goal is locally exponentially stable. Then: (i) $\calA = \{\tilde g : W_{\mathrm{z}}(\tilde g) < 1\}$, and $\calA$ is open; (ii) $\Vm$ is lower semicontinuous on $\mathsf{G}$, $\Vm$ and $W_{\mathrm{z}}$ are continuous on $\calA$, and $W_{\mathrm{z}}(\tilde g_j) \to 1$ for every sequence $(\tilde g_j)$ in $\calA$ that converges to a point outside $\calA$; (iii) if $W' : \calA \to [0, 1)$ is continuous at $\Id$, vanishes at $\Id$ and satisfies \eqref{eq:zuboveq} on $\calA$, then $W' = W_{\mathrm{z}}$.
\end{proposition}

\begin{proof}
Let $r_0$, $c_0$ and $q$ be as in the definition of local exponential stability, and put $\mathcal{O} = \Exp(B_{r_0})$, an open neighborhood of $\Id$.

(i) Since $t \mapsto 1 - e^{-t}$ is increasing with value $1$ only at $t = \infty$, the first paragraph of the proof of Proposition~\ref{prop:maximal} gives $\calA = \{W_{\mathrm{z}} < 1\}$. A trajectory converges to $\Id$ exactly when it enters $\mathcal{O}$, by the definition of local exponential stability and the same paragraph. Hence $\calA = \bigcup_{k \ge 0}\Phi_u^{-k}(\mathcal{O})$, which is open because every $\Phi_u^{k}$ is continuous.

(ii) Let $\tilde g$ be in $\calA$, and let $k_0$ satisfy $\Phi_u^{k_0}(\tilde g) \in \mathcal{O}$. By continuity, a neighborhood $\mathcal{V}$ of $\tilde g$ satisfies $\Phi_u^{k_0}(\mathcal{V}) \subseteq \mathcal{O}$. As in the proof of Proposition~\ref{prop:maximal}, for $\tilde g'$ in $\mathcal{V}$ and $K \ge k_0$,
\[
  \textstyle\sum_{k \ge K}\ell\bigl(\Phi_u^{k}(\tilde g')\bigr)^{2} \le \Cl c_0^{2}r_0^{2}\,q^{2(K - k_0)}/(1 - q^{2}) .
\]
So the partial sums of the accumulated cost, which are continuous, converge to $\Vm$ uniformly on $\mathcal{V}$, and $\Vm$ is continuous on $\mathcal{V}$. The same holds for $W_{\mathrm{z}}$. Next, $\Vm$ is the supremum of its continuous partial sums, so it is lower semicontinuous on $\mathsf{G}$. If $\tilde g_j \to \tilde g^{*}$ with $\tilde g^{*}$ not in $\calA$, then $\Vm(\tilde g^{*}) = \infty$ by (i), so lower semicontinuity gives $\Vm(\tilde g_j) \to \infty$ and $W_{\mathrm{z}}(\tilde g_j) \to 1$.

(iii) Put $V' = -\log(1 - W')$, a function from $\calA$ to $[0, \infty)$. It is continuous at $\Id$, vanishes there, and satisfies $V'(\tilde g) = \ell(\tilde g)^{2} + V'(\Phi_u(\tilde g))$ on $\calA$. The set $\calA$ satisfies $\Phi_u(\calA) \subseteq \calA$ and contains $\Id$, so Corollary~\ref{cor:unique} gives $V' = \Vm$, hence $W' = W_{\mathrm{z}}$.
\end{proof}

This is the characterization of \citet[Theorems 2, 7 and 11 and Corollary 1]{serry2025safe}, which is stated there on $\R^{n}$ with a safe set. The proof above uses the continuity of $\Phi_u$, the exponential estimate and the cost bound. The application of Proposition~\ref{prop:anchored} to $\Wm$, after its proof and in Corollary~\ref{cor:exists}, needs the continuity of $\Vm$ on $\calA$ and the bound $\Vm \ge \ell^{2}$.

\section{The trained policies}
\label{app:policies}
\paragraph{The decision process.} The observation is $\vecop\tilde g$, the entries of the rotation and translation parts of the error, of size $n_o = 4, 9, 6, 12$ on the four groups. The scale is $\bar\phi = \pi$ on the rotation groups and $\bar\phi = (\pi^{2} + (\pi T_{\max}/2)^{2})^{1/2}$ on the rigid-motion groups, with $T_{\max} = 1$. The clip length is $\amax = 0.1\,\bar\phi$. Episodes have 50 steps and no terminal state. On $\SE(n)$ the reward is not bounded below, since the translation may leave the ball of radius $T_{\max}$ during an episode. The cost $\phi(\tilde g) = \|\Log\tilde g\|$ is undefined at the errors whose rotation part has rotation angle $\pi$, a set of measure zero. There the implementation computes the logarithm of the rotation part from the angle on $\SO(2)$ and $\SE(2)$ and from the unit quaternion on $\SO(3)$ and $\SE(3)$. The distribution $\mu_0$ has a uniform (Haar) rotation part and, on $\SE(n)$, a translation uniform on the ball of radius $T_{\max}$, so on $\SE(n)$ the coordinate ball combines radians with these units. The task is solvable from every initial error in $\Exp(B_{\rG})$. The reference policy $u_\star(\tilde g) = -\Log\tilde g$ reaches the goal in $\lceil\|\Log\tilde g\|/\amax\rceil$ steps, since $\tilde g\Exp(-t\Log\tilde g) = \Exp((1 - t)\Log\tilde g)$ \citep[Theorem 2.14]{hall2015lie} and the segment from $\Log\tilde g$ to $0$ lies in $B_{\rG} \subseteq \Dch$. It is realized by the action $a = -\Log\tilde g/\max\{\amax, \|\Log\tilde g\|\}$.

\begin{figure}[!ht]
\centering
\includegraphics[width=\linewidth]{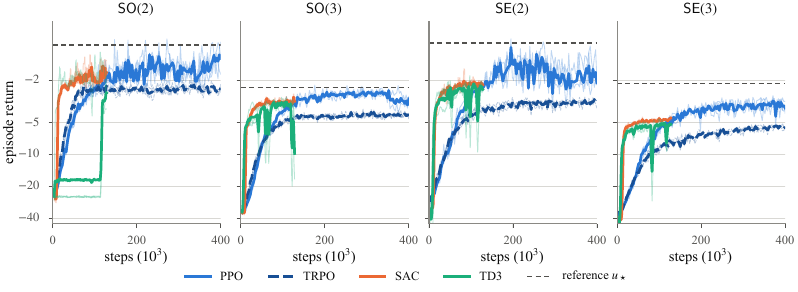}
\vspace{-14pt}
\caption{Learning curves: mean undiscounted return against environment steps, smoothed over 3200 steps. Thin lines are the seeds, thick lines their mean, and the horizontal dashed line the reference policy $u_\star$.}
\label{fig:curves}
\end{figure}

\paragraph{The algorithms.} PPO \citep{schulman2017ppo} maximizes a clipped surrogate of the expected discounted return, with generalized advantage estimation. TRPO \citep{schulman2015trpo} takes natural-gradient steps within a trust region of mean Kullback--Leibler divergence $0.01$, computed by conjugate gradients and a backtracking line search. Both use a Gaussian policy with mean $\mu(\vecop\tilde g)$ and a learned diagonal scale that does not depend on the state, with its logarithm initialized at $-1.5$. Both use the discount $\beta = 0.916$, the advantage parameter $0.8$, 32 environments with rollouts of 100 steps, and $4\times10^{5}$ environment steps. SAC \citep{haarnoja2018sac} and TD3 \citep{fujimoto2018td3} use $1.28\times10^{5}$ environment steps and a replay buffer of the same size, with discounts $0.975$ and $0.823$. The policy of SAC is a Gaussian with mean $\mu(\vecop\tilde g)$ and a scale that depends on the state, whose samples are squashed by $\tanh$. The scalar hyperparameters of PPO, SAC and TD3 (learning rates, discounts, batch sizes and update schedules) are those that \citet{schuck2026primer} tuned for the tangent-space action, and TRPO uses the critic and advantage estimation of PPO. A few choices differ from theirs. The actors of SAC and TD3 have two hidden layers of width 64, and all actor and critic networks use $\tanh$ activations. TD3 explores with Gaussian noise of standard deviation $0.1$ and updates its actor every second step, as in \citet{fujimoto2018td3}. The critics of PPO and TRPO have two hidden layers of width 64, and those of SAC and TD3 have width 256. Every actor has a $\tanh$ output. The environment saturates and clips every sampled action, so exploration always produces a valid step. The deterministic policy is the mean network $\mu$ of PPO, TRPO and SAC and the actor of TD3. Figure~\ref{fig:curves} shows the learning curves and Table~\ref{tab:policies} the trained policies.

\begin{table}[!ht]
\centering
\caption{The trained policies on 4096 initial errors from $\mu_0$: return over 50 steps against the reference policy, median final error before the correction of Section~\ref{sec:loop}, the offset $\|\mu(\vecop\Id)\|$, and $\spr(A)$ after it.}
\label{tab:policies}
\scriptsize\setlength{\tabcolsep}{3pt}
\begin{tabular}{lllrrrrrr}
\toprule
group & algorithm & seed & return & reference & rotation & translation & offset & $\spr(A)$ \\
\midrule
$\SO(2)$ & PPO & 0 & -1.607 & -1.411 & 0.0020 & 0.0000 & 0.0070 & 0.083 \\
&  & 1 & -1.488 & -1.411 & 0.0013 & 0.0000 & 0.0044 & 0.088 \\
&  & 2 & -1.538 & -1.411 & 0.0023 & 0.0000 & 0.0069 & 0.068 \\
& TRPO & 0 & -1.467 & -1.411 & 0.0007 & 0.0000 & 0.0025 & 0.112 \\
&  & 1 & -1.440 & -1.411 & 0.0016 & 0.0000 & 0.0052 & 0.023 \\
&  & 2 & -1.591 & -1.411 & 0.0057 & 0.0000 & 0.0165 & 0.083 \\
& SAC & 0 & -1.485 & -1.411 & 0.0033 & 0.0000 & 0.0080 & 0.211 \\
&  & 1 & -1.447 & -1.411 & 0.0008 & 0.0000 & 0.0031 & 0.161 \\
&  & 2 & -1.467 & -1.411 & 0.0020 & 0.0000 & 0.0076 & 0.315 \\
& TD3 & 0 & -3.701 & -1.411 & 0.1569 & 0.0000 & 1.0000 & 1.000 \\
&  & 1 & -1.815 & -1.411 & 0.0282 & 0.0000 & 0.0273 & 0.873 \\
&  & 2 & -1.657 & -1.411 & 0.0169 & 0.0000 & 0.0541 & 0.048 \\
\midrule
$\SO(3)$ & PPO & 0 & -3.319 & -2.385 & 0.0495 & 0.0000 & 0.1703 & 0.139 \\
&  & 1 & -2.832 & -2.385 & 0.0134 & 0.0000 & 0.0470 & 0.088 \\
&  & 2 & -3.261 & -2.385 & 0.0459 & 0.0000 & 0.1567 & 0.088 \\
& TRPO & 0 & -2.667 & -2.385 & 0.0010 & 0.0000 & 0.0033 & 0.085 \\
&  & 1 & -2.642 & -2.385 & 0.0005 & 0.0000 & 0.0017 & 0.094 \\
&  & 2 & -2.650 & -2.385 & 0.0027 & 0.0000 & 0.0091 & 0.110 \\
& SAC & 0 & -2.607 & -2.385 & 0.0133 & 0.0000 & 0.0459 & 0.123 \\
&  & 1 & -2.525 & -2.385 & 0.0076 & 0.0000 & 0.0224 & 0.081 \\
&  & 2 & -2.782 & -2.385 & 0.0262 & 0.0000 & 0.0777 & 0.079 \\
& TD3 & 0 & -9.206 & -2.385 & 0.4731 & 0.0000 & 1.4293 & 1.000 \\
&  & 1 & -3.539 & -2.385 & 0.0580 & 0.0000 & 0.1186 & 0.868 \\
&  & 2 & -8.001 & -2.385 & 0.3323 & 0.0000 & 1.7292 & 1.000 \\
\midrule
$\SE(2)$ & PPO & 0 & -1.745 & -1.353 & 0.0005 & 0.0015 & 0.0046 & 0.196 \\
&  & 1 & -1.918 & -1.353 & 0.0048 & 0.0062 & 0.0218 & 0.081 \\
&  & 2 & -1.560 & -1.353 & 0.0002 & 0.0021 & 0.0061 & 0.141 \\
& TRPO & 0 & -1.429 & -1.353 & 0.0036 & 0.0019 & 0.0123 & 0.110 \\
&  & 1 & -1.508 & -1.353 & 0.0064 & 0.0019 & 0.0209 & 0.085 \\
&  & 2 & -1.497 & -1.353 & 0.0007 & 0.0078 & 0.0230 & 0.169 \\
& SAC & 0 & -1.600 & -1.353 & 0.0205 & 0.0076 & 0.0415 & 0.394 \\
&  & 1 & -1.559 & -1.353 & 0.0178 & 0.0053 & 0.0416 & 0.335 \\
&  & 2 & -1.517 & -1.353 & 0.0048 & 0.0141 & 0.0285 & 0.374 \\
& TD3 & 0 & -1.740 & -1.353 & 0.0036 & 0.0308 & 0.0742 & 0.166 \\
&  & 1 & -1.914 & -1.353 & 0.0021 & 0.0435 & 0.1033 & 0.183 \\
&  & 2 & -1.661 & -1.353 & 0.0027 & 0.0238 & 0.0618 & 0.089 \\
\midrule
$\SE(3)$ & PPO & 0 & -3.308 & -2.170 & 0.0088 & 0.0114 & 0.0423 & 0.132 \\
&  & 1 & -3.048 & -2.170 & 0.0103 & 0.0030 & 0.0328 & 0.181 \\
&  & 2 & -3.023 & -2.170 & 0.0057 & 0.0043 & 0.0212 & 0.102 \\
& TRPO & 0 & -2.946 & -2.170 & 0.0049 & 0.0034 & 0.0143 & 0.394 \\
&  & 1 & -2.959 & -2.170 & 0.0056 & 0.0035 & 0.0175 & 0.426 \\
&  & 2 & -3.010 & -2.170 & 0.0033 & 0.0051 & 0.0123 & 0.485 \\
& SAC & 0 & -4.547 & -2.170 & 0.1468 & 0.0971 & 0.1996 & 4.006 \\
&  & 1 & -4.548 & -2.170 & 0.1547 & 0.0846 & 0.2361 & 3.452 \\
&  & 2 & -4.541 & -2.170 & 0.1451 & 0.1016 & 0.2653 & 4.637 \\
& TD3 & 0 & -5.050 & -2.170 & 0.1546 & 0.1186 & 0.7267 & 4.984 \\
&  & 1 & -4.975 & -2.170 & 0.1729 & 0.0652 & 0.4753 & 5.428 \\
&  & 2 & -5.438 & -2.170 & 0.1358 & 0.1671 & 1.2136 & 14.345 \\
\bottomrule
\end{tabular}~
\end{table}

\section{Training the candidate}
\label{app:training}
\looseness -1 Algorithm~\ref{alg:train-full} gives the training and refinement in full. The trained parameters are the weights of $\psi$ and a scalar $\tilde\kappa$, with $\kappa = \kappa_0 + \log(1 + e^{\tilde\kappa})$, so $\kappa$ stays above $\kappa_0$. The trained values of $\kappa$ lie between $0.24$ and $0.49$ over the 65 closed loops trained from scratch, the 39 policies and the 26 closed loops with an actuator error of the first screen of Appendix~\ref{app:screen}, and between $0.24$ and $0.41$ before the refinement. The matrix $A$ is estimated by central differences with step $10^{-6}$ and the matrix $Q$ by second-order central differences with step $10^{-4}$. Both are computed in double precision, as is $\Pm$, because central differences with these steps are inaccurate in single precision, while the network is trained in single precision. On $\SO(3)$ and $\SE(3)$ the computed $Q$ is $\operatorname{diag}(2, 2, 2)$ and $\operatorname{diag}(2, 2, 2, 1, 1, 1)$, as Appendix~\ref{app:lie} predicts. In the counterexample search, each ascent step moves every coordinate of a point by the same length, in the direction of the sign of the gradient of $\delta V$, and the point is then projected onto $\bar B_\gamma$. In six dimensions a volume-uniform point of the ball lies within radius $0.1$ of the goal with probability $2\times10^{-9}$, which is why half of the pool is drawn with a radius uniform on $[0, \gamma]$. The weights of $\psi$ are initialized from the LeCun normal distribution. In the third stage of the level rule, from $32768$ sample points in $\Omega_\rho \setminus N_b$, gradient ascent (Adam, 1000 steps) maximizes $\delta V$, and separately $\|F(z)\| - \gamma$, with a penalty for leaving $\Omega_\rho \setminus N_b$. The search runs for at most 50 rounds, and the search points are checked every 25 Adam steps. It is then run a second time from its final level with $32768$ new starting points. The entry sample of $N_b$ is simulated for 600 steps.

\begin{algorithm2e}[!ht]
\caption{Training and refinement of the candidate \eqref{eq:cert}, in full.}
\label{alg:train-full}
\DontPrintSemicolon
\SetKwFunction{Train}{Train}
\SetKwProg{Fn}{procedure}{}{}
\KwIn{the closed loop $F$ on $\bar B_\gamma$ and the constants of Table~\ref{tab:consts}}
\KwOut{the network $\psi$ and the coefficient $\kappa$}
estimate $A = DF(0)$ and $Q$ by central differences, compute $\Pm$ from \eqref{eq:stein}, initialize $\psi$ and $\kappa$, and set the pinned sets $\calC_{\mathrm{pin}}, \calX_{\mathrm{pin}} \gets \emptyset$\;
\Fn{\Train{$T_0$, $\mathit{ramp}$}}{
  set $\lambda \gets 1$, $\calC, \calX \gets \emptyset$ and a fresh Adam state, and fill a pool $\mathcal{P}$ with $32768$ points of $B_\gamma$, half volume-uniform and half radius-uniform, a quarter of which is replaced every $200$ iterations\;
    \For{$i = 0, 1, \dots, T_0 - 1$}{
    \If{$i$ is a multiple of $10$}{
      draw $2^{18}$ uniform points $U$ and $1024$ uniform points $U'$ of $B_\gamma$, and move $U'$ by $30$ sign-gradient ascent steps of length $\gamma/30$ on $\delta V$, each projected onto $\bar B_\gamma$\;
      add to $\calC$ the points of $U'$ at which $\delta V > -\zeta$ and those among the $2048$ points of $U$ with the largest $\delta V$ at which $\delta V > -\zeta$, and to $\calX$ the points $z$ of $U$ with $\ell(\Exp z) \ge r_b$ and $\|F(z)\| > \gamma$, keeping the newest $16384$ in each\;
    }
    draw $Z_i$ of size $4096$ from $\mathcal{P}$, $C_i$ and $X_i$ of size $1024$ from $\calC$ and $\calX$ (with $C_i$ from $Z_i$ while $\calC$ has fewer than $1024$ points, half of $C_i$ from $\calC_{\mathrm{pin}}$ and half of $X_i$ from $\calX_{\mathrm{pin}}$ when nonempty, and no exit term while $\calX$ and $\calX_{\mathrm{pin}}$ are empty), and $Y_i$ of size $2048$ from the shell $\gamma \le \|z\| \le \gamma + \amax$\;
    with $w_i = \min\{1, 5i/T_0\}$ if $\mathit{ramp}$ and $w_i = 1$ otherwise, take one Adam step on $\calL_i$ of \eqref{eq:objective}, skipped if the loss is not finite, then set $\lambda \gets \min\{\bar\lambda, (\lambda + \eta_\lambda\calL_{\mathrm{dec}}(Z_i))_{+}\}$, with the cap $\bar\lambda$ and the step $\eta_\lambda$ of Table~\ref{tab:consts}\;
  }
}
\Train{$T$, true}\;
\For{$j = 0, 1, \dots, 6$}{
  on $2\times10^{5}$ uniform and $2\times10^{5}$ radius-uniform fresh points of $\bar B_\gamma$, set the level to the smaller of the largest $V$ on the uniform points and $(1 - 10^{-9})$ times the smallest $V$ at a point where decrease or containment fails, then lower it by one round of the search of the level rule from $32768$ points of the region\;
    record the candidate and its quick estimate of the region, the fraction of $5\times10^{4}$ draws from $\mu_0$ in the resulting $\Omega_\rho$\;
  \lIf{no sampled point fails decrease and the search does not lower the level, or $j = 6$}{stop}
  add to $\calC_{\mathrm{pin}}$ the at most $4000$ counterexamples to decrease with the smallest $V$, each with $8$ copies $z \odot (\mathbf{1} + 0.05\,\upsilon) + 0.01\,\|z\|\,\upsilon'$ for standard normal vectors $\upsilon, \upsilon'$, with $\odot$ the entrywise product and $\mathbf{1}$ the vector of ones, kept outside $N_b$ and in $\bar B_\gamma$, and to $\calX_{\mathrm{pin}}$ at most $4000$ counterexamples to containment, keeping the newest $40000$ and $20000$\;
  \Train{$T'$, false}, starting from the current $\psi$ and $\kappa$\;
}
\Return{the recorded candidate with the largest quick estimate}
\end{algorithm2e}

\paragraph{The loss terms.} A counterexample is a point at which a condition fails. Fix a margin $\zeta$ with $0 < \zeta \le r_b^{2}$. The buffer $\calC$ collects the points with $\delta V > -\zeta$ found by the counterexample search of Algorithm~\ref{alg:train-full}, at which decrease fails or nearly fails. The buffer $\calX$ collects the points $z$ with $\ell(\Exp z) \ge r_b$ whose next point $F(z)$ leaves $\bar B_\gamma$, which lie in no region that satisfies containment. Let $\operatorname{mean}_Z$ be the average over a finite set $Z$, $Y$ a sample of the shell $\gamma \le \|z\| \le \gamma + \amax$ just outside the ball, $\nu_Y$ the median of $V$ on $Y$, and $c_{\mathrm{ex}} > 1$ a factor. The residual is $\calL_{\mathrm{rec}}(Z) = \operatorname{mean}_{Z}\bigl(\delta V(z) + \ell(\Exp z)^{2}\bigr)^{2}$. The hinge is $\calL_{\mathrm{dec}}(Z) = \operatorname{mean}\bigl(\delta V(z) + \zeta\bigr)_{+}$ over the points $z$ of $Z$ with $\ell(\Exp z) \ge r_b$. The shell term is $\calL_{\mathrm{sh}}(Y) = \operatorname{mean}_{Y}\bigl(1 - V(z)/\nu_Y\bigr)_{+}^{2}$. The exit term is $\calL_{\mathrm{ex}}(X) = \operatorname{mean}_{X}\bigl(1 - V(z)/(c_{\mathrm{ex}}\,\bar V_i)\bigr)_{+}^{2}$, with $\bar V_i = \max_{Z_i} V$ the largest value of $V$ on the batch $Z_i$ of iteration $i$. The spectral term is $\calR_{\mathrm{sp}} = \sum_j(\|\Theta_j\|_2 - \bar\sigma)_{+}^{2}$, summed over the weight matrices $\Theta_j$ of $\psi$, with a bound $\bar\sigma > 0$ on their spectral norms. Iteration $i$ draws a batch $Z_i$ from the pool $\mathcal{P}$ of Algorithm~\ref{alg:train-full}, batches $C_i$ from $\calC$ and $X_i$ from $\calX$, and a shell sample $Y_i$, and takes one Adam step on
\begin{equation}
  \calL_i = \lambda_i\calL_{\mathrm{dec}}(Z_i) + w_i\calL_{\mathrm{rec}}(Z_i) + w_{\mathrm{sh}}\calL_{\mathrm{sh}}(Y_i) + w_{\mathrm{ce}}\calL_{\mathrm{dec}}(C_i) + w_{\mathrm{ex}}\calL_{\mathrm{ex}}(X_i) + w_{\mathrm{sp}}\calR_{\mathrm{sp}} ,
  \label{eq:objective}
\end{equation}
with positive weights $w_{\mathrm{sh}}$, $w_{\mathrm{ce}}$, $w_{\mathrm{ex}}$ and $w_{\mathrm{sp}}$, and with $\lambda_i$ the value at iteration $i$ of the hinge weight $\lambda$. The median $\nu_Y$ and the maximum $\bar V_i$ are treated as constants when the gradient is computed. Each term has its own purpose. Besides pulling $V$ toward $\Wm$ (Section~\ref{sec:learning}), the residual $\calL_{\mathrm{rec}}$ constrains the one-step change $\delta V$, which the conditions test, and needs no derivative of $\Log$ and no simulation. The hinge $\calL_{\mathrm{dec}}$ asks for a decrease of at least $\zeta$ outside $N_b$. There $\Wm$ decreases by $\ell^{2} \ge r_b^{2} \ge \zeta$, so the hinge does not conflict with the residual. Its second copy acts on the counterexamples in $\calC$, because counterexamples are rare: for the final candidates, decrease or containment fails on at most $0.90$ percent of a uniform sample of the ball outside $N_b$ (Appendix~\ref{app:results}). The exit term raises $V$ on the points of $\calX$, the counterexamples to containment found during training, above $c_{\mathrm{ex}}$ times the largest value of $V$ in the batch. The shell term raises the smallest values of $V$ just outside the ball, so small sublevel sets stay inside the ball. Both are ratios, so they do not change the scale of $V$. The weight $w_i$ of the residual rises from zero to one over the first fifth of the iterations, and the weight $\lambda_i$ of the hinge is adapted by dual ascent \citep[Section 2.1]{boyd2011distributed}.

\paragraph{The refinement.} The search of the level rule (Section~\ref{sec:learning}) runs inside the region and finds sets of small measure on which decrease or containment fails, which the counterexample search of training misses. The refinement rounds of Algorithm~\ref{alg:train-full} feed them back. Each refinement round computes a quick estimate of the region on fresh samples, with random seeds that the level rule never uses. The refinement stops if no sampled point fails decrease and the search does not lower the level. Otherwise the counterexamples are pinned, and training continues for $T'$ iterations with $w_i = 1$ and half of every batch $C_i$ and $X_i$ drawn from the pinned points. The candidate with the largest quick estimate over the refinement rounds is kept. On 44 of the 65 closed loops trained from scratch the candidate before refinement was kept, and on the other 21 a refined one, after 1 to 6 refinement rounds. The quick estimate of a refinement round sets the level from $4\times10^{5}$ fresh points of the ball and one round of the search from $32768$ points, and measures the region on $5\times10^{4}$ draws from $\mu_0$, so it can miss counterexamples that the full level rule finds. On PPO seed 1 with gain $1.25$ on $\SE(3)$, the seven quick estimates, before the refinement and after each of its six refinement rounds, are $53.5$ to $64.8$ percent, and the full level rule gives $47.3$ percent for the kept candidate and $47.5$ percent for the candidate before refinement. Using the full level rule inside the refinement would let it see these counterexamples, at several times the cost, and is left for future work.

\paragraph{Closed loops of the later screens.} On the 29 closed loops with an actuator error from the second and third screens of Appendix~\ref{app:screen}, the candidate of the same policy without the actuator error is first given a new level by the level rule on the closed loop with the actuator error. If its region is within 2 percentage points of the oracle region, it is kept, as on 14 closed loops. Otherwise training restarts from it for 1000 iterations, with $\Pm$ recomputed for the closed loop with the actuator error, and is followed by the refinement. The larger of the two regions is kept, which is the refitted one on all 15. The refit takes $2.5$ to $4.3$ hours, against $2.9$ to $6.7$ hours from scratch on the first screen, and keeping the candidate costs only its evaluation, 41 to 156 minutes. On five closed loops of the first screen, this procedure, with one run of the search, gives a region within $0.5$ percentage points of the reported region of the candidate trained from scratch on four. On the fifth, PPO seed 1 with misalignment $15^\circ$ on $\SE(3)$, it gives $4.1$ points less, and the refinement again keeps a candidate whose quick estimate, $64.3$ percent, exceeds the $44.5$ percent given by the level rule. On PPO seed 1 on $\SE(3)$ with the first per-axis gains, the refinement produced non-finite losses from early in its third round on, and these steps were skipped, so the kept candidate of its last round has the parameters reached early in the third. A candidate kept without a refit has the quadratic term of the closed loop without the actuator error, so it does not agree with the accumulated cost to second order. Near the goal its region rests on the sampled entry check, which holds on all 14, and on an eigenvalue check. On all 14, the largest eigenvalue of $A^{\top}\Pm A - \Pm$, with the new $A$ and the kept $\Pm$, lies between $-0.997$ and $-0.529$. The leading term of $\delta V$ is therefore negative definite, and the proof of Lemma~\ref{lem:localdec} applies with this matrix in place of $-Q$, so the candidate still decreases strictly near the goal. Its fit, $2.6$ to $45.4$ percent, compares it with the accumulated cost of a closed loop on which it was not trained.

\begin{table}[!ht]
\centering
\caption{Constants of Algorithm~\ref{alg:train-full} and of the level rule, the same on every group except the number of iterations.}
\label{tab:consts}
\scriptsize\setlength{\tabcolsep}{1.8pt}
\begin{tabular}{@{}llll@{}}
\toprule
candidate & $\psi$: ReLU, widths $256, 256$, $n_\psi = 16$ & $\kappa_0 = 0.05$, initial $\kappa = 0.3$ & $r_a = 0.3$, $\gamma = \pi - \pi/8$ \\
optimizer & Adam, rate $10^{-3}$ & $T = 8000$ ($\SO(n)$), $4000$ ($\SE(n)$) & batch $|Z_i| = 4096$ \\
training points & pool of $32768$, half radius-uniform & a quarter replaced every $200$ iterations & $|Y_i| = 2048$ \\
loss & $\zeta = 0.01$, $r_b = 0.1$ & $w_{\mathrm{sh}} = 0.3$, $w_{\mathrm{sp}} = 10^{-3}$, $\bar\sigma = 2$ & $\eta_\lambda = 0.02$, $\bar\lambda = 50$ \\
counterexamples & every $10$ iterations, $|U| = 2^{18}$ & top $2048$ kept, $|U'| = 1024$ & $30$ ascent steps, $16384$ kept \\
& $|C_i| = |X_i| = 1024$ & $w_{\mathrm{ce}} = 5$, $w_{\mathrm{ex}} = 1$, $c_{\mathrm{ex}} = 1.2$ & \\
refinement & at most $6$ rounds, $T' = 1000$ & $4\times10^{5}$ fresh samples per round & $8$ copies per counterexample \\
level & $|\calZ_1| = 8\times10^{5}$, $|\calZ_2| = 4\times10^{5}$ & up to $10^{5}$ samples of $N_b$ & simulations of $600$ steps \\
& third stage, twice: $32768$ starts & $1000$ Adam steps, rate $3\times10^{-3}$ & penalty weight $100$, at most $50$ rounds \\
oracle & $3000$ simulated steps per point & $2\times10^{5}$ sample points & converged if error below $0.05$ \\
\bottomrule
\end{tabular}
\end{table}

\section{Results of every closed loop}
\label{app:results}
\begin{figure}[!ht]
\centering
\includegraphics[width=0.9\linewidth]{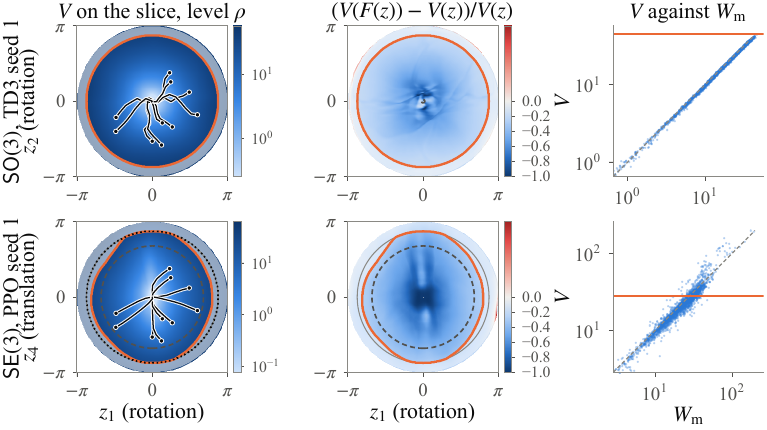}
\vspace{2pt}
\caption{Two policies on a coordinate slice through the goal: TD3 seed 1 on $\SO(3)$ (top, $z_1$ and $z_2$) and PPO seed 1 on $\SE(3)$ (bottom, $z_1$ and $z_4$). Left: the candidate on a logarithmic scale, with the boundaries of $\Omega_\rho$ (orange), the oracle region (dotted) and the region of the squared-norm baseline (dashed), the ball $\bar B_\gamma$ (gray) and trajectories. On top, $\|z\|^{2}$ increases along some trajectories, and the region of the squared-norm baseline is too small to draw. Middle: the one-step change of $V$ relative to $V$, blue for decrease and red for increase. Right: the candidate against the accumulated cost, with the level $\rho$ (orange).}
\label{fig:slices}
\end{figure}

Table~\ref{tab:main} summarizes the 39 policies that pass the stability test by group and algorithm, Figure~\ref{fig:slices} shows two of them on a slice through the goal, Table~\ref{tab:all} gives each of them, Table~\ref{tab:fail} the nine that fail the test, and Table~\ref{tab:pert-all} every evaluated closed loop with an actuator error. The level $\rho$ is the final level of the candidate, and the rounds are those of the first and the second run of the adversarial search, each ending with a round that found nothing. The fit of a candidate is the median of $|V - \Wm|/\Wm$ over the points $z$ with $\ell(\Exp z) \ge r_b$ among $2\times10^{4}$ points drawn in the ball, with $\Wm$ from 600 simulated steps. The three TD3 policies on the rotation groups that fail the stability test have offsets $\|\mu(\vecop\Id)\|$ of $1.00$, $1.43$ and $1.73$ (Table~\ref{tab:policies}), close to $1$, $\sqrt2$ and $\sqrt3$. The output of the $\tanh$ layer of the actor on $\SO(2)$, and two or three of its three outputs on $\SO(3)$, are then close to $\pm1$ at the goal, and the derivative of $\tanh$ is nearly zero there. The derivative of $u$ at the goal is therefore nearly zero, and $A$ is close to the identity. This explanation is inferred from the offsets and the spectral radii and is not proved.

\begin{table}[!ht]
\centering
\caption{Regions of the 39 policies by group and algorithm, all entries except the number of seeds in percent: the ball, the oracle region, the regions of $V$, $V_d$ and $V_P$, the same three regions as fractions of the oracle region, and the fit, the median of $|V - \Wm|/\Wm$ over the points with $\ell(\Exp z) \ge r_b$ among $2\times10^{4}$ points drawn in the ball. A range covers the seeds. The ball is measured on the samples of $\mu_0$, and on $\SO(2)$ it holds exactly $87.5$ percent.}
\label{tab:main}
\scriptsize\setlength{\tabcolsep}{3.5pt}
\begin{tabular}{@{}llrrrrrrrrrr@{}}
\toprule
group & algorithm & seeds & ball & oracle & $V$ & $V_d$ & $V_P$ & $V$/oracle & $V_d$/oracle & $V_P$/oracle & fit \\
\midrule
$\SO(2)$ & PPO & 3 & 87.6 & 87.6 & 87.6 & 87.6 & 87.6 & 100.0 & 100.0 & 100.0 & 0.0--0.2 \\
& TRPO & 3 & 87.6 & 86.8--87.6 & 86.1--87.6 & 85.9--87.6 & 85.9--87.6 & 99.3--100.0 & 99.0--100.0 & 99.0--100.0 & 0.1--2.1 \\
& SAC & 3 & 87.6 & 87.6 & 87.6 & 87.6 & 87.6 & 100.0 & 100.0 & 100.0 & 0.1--1.3 \\
& TD3 & 2 & 87.6 & 87.6 & 87.6 & 87.6 & 87.6 & 100.0 & 100.0 & 100.0 & 0.2--2.0 \\
\midrule
$\SO(3)$ & PPO & 3 & 75.4 & 75.4 & 75.4 & 75.4 & 75.4 & 100.0 & 100.0 & 100.0 & 0.3--0.5 \\
& TRPO & 3 & 75.4 & 75.4 & 75.4 & 75.4 & 75.4 & 100.0 & 100.0 & 100.0 & 0.2--0.5 \\
& SAC & 3 & 75.4 & 75.4 & 75.4 & 75.4 & 75.4 & 100.0 & 100.0 & 100.0 & 0.2--0.3 \\
& TD3 & 1 & 75.4 & 75.4 & 75.4 & 0.0 & 0.0 & 100.0 & 0.0 & 0.0 & 1.2 \\
\midrule
$\SE(2)$ & PPO & 3 & 82.0 & 82.0 & 82.0 & 60.0--82.0 & 60.4--82.0 & 100.0 & 73.3--100.0 & 73.7--100.0 & 0.8--3.6 \\
& TRPO & 3 & 82.0 & 82.0 & 82.0 & 76.0--82.0 & 75.6--82.0 & 100.0 & 92.7--100.0 & 92.3--100.0 & 0.5--3.7 \\
& SAC & 3 & 82.0 & 82.0 & 82.0 & 82.0 & 82.0 & 100.0 & 100.0 & 100.0 & 0.4 \\
& TD3 & 3 & 82.0 & 82.0 & 82.0 & 82.0 & 82.0 & 100.0 & 100.0 & 100.0 & 0.7--1.0 \\
\midrule
$\SE(3)$ & PPO & 3 & 65.1 & 65.1 & 52.5--65.1 & 32.7--45.1 & 32.2--46.3 & 80.6--100.0 & 50.2--69.3 & 49.5--71.0 & 2.4--3.9 \\
& TRPO & 3 & 65.1 & 65.1 & 65.1 & 64.3--65.1 & 62.3--65.1 & 100.0 & 98.7--100.0 & 95.6--100.0 & 1.1--1.3 \\
\bottomrule
\end{tabular}~
\end{table}

\begin{table}[!ht]
\centering
\caption{The nine policies that fail the stability test, with $\spr(A)$ after the correction of Section~\ref{sec:loop}. They are not evaluated.}
\label{tab:fail}
\scriptsize
\begin{tabular}{@{}lllr@{}}
\toprule
group & algorithm & seed & $\spr(A)$ \\
\midrule
$\SO(2)$ & TD3 & 0 & 1.000 \\
$\SO(3)$ & TD3 & 0 & 1.000 \\
$\SO(3)$ & TD3 & 2 & 1.000 \\
$\SE(3)$ & SAC & 0 & 4.006 \\
$\SE(3)$ & SAC & 1 & 3.452 \\
$\SE(3)$ & SAC & 2 & 4.637 \\
$\SE(3)$ & TD3 & 0 & 4.984 \\
$\SE(3)$ & TD3 & 1 & 5.428 \\
$\SE(3)$ & TD3 & 2 & 14.345 \\
\bottomrule
\end{tabular}~
\end{table}

\paragraph{Checks and fixed points.} Over the 94 closed loops, no region of the candidate, of the two quadratic baselines or of the oracle contains a sampled initial error that fails to converge in 600 steps. For the candidate, the entry condition holds at every sampled point of $N_b$, and every such point converges in simulation. For the final candidates, decrease or containment fails on at most $0.90$ percent of a uniform sample of the ball outside $N_b$ before the level is chosen. For the candidate and both quadratic baselines, each run of the adversarial search ended with a search round that found nothing on all 94 closed loops, after at most 13 of the 50 search rounds allowed in the first run and 6 in the second. The first run lowered the level of the candidate below the value set by the samples on 59 closed loops, all on $\SE(2)$ and $\SE(3)$, and the second lowered it further on 26, by more than one percentage point of $\mu_0$ on two. On two closed loops, PPO seed 1 on $\SE(3)$ without an actuator error and PPO seed 0 on $\SE(3)$ with misalignment $30^\circ$, a run of the search with $r_b = 0.2$ (Appendix~\ref{app:ablations}) ended below the level of the first run, and on PPO seed 1 below that of the second run as well, at points that also fail the conditions with $r_b = 0.1$. On PPO seed 1 it found two points at the edge of the ball whose next step leaves it, which neither run found. On both closed loops, the second run starts from this level instead, and on PPO seed 0 it lowers the level further. Every policy on the compact groups $\SO(2)$ and $\SO(3)$ that passes the stability test has a second fixed point (Proposition~\ref{prop:lefschetz} in Appendix~\ref{app:proofs}). For TRPO seed 2 on $\SO(2)$ it lies inside the ball, at $\|z\| = 2.696 < \gamma = 2.749$, and the region of the candidate stops before it (Table~\ref{tab:all}). The regions of the two quadratic baselines extend $8.4\times10^{-6}$ in $\|z\|$ past it, so they contain this fixed point and are not inner estimates. The excess has $\mu_0$-measure about $1.3\times10^{-6}$, which the $2\times10^{5}$ samples do not resolve. For the other policies on the compact groups it lies outside the ball. Every initial error drawn from $\mu_0$ converges on 37 of the 39 policies, and $99.978$ and $99.987$ percent on PPO seeds 0 and 1 on $\SE(3)$, and the initial errors that do not converge lie outside the ball. The fit of every candidate is the same with $\Wm$ from 3000 simulated steps as with 600.

The oracle decreases by exactly $\ell^{2}$ at every step, so only containment and convergence can fail, and its level is set by the first two stages alone, on a fixed subsample of $2\times10^{5}$ of their points and without the search. The oracle region is therefore checked less stringently than the regions of the candidate and the two quadratic baselines and can be slightly optimistic. It can also be smaller than the part of the ball that converges, since converging trajectories may leave the ball. On PPO seed 0 with misalignment $40^\circ$ on $\SE(3)$, every sampled initial error in the ball converges, but the oracle region is $63.5$ percent against $65.1$ for the ball, since some trajectories that converge leave the ball on the way.

For the candidate, the entry condition was checked on $10^{5}$ sampled points of $N_b$ on $\SO(2)$ and $\SE(2)$, on $88577$ on $\SO(3)$ and on $33064$ on $\SE(3)$, where fewer points of the sample fall in $N_b$. The decrease near the goal was also sampled. On 120 shells of radius between $10^{-4}$ and $\gamma$, with 20000 directions each and a gradient ascent of $\delta V(z)/\|z\|^{2}$ on each shell, or on $\SO(2)$ on a grid of $1.1\times10^{6}$ radii, $\delta V(z) \le -0.5\,\lambda_{\min}(Q)\|z\|^{2}$ holds at every tested point with $\|z\| \le 0.105$. This holds on all 94 closed loops, including the 14 with a kept candidate. For comparison, $N_b$ lies in the ball of radius $0.071$ on $\SO(n)$ and $0.100$ on $\SE(n)$. On 89 of the 94 the closed loop also contracts in the norm of $P$ at every tested point of the smallest sublevel set of $V_P$ that contains $N_b$. The radius $0.105$ is sampled, not a bound. Away from the goal, on the 60 closed loops on which a run of the search lowered the level, it found counterexamples missed by the $1.2\times10^{6}$ samples of the first two stages. On all 94 closed loops every initial error drawn from $\mu_0$ that lies in the ball converges, and on 69 of them so does every point of a uniform sample of the ball, so there the convergence checks show consistency rather than test the regions. On the other 25, all with PPO seeds 1 and 2 or TRPO seed 1 on $\SE(2)$ and $\SE(3)$, between $0.002$ and $0.64$ percent of the uniform sample fails to converge, and no region contains any of these points. The sampled and searched Lyapunov conditions are the checks that limit the regions. Two further results are recorded for completeness. First, the entry condition fails for the linearization-based baseline on TD3 seed 1 on $\SO(3)$, at 93 of 39145 sampled points of $N_b$, and its region there contains one of the $2\times10^{5}$ sampled initial errors. This region is reported as computed by the level rule, with the counterexamples. Second, a function built from the state-value critic $v_{\mathrm c}$ of a PPO or TRPO policy, $v_{\mathrm c}(\Id) - v_{\mathrm c}(\Exp z)$, was put through the level rule, with one run of the search, on the 50 closed loops trained from scratch whose policy has such a critic. It gives no region on 48 of them, one of the $2\times10^{5}$ sampled initial errors on PPO seed 1 on $\SO(2)$, and $0.3$ percent on PPO seed 2 on $\SO(2)$, where the entry condition fails at all 13469 sampled points of $N_b$ in it. On five closed loops on $\SO(2)$ its search used all 50 rounds, and the last round still found points. It is not a baseline of the paper.

\begin{table}[!ht]
\centering
\caption{Every policy that passes the stability test: $\spr(A)$, the $\mu_0$-measure in percent of the oracle region and of the regions of $V$, $V_d$ and $V_P$, the region of $V$ as a percentage of the oracle region and of the uniform points of the ball, the level $\rho$ of $V$, the rounds of the first and the second run of the adversarial search for $V$, and the fit of $V$ in percent, the median of $|V - \Wm|/\Wm$ over the points with $\ell(\Exp z) \ge r_b$ among $2\times10^{4}$ points drawn in the ball, all percentages rounded to one decimal. Differences quoted in the text are computed before rounding. The dagger marks the closed loop whose second run starts from a lower level, as explained in the text.}
\label{tab:all}
\scriptsize\setlength{\tabcolsep}{3pt}\renewcommand{\arraystretch}{0.97}
\begin{tabular}{@{}lllrrrrrrrrrr@{}}
\toprule
group & algorithm & seed & $\spr(A)$ & oracle & $V$ & $V_d$ & $V_P$ & $V$/oracle & $V$ (ball) & $\rho$ & rounds & fit \\
\midrule
$\SO(2)$ & PPO & 0 & 0.083 & 87.6 & 87.6 & 87.6 & 87.6 & 100.0 & 100.0 & 36.2 & 1+1 & 0.0 \\
& PPO & 1 & 0.088 & 87.6 & 87.6 & 87.6 & 87.6 & 100.0 & 100.0 & 34.4 & 1+1 & 0.1 \\
& PPO & 2 & 0.068 & 87.6 & 87.6 & 87.6 & 87.6 & 100.0 & 100.0 & 35.6 & 1+1 & 0.2 \\
& TRPO & 0 & 0.112 & 87.6 & 87.6 & 87.6 & 87.6 & 100.0 & 100.0 & 34.2 & 1+1 & 0.1 \\
& TRPO & 1 & 0.023 & 87.6 & 87.6 & 87.6 & 87.6 & 100.0 & 100.0 & 34.1 & 1+1 & 0.3 \\
& TRPO & 2 & 0.083 & 86.8 & 86.1 & 85.9 & 85.9 & 99.3 & 98.3 & 34.2 & 1+1 & 2.1 \\
& SAC & 0 & 0.211 & 87.6 & 87.6 & 87.6 & 87.6 & 100.0 & 100.0 & 34.6 & 1+1 & 0.1 \\
& SAC & 1 & 0.161 & 87.6 & 87.6 & 87.6 & 87.6 & 100.0 & 100.0 & 34.4 & 1+1 & 0.2 \\
& SAC & 2 & 0.315 & 87.6 & 87.6 & 87.6 & 87.6 & 100.0 & 100.0 & 35.1 & 1+1 & 1.3 \\
& TD3 & 1 & 0.873 & 87.6 & 87.6 & 87.6 & 87.6 & 100.0 & 100.0 & 35.1 & 1+1 & 2.0 \\
& TD3 & 2 & 0.048 & 87.6 & 87.6 & 87.6 & 87.6 & 100.0 & 100.0 & 35.8 & 1+1 & 0.2 \\
\midrule
$\SO(3)$ & PPO & 0 & 0.139 & 75.4 & 75.4 & 75.4 & 75.4 & 100.0 & 100.0 & 44.4 & 1+1 & 0.3 \\
& PPO & 1 & 0.088 & 75.4 & 75.4 & 75.4 & 75.4 & 100.0 & 100.0 & 42.6 & 1+1 & 0.3 \\
& PPO & 2 & 0.088 & 75.4 & 75.4 & 75.4 & 75.4 & 100.0 & 100.0 & 41.8 & 1+1 & 0.5 \\
& TRPO & 0 & 0.085 & 75.4 & 75.4 & 75.4 & 75.4 & 100.0 & 100.0 & 39.2 & 1+1 & 0.5 \\
& TRPO & 1 & 0.094 & 75.4 & 75.4 & 75.4 & 75.4 & 100.0 & 100.0 & 44.5 & 1+1 & 0.3 \\
& TRPO & 2 & 0.110 & 75.4 & 75.4 & 75.4 & 75.4 & 100.0 & 100.0 & 38.9 & 1+1 & 0.2 \\
& SAC & 0 & 0.123 & 75.4 & 75.4 & 75.4 & 75.4 & 100.0 & 100.0 & 35.6 & 1+1 & 0.2 \\
& SAC & 1 & 0.081 & 75.4 & 75.4 & 75.4 & 75.4 & 100.0 & 100.0 & 35.0 & 1+1 & 0.3 \\
& SAC & 2 & 0.079 & 75.4 & 75.4 & 75.4 & 75.4 & 100.0 & 100.0 & 35.9 & 1+1 & 0.2 \\
& TD3 & 1 & 0.868 & 75.4 & 75.4 & 0.0 & 0.0 & 100.0 & 100.0 & 45.0 & 1+1 & 1.2 \\
\midrule
$\SE(2)$ & PPO & 0 & 0.196 & 82.0 & 82.0 & 82.0 & 82.0 & 100.0 & 100.0 & 77.6 & 1+1 & 0.9 \\
& PPO & 1 & 0.081 & 82.0 & 82.0 & 60.0 & 60.4 & 100.0 & 99.5 & 109.0 & 4+1 & 3.6 \\
& PPO & 2 & 0.141 & 82.0 & 82.0 & 82.0 & 82.0 & 100.0 & 100.0 & 47.7 & 1+1 & 0.8 \\
& TRPO & 0 & 0.110 & 82.0 & 82.0 & 82.0 & 82.0 & 100.0 & 100.0 & 32.5 & 1+1 & 0.5 \\
& TRPO & 1 & 0.085 & 82.0 & 82.0 & 76.0 & 75.6 & 100.0 & 97.3 & 37.4 & 2+2 & 3.7 \\
& TRPO & 2 & 0.169 & 82.0 & 82.0 & 82.0 & 82.0 & 100.0 & 100.0 & 43.1 & 1+1 & 0.6 \\
& SAC & 0 & 0.394 & 82.0 & 82.0 & 82.0 & 82.0 & 100.0 & 100.0 & 32.2 & 1+1 & 0.4 \\
& SAC & 1 & 0.335 & 82.0 & 82.0 & 82.0 & 82.0 & 100.0 & 100.0 & 32.1 & 1+1 & 0.4 \\
& SAC & 2 & 0.374 & 82.0 & 82.0 & 82.0 & 82.0 & 100.0 & 100.0 & 31.1 & 1+1 & 0.4 \\
& TD3 & 0 & 0.166 & 82.0 & 82.0 & 82.0 & 82.0 & 100.0 & 100.0 & 35.1 & 1+1 & 0.7 \\
& TD3 & 1 & 0.183 & 82.0 & 82.0 & 82.0 & 82.0 & 100.0 & 100.0 & 38.6 & 1+1 & 0.8 \\
& TD3 & 2 & 0.089 & 82.0 & 82.0 & 82.0 & 82.0 & 100.0 & 100.0 & 34.7 & 1+1 & 1.0 \\
\midrule
$\SE(3)$ & PPO & 0 & 0.132 & 65.1 & 65.0 & 45.1 & 46.3 & 99.8 & 99.8 & 44.1 & 5+1 & 2.6 \\
& PPO & 1 & 0.181 & 65.1 & 52.5 & 32.7 & 32.2 & 80.6 & 67.2 & 27.7 & 7+1$^{\dagger}$ & 3.9 \\
& PPO & 2 & 0.102 & 65.1 & 65.1 & 40.9 & 40.8 & 100.0 & 99.5 & 43.9 & 6+1 & 2.4 \\
& TRPO & 0 & 0.394 & 65.1 & 65.1 & 65.1 & 65.1 & 100.0 & 100.0 & 46.4 & 1+1 & 1.3 \\
& TRPO & 1 & 0.426 & 65.1 & 65.1 & 64.3 & 62.3 & 100.0 & 100.0 & 43.1 & 2+1 & 1.1 \\
& TRPO & 2 & 0.485 & 65.1 & 65.1 & 65.1 & 65.1 & 100.0 & 100.0 & 48.0 & 1+1 & 1.1 \\
\bottomrule
\end{tabular}~
\end{table}

\begin{table}[p]
\centering
\caption{The 55 closed loops with an actuator error of Table~\ref{tab:pert}, with $\spr(A)$ after the actuator error. The regions are in percent of $\mu_0$, and the four ratios in percent of the oracle region. The column $V_h$ is the neural Lyapunov baseline of Section~\ref{sec:learning}, and an asterisk marks a closed loop on which it gives no region. Differences quoted in the text are computed before rounding. The level $\rho$ and the rounds of the two runs of the search are those of $V$. A superscript $c$ marks a candidate kept from the closed loop without the actuator error with a new level, and $w$ one refitted from it (Appendix~\ref{app:training}). The others are trained from scratch. For a candidate marked $c$, the fit is measured on a closed loop on which it was not trained. A dagger marks a closed loop whose second run of the search starts from a lower level, as explained in the text.}
\label{tab:pert-all}
\scriptsize\setlength{\tabcolsep}{1.4pt}\renewcommand{\arraystretch}{0.87}
\begin{tabular}{@{}lllrrrrrrrrrrrrr@{}}
\toprule
group & policy & actuator error & $\spr(A)$ & oracle & $V$ & $V_d$ & $V_P$ & $V_h$ & $V$/oracle & $V_d$/oracle & $V_P$/oracle & $V_h$/oracle & $\rho$ & rounds & fit \\
\midrule
$\SE(2)$ & PPO 1 & gain 0.5 & 0.52 & 82.0 & 82.0$^{w}$ & 63.7 & 66.0 & 79.6 & 100.0 & 77.8 & 80.5 & 97.1 & 129.6 & 3+1 & 6.7 \\
& PPO 1 & gain 0.75 & 0.28 & 82.0 & 82.0$^{c}$ & 62.0 & 63.5 & 79.5 & 100.0 & 75.7 & 77.4 & 97.0 & 109.0 & 3+1 & 20.8 \\
& PPO 1 & gain 1.25 & 0.35 & 82.0 & 82.0 & 57.9 & 57.0 & 81.2 & 100.0 & 70.7 & 69.6 & 99.0 & 88.0 & 4+1 & 4.4 \\
& TRPO 1 & gain 1.25 & 0.36 & 82.0 & 82.0 & 75.7 & 72.9 & 66.5 & 100.0 & 92.3 & 88.9 & 81.1 & 36.9 & 3+1 & 1.9 \\
& PPO 1 & gain 1.5 & 0.62 & 82.0 & 82.0 & 57.5 & 53.2 & 81.5 & 100.0 & 70.2 & 64.9 & 99.4 & 106.4 & 1+1 & 3.7 \\
& TRPO 1 & gain 1.5 & 0.63 & 82.0 & 82.0 & 75.3 & 67.1 & 64.7 & 100.0 & 91.8 & 81.9 & 78.9 & 36.6 & 3+1 & 2.1 \\
& PPO 1 & per-axis gains 1 & 0.31 & 82.0 & 82.0$^{w}$ & 55.3 & 59.6 & 71.4 & 100.0 & 67.5 & 72.7 & 87.1 & 265.5 & 3+1 & 2.8 \\
& TRPO 1 & per-axis gains 1 & 0.29 & 82.0 & 82.0$^{c}$ & 75.8 & 75.9 & 64.7 & 100.0 & 92.5 & 92.6 & 78.9 & 52.6 & 3+1 & 10.2 \\
& PPO 1 & per-axis gains 2 & 0.12 & 82.0 & 81.9$^{c}$ & 57.9 & 59.4 & 74.5 & 99.9 & 70.7 & 72.4 & 90.8 & 50.4 & 4+1 & 5.6 \\
& PPO 1 & misalignment 5$^\circ$ & 0.07 & 82.0 & 82.0$^{c}$ & 63.6 & 63.9 & 82.0 & 100.0 & 77.5 & 78.0 & 100.0 & 110.1 & 3+1 & 4.2 \\
& TRPO 1 & misalignment 5$^\circ$ & 0.10 & 82.0 & 82.0$^{c}$ & 74.6 & 74.3 & 68.9 & 100.0 & 91.0 & 90.7 & 84.1 & 36.1 & 3+1 & 4.0 \\
& PPO 0 & misalignment 10$^\circ$ & 0.20 & 82.0 & 82.0$^{c}$ & 72.3 & 74.6 & 82.0 & 100.0 & 88.2 & 91.1 & 100.0 & 65.2 & 4+1 & 3.3 \\
& PPO 1 & misalignment 10$^\circ$ & 0.14 & 82.0 & 82.0$^{c}$ & 67.3 & 68.3 & 82.0 & 100.0 & 82.1 & 83.3 & 100.0 & 111.3 & 2+4 & 5.5 \\
& TRPO 1 & misalignment 10$^\circ$ & 0.17 & 82.0 & 82.0$^{c}$ & 73.2 & 73.6 & 68.2 & 100.0 & 89.3 & 89.8 & 83.2 & 35.1 & 4+2 & 4.6 \\
& PPO 0 & misalignment 15$^\circ$ & 0.28 & 82.0 & 82.0 & 62.7 & 63.6 & 72.3 & 100.0 & 76.5 & 77.6 & 88.2 & 111.5 & 3+1 & 4.0 \\
& PPO 1 & misalignment 15$^\circ$ & 0.23 & 82.0 & 82.0 & 71.5 & 73.8 & 82.0 & 100.0 & 87.3 & 90.0 & 100.0 & 129.4 & 4+1 & 3.1 \\
& TRPO 1 & misalignment 15$^\circ$ & 0.26 & 82.0 & 82.0 & 71.8 & 73.3 & 57.6 & 100.0 & 87.6 & 89.5 & 70.2 & 43.4 & 3+2 & 1.5 \\
& PPO 0 & misalignment 30$^\circ$ & 0.54 & 81.7 & 81.7 & 48.3 & 51.3 & 66.6 & 99.9 & 59.1 & 62.8 & 81.5 & 64.0 & 2+1 & 4.2 \\
& TD3 2 & misalignment 30$^\circ$ & 0.51 & 82.0 & 82.0$^{w}$ & 62.4 & 63.3 & 82.0 & 100.0 & 76.2 & 77.3 & 100.0 & 51.1 & 4+2 & 1.9 \\
& TRPO 1 & misalignment 30$^\circ$ & 0.52 & 82.0 & 82.0 & 67.2 & 76.5 & 41.4 & 100.0 & 82.0 & 93.4 & 50.6 & 36.7 & 4+2 & 2.6 \\
& PPO 0 & misalignment 40$^\circ$ & 0.71 & 81.5 & 81.3 & 40.5 & 47.2 & 61.2 & 99.7 & 49.7 & 57.9 & 75.1 & 66.4 & 2+1 & 3.0 \\
& TD3 0 & misalignment 40$^\circ$ & 0.63 & 82.0 & 82.0$^{w}$ & 31.7 & 43.3 & 76.5 & 100.0 & 38.6 & 52.9 & 93.3 & 52.5 & 4+2 & 3.0 \\
& TD3 1 & misalignment 40$^\circ$ & 0.62 & 82.0 & 80.6$^{c}$ & 54.3 & 0.6 & 81.3 & 98.3 & 66.2 & 0.7 & 99.2 & 29.1 & 3+2 & 12.1 \\
& TD3 2 & misalignment 40$^\circ$ & 0.67 & 82.0 & 82.0$^{w}$ & 43.3 & 53.1 & 80.8 & 100.0 & 52.8 & 64.8 & 98.6 & 38.9 & 2+1 & 4.1 \\
\midrule
$\SE(3)$ & PPO 1 & gain 0.5 & 0.52 & 65.1 & 51.8$^{w}$ & 34.0 & 37.2 & 8.5 & 79.5 & 52.2 & 57.0 & 13.1 & 51.0 & 5+3 & 3.7 \\
& PPO 2 & gain 0.5 & 0.53 & 65.1 & 64.6$^{c}$ & 41.4 & 42.6 & 64.1 & 99.1 & 63.5 & 65.4 & 98.5 & 35.9 & 2+5 & 45.4 \\
& PPO 1 & gain 0.75 & 0.28 & 65.1 & 56.1$^{w}$ & 33.5 & 35.2 & 4.5 & 86.1 & 51.5 & 54.0 & 6.9 & 38.1 & 5+1 & 3.6 \\
& PPO 2 & gain 0.75 & 0.30 & 65.1 & 65.1$^{c}$ & 41.1 & 41.9 & 36.1 & 100.0 & 63.1 & 64.3 & 55.4 & 42.3 & 5+2 & 22.0 \\
& PPO 1 & gain 1.25 & 0.48 & 65.1 & 47.3 & 31.6 & 26.5 & 8.7 & 72.6 & 48.5 & 40.7 & 13.3 & 25.4 & 3+1 & 3.8 \\
& PPO 2 & gain 1.25 & 0.35 & 65.1 & 65.1 & 40.8 & 39.2 & 63.1 & 100.0 & 62.6 & 60.2 & 96.9 & 44.0 & 7+2 & 3.4 \\
& PPO 1 & gain 1.5 & 0.77 & 65.1 & 58.6 & 30.3 & 16.4 & 0.0$^{*}$ & 90.0 & 46.6 & 25.2 & 0.0 & 30.6 & 3+1 & 3.6 \\
& PPO 2 & gain 1.5 & 0.62 & 65.1 & 65.1 & 40.8 & 35.1 & 59.7 & 100.0 & 62.6 & 54.0 & 91.7 & 43.8 & 6+4 & 2.4 \\
& PPO 0 & per-axis gains 1 & 0.26 & 65.1 & 64.2$^{w}$ & 32.1 & 34.0 & 22.6 & 98.6 & 49.4 & 52.1 & 34.7 & 45.0 & 3+3 & 3.1 \\
& PPO 1 & per-axis gains 1 & 0.25 & 65.1 & 51.3$^{w}$ & 22.5 & 23.1 & 7.5 & 78.8 & 34.5 & 35.4 & 11.5 & 31.2 & 3+1 & 4.0 \\
& PPO 2 & per-axis gains 1 & 0.24 & 65.1 & 65.0$^{w}$ & 38.2 & 37.8 & 53.7 & 99.7 & 58.6 & 58.0 & 82.4 & 47.0 & 3+4 & 2.3 \\
& PPO 0 & per-axis gains 2 & 0.18 & 65.1 & 65.0$^{c}$ & 42.3 & 42.3 & 47.5 & 99.7 & 64.9 & 64.9 & 73.0 & 43.5 & 3+1 & 4.6 \\
& PPO 1 & per-axis gains 2 & 0.27 & 65.1 & 48.6$^{w}$ & 31.5 & 31.5 & 14.5 & 74.6 & 48.4 & 48.3 & 22.2 & 26.7 & 6+1 & 3.7 \\
& PPO 2 & per-axis gains 2 & 0.22 & 65.1 & 65.1$^{c}$ & 40.8 & 40.7 & 62.7 & 99.9 & 62.6 & 62.4 & 96.2 & 39.4 & 5+4 & 4.5 \\
& PPO 1 & misalignment 5$^\circ$ & 0.17 & 65.1 & 50.7$^{w}$ & 30.8 & 30.2 & 18.6 & 77.8 & 47.3 & 46.4 & 28.5 & 27.3 & 3+4 & 3.6 \\
& PPO 2 & misalignment 5$^\circ$ & 0.15 & 65.1 & 65.1$^{c}$ & 38.0 & 38.1 & 39.1 & 100.0 & 58.4 & 58.5 & 60.1 & 42.9 & 6+1 & 2.6 \\
& PPO 0 & misalignment 10$^\circ$ & 0.24 & 65.1 & 64.9$^{w}$ & 38.2 & 38.0 & 36.3 & 99.7 & 58.7 & 58.4 & 55.7 & 45.7 & 3+6 & 2.5 \\
& PPO 1 & misalignment 10$^\circ$ & 0.26 & 65.1 & 43.2$^{w}$ & 28.5 & 27.9 & 17.8 & 66.3 & 43.8 & 42.8 & 27.3 & 23.8 & 2+1 & 3.7 \\
& PPO 2 & misalignment 10$^\circ$ & 0.23 & 65.1 & 64.8$^{w}$ & 35.0 & 35.1 & 46.7 & 99.5 & 53.7 & 53.9 & 71.8 & 38.6 & 4+2 & 2.7 \\
& PPO 0 & misalignment 15$^\circ$ & 0.32 & 65.1 & 64.9 & 33.2 & 33.2 & 37.5 & 99.7 & 50.9 & 51.0 & 57.5 & 44.7 & 7+1 & 2.7 \\
& PPO 1 & misalignment 15$^\circ$ & 0.34 & 65.1 & 48.6 & 26.0 & 25.4 & 15.5 & 74.6 & 40.0 & 39.0 & 23.9 & 26.7 & 5+2 & 3.8 \\
& PPO 2 & misalignment 15$^\circ$ & 0.32 & 65.1 & 63.7 & 31.6 & 31.8 & 35.4 & 97.7 & 48.5 & 48.8 & 54.4 & 36.0 & 5+1 & 3.4 \\
& PPO 0 & misalignment 30$^\circ$ & 0.58 & 65.1 & 64.3 & 16.3 & 18.6 & 0.0$^{*}$ & 98.8 & 25.0 & 28.5 & 0.0 & 44.8 & 4+3$^{\dagger}$ & 3.2 \\
& PPO 1 & misalignment 30$^\circ$ & 0.60 & 65.1 & 45.2 & 15.1 & 16.5 & 0.0$^{*}$ & 69.4 & 23.2 & 25.3 & 0.0 & 26.7 & 3+1 & 5.1 \\
& PPO 2 & misalignment 30$^\circ$ & 0.58 & 65.1 & 61.0 & 20.9 & 20.4 & 0.0$^{*}$ & 93.6 & 32.1 & 31.3 & 0.0 & 36.6 & 4+1 & 4.1 \\
& PPO 0 & misalignment 40$^\circ$ & 0.75 & 63.5 & 59.5 & 7.9 & 5.5 & 0.0$^{*}$ & 93.7 & 12.4 & 8.7 & 0.0 & 40.8 & 3+2 & 3.5 \\
& PPO 1 & misalignment 40$^\circ$ & 0.76 & 64.9 & 44.8 & 9.7 & 7.2 & 0.0$^{*}$ & 69.0 & 14.9 & 11.1 & 0.0 & 29.8 & 3+1 & 4.9 \\
& PPO 2 & misalignment 40$^\circ$ & 0.74 & 65.0 & 46.7 & 13.0 & 9.9 & 0.0$^{*}$ & 71.8 & 20.0 & 15.3 & 0.0 & 30.2 & 5+2 & 3.6 \\
& TRPO 0 & misalignment 40$^\circ$ & 0.67 & 65.1 & 65.1 & 45.9 & 38.2 & 46.5 & 99.9 & 70.5 & 58.6 & 71.4 & 60.6 & 2+2 & 1.8 \\
& TRPO 1 & misalignment 40$^\circ$ & 0.70 & 65.1 & 64.3 & 28.3 & 20.8 & 0.0$^{*}$ & 98.7 & 43.5 & 32.0 & 0.0 & 56.7 & 13+2 & 2.0 \\
& TRPO 2 & misalignment 40$^\circ$ & 0.68 & 65.1 & 65.1 & 37.9 & 25.3 & 40.3 & 100.0 & 58.1 & 38.8 & 61.9 & 72.9 & 1+2 & 1.3 \\
\bottomrule
\end{tabular}~
\end{table}

\paragraph{The policy correction.} Over $10^{5}$ points drawn uniformly in the ball, subtracting $\mu(\vecop\Id)$ changes the applied increment $\proj(\amax\sat(\cdot))$ by at most $0.17\,\amax$ over the 39 policies, with the largest change for PPO seed 0 on $\SO(3)$, and its median over the ball is below $0.14\,\amax$ on every policy and below $0.02\,\amax$ on 24 of them. On the 39 policies, the norm $\|\mu(\vecop\Id)\|$ of the subtracted constant is at most $0.17$.

\paragraph{The stability test.}\looseness -1 With central-difference steps $10^{-4}$, $10^{-5}$, $10^{-6}$ (the step used) and $10^{-7}$, the estimate of $\spr(A)$ changes by at most $1.4\times10^{-4}$ over the 39 policies, on SAC seed 0 on $\SO(2)$, and by less than $10^{-6}$ on 36 of them. The largest value among the 94 evaluated closed loops is $0.873$.

\paragraph{Cost.} On two CPU threads per closed loop, training and refinement take 0.5 to 1.3 hours on $\SO(2)$, 0.8 to 1.3 on $\SO(3)$, 0.5 to 3.7 on $\SE(2)$ and 0.8 to 5.0 on $\SE(3)$. The evaluation, including the search and all checks for the candidate, the two quadratic baselines and the critic function, takes 7 to 22 minutes on $\SO(2)$, 10 to 28 on $\SO(3)$, 19 to 73 on $\SE(2)$ and 22 to 145 on $\SE(3)$. These times are for the 39 policies and the 8 closed loops with an actuator error that ran first. The other 18 closed loops with an actuator error ran four at a time on the same machine, and their training and refinement took 4.0 to 6.7 hours. Over the 26 closed loops with an actuator error of the first screen, training and refinement take 2.9 to 6.7 hours, and Appendix~\ref{app:training} gives the cost for the other 29. Simulating the oracle, 3000 steps from each of its $2\times10^{5}$ level samples and from the initial errors that lie in the ball, up to $4\times10^{5}$ points in all, takes 18 to 48 minutes per closed loop. For the first 47, on two threads, evaluating the candidate takes 1.7 to 3.4 microseconds per point, and computing $\Vm$ by simulation takes 0.5 to 1.0 milliseconds per point with 600 steps and 2.6 to 4.9 with 3000 steps.

\paragraph{The neural Lyapunov baseline.} On the 55 closed loops with an actuator error, the neural Lyapunov baseline $V_h$ of Section~\ref{sec:learning} is trained from scratch by Algorithm~\ref{alg:train-full} with the weight of the residual set to zero and without the refinement, so with the iterations, constants and counterexample search of the candidate. Its level comes from the level rule, with both runs of the search, as for $V$. The second search lowered its level on 12 closed loops. The candidate exceeds $V_h$ by more than $0.01$ percentage points on 49 closed loops, by at most $64.3$ points, the two agree within $0.01$ on 5, and $V_h$ is larger on one, TD3 seed 1 on $\SE(2)$ with misalignment $40^\circ$, by $0.8$ points. Over all 55 closed loops, the median of the difference is $18.1$ points. Against the larger quadratic baseline, $V_h$ is larger on 29 closed loops and smaller on 26, and it recovers a median of $73.0$ percent of the oracle region. On 8 closed loops on $\SE(3)$, PPO seeds 0 to 2 at $30^\circ$ and $40^\circ$, PPO seed 1 with gain $1.5$ and TRPO seed 1 at $40^\circ$, the search leaves a level between $0.015$ and $7.8$, under which no sampled initial error lies, and the entry condition fails there at 12 to 72 sampled points of $N_b$, so $V_h$ gives no region. On all 55, every simulated point of $N_b$ converges and no region of $V_h$ contains a sampled initial error that fails to converge.

\paragraph{Search of the oracle region.} On PPO seed 0 on $\SE(3)$ with misalignment $40^\circ$, PPO seed 1 on $\SE(3)$ without an actuator error and PPO seed 0 on $\SE(2)$ with misalignment $40^\circ$, we searched the oracle region for points whose next step leaves the ball, with $\Wm$ from 3000 simulated steps. Since $\Wm$ decreases exactly along trajectories, decrease and entry hold for it, and every point of the oracle region converges, so the search looks only for points at which containment fails. It drew $2\times10^{7}$ fresh points of the ball, and then moved the counterexamples with the least $\Wm$ by random steps, kept only when the next step still leaves the ball and $\Wm$ falls. On all three closed loops it found counterexamples at which $\Wm$ is below the level of the oracle region, with least values of $0.885$, $0.559$ and $0.932$ times that level, in this order. Most of the reduction comes from the $2\times10^{7}$ points, 100 times the $2\times10^{5}$ samples that set the level of the oracle region, and not from the random steps. With the level lowered below the least value of $\Wm$ at a counterexample, the oracle region falls from $63.47$ to $60.53$ percent of $\mu_0$ on the first closed loop, and the candidate recovers $98.2$ percent of it instead of $93.7$. It falls by $0.06$ points on the third and not at all on the second. The oracle region is therefore optimistic on these closed loops, and the fractions of it that the tables report would not fall if its level were lowered. The search took 21 to 32 minutes per closed loop on one thread. Stored in full, the trajectory needs about 9 megabytes per point for a gradient of $\Wm$ through 3000 steps, so we recompute it during the backward pass. Once a trajectory reaches the goal to machine precision, automatic differentiation of the norm of the increment and of the formula for $\Exp$ gives $0/0$, and we hold the state fixed from then on, which leaves $\Wm$ unchanged. A gradient then takes about 20 milliseconds per point on one thread, against $0.024$ milliseconds for a gradient of $\delta V$. One round of the third stage of the level rule takes about 26 minutes on the candidate. On the oracle only the ascent on $\|F(z)\| - \gamma$ is needed, but its penalty for leaving $\Omega_\rho$ needs a gradient of $\Wm$ at every step, so one round would take about 8 days.

\FloatBarrier
\section{Held-out checks, ablations and the region without the actuator error}
\label{app:ablations}
\paragraph{Held-out samples.} For each of the 94 closed loops, $10^{6}$ further points were drawn uniformly in the ball, with a random seed used nowhere else, and the points inside the region of the candidate at the level of the first run of the search, which contains the reported region, were checked for decrease outside $N_b$ and for containment. On 93 of the closed loops no point fails. Suppose that the points of the region at which a condition fails make up a fraction $p$ of the ball. Then all $10^{6}$ points miss them with probability $(1 - p)^{10^{6}} \le e^{-10^{6}p}$, so with confidence $1 - 10^{-6}$ they make up at most $\ln(10^{6})/10^{6} = 1.4\times10^{-5}$ of the ball on each of the 93 closed loops. On TRPO seed 2 on $\SO(2)$ one point fails containment. A grid of $2\times10^{7}$ points with spacing $2.7\times10^{-7}$ on the ball of $\SO(2)$ finds five such points, all within $1.1\times10^{-6}$ of the rim $\|z\| = \gamma$. On this grid, multiplying the level by $1 - 1.1\times10^{-6}$ removes them and leaves the reported size of the region unchanged. The same grid finds no failing point on the other ten closed loops on $\SO(2)$.

\paragraph{The quadratic term alone.} The function $z^{\top}\Pm z$, which is the second-order term of the accumulated cost and the quadratic term of the candidate \eqref{eq:cert}, was put through the same level rule on 16 closed loops: the 8 policies of Section~\ref{sec:sims} on which the candidate exceeds both baselines and 8 closed loops with an actuator error. On 14 of them its region is smaller than that of the candidate by 27.9 to 75.4 percentage points. On the other two, TRPO seed 1 on $\SE(3)$ and TRPO seed 2 on $\SO(2)$, it is within 0.3 points of it, and there every function gives nearly the whole ball. On $\SO(n)$, where $Q = 2\Id$, it is twice the linearization-based baseline and has the same region. The region away from the goal therefore comes from the network part of the candidate.

\paragraph{Ablations.} Table~\ref{tab:ablation} changes one part of the method at a time on five of these closed loops: four policies without an actuator error and PPO seed 1 on $\SE(3)$ with gain $1.5$. The variant without the quadratic term is the network part $V_0$ of \eqref{eq:cert} alone. The hinge-only variant drops the residual $\calL_{\mathrm{rec}}$ from \eqref{eq:objective}, which leaves the usual loss of a neural Lyapunov function, the loss of the neural Lyapunov baseline of Section~\ref{sec:learning}, here followed by the refinement. The regression variant replaces the residual by the mean squared difference between $\log(1 + V)$ and $\log(1 + \Wm)$ at 16384 fixed points, with $\Wm$ from 600 simulated steps. The variant without refinement is the candidate after training, before the refinement rounds. Every other part and constant is as in Table~\ref{tab:consts}, and each variant goes through the same level rule and, where it has one, the same refinement, whose rounds train with the loss of the candidate, residual included. The hinge-only and regression variants keep their candidate from before the refinement, except the regression variant on the two closed loops of PPO seed 1 on $\SE(3)$ without an actuator error and with gain $1.5$. On TD3 seed 1 on $\SO(3)$ and PPO seed 1 on $\SE(2)$ every variant gives the whole ball. On the other three closed loops, all on $\SE(3)$, the hinge-only and regression variants each lose more than 50 percentage points on at least one of them, up to 58.6 and 58.0 points, and the variant without refinement loses 3.1 to 15.0 points. With gain $1.5$, the region of the hinge-only variant contains none of the samples of $\mu_0$ and fails the entry check. The variant without the quadratic term loses 0.3 points on PPO seed 2 and 9.3 points with gain $1.5$, and gains 1.1 points on PPO seed 1. There the level of the candidate comes from the run of the search with $r_b = 0.2$ (Appendix~\ref{app:results}), which the variants did not get. The entry marked with an asterisk comes from a refinement that produced non-finite losses from its fourth round on, whose steps were skipped. The kept candidate is from before the first round, with the region of the table, since the rounds in between gave smaller regions. The final candidates of Tables~\ref{tab:all} and \ref{tab:pert-all} met no non-finite loss, except PPO seed 1 on $\SE(3)$ with the first per-axis gains (Appendix~\ref{app:training}).

\begin{table}[!ht]
\centering
\caption{Regions in percent when one part of the method is changed, on five closed loops. The columns are explained in the text, and an asterisk marks a refinement with non-finite losses. Differences quoted in the text are computed before rounding.}
\label{tab:ablation}
\scriptsize\setlength{\tabcolsep}{3.5pt}
\begin{tabular}{@{}lllrrrrrr@{}}
\toprule
group & policy & actuator error & candidate & no quadratic term & hinge-only & regression & no refinement & oracle \\
\midrule
$\SO(3)$ & TD3 1 & none & 75.4 & 75.4 & 75.4 & 75.4 & 75.4 & 75.4 \\
$\SE(2)$ & PPO 1 & none & 82.0 & 82.0 & 82.0 & 82.0 & 82.0 & 82.0 \\
$\SE(3)$ & PPO 1 & none & 52.5 & 53.6 & 9.9 & 8.5 & 49.4 & 65.1 \\
$\SE(3)$ & PPO 2 & none & 65.1 & 64.9 & 63.8 & 7.1$^{*}$ & 51.2 & 65.1 \\
$\SE(3)$ & PPO 1 & gain 1.5 & 58.6 & 49.3 & 0.0 & 8.0 & 43.6 & 65.1 \\
\bottomrule
\end{tabular}~
\end{table}

\paragraph{Effect of the refinement.} Over the 16 closed loops, the refinement changes the region by more than one percentage point on 4, all PPO on $\SE(3)$, where it enlarges it by 2.5 to 15.0 points. The largest decrease it causes is 0.2 points.

\paragraph{Sensitivity to $r_a$ and $r_b$.} Here each region comes from one run of the search, at every value compared. On six closed loops the candidate was put through the level rule with $r_b = 0.05$ and $0.2$ as well as $0.1$, and all checks pass. On TD3 seed 1 on $\SO(3)$, PPO seed 1 on $\SE(2)$ and PPO seed 2 on $\SE(3)$ the region is the same at all three values to within $0.003$ percentage points, and on PPO seed 0 on $\SE(3)$ with misalignment $30^\circ$ and PPO seed 1 with gain $1.5$ it changes by at most $0.23$ points. On PPO seed 1 on $\SE(3)$ without an actuator error it is $62.0$, $60.5$ and $52.5$ percent at $r_b = 0.05$, $0.1$ and $0.2$. There the run with $r_b = 0.2$ found two points at the edge of the ball, below the level of the run with $r_b = 0.1$, at which containment fails. Containment does not depend on $r_b$, so the reported level of this closed loop is set below them (Appendix~\ref{app:results}). The starting points of the search change with $r_b$, and its result depends on them. Trained from scratch with $r_a = 0.15$ instead of $0.3$, the candidate gives $64.2$, $59.9$ and $58.2$ percent on PPO seed 0 on $\SE(3)$ with misalignment $30^\circ$, PPO seed 1 without an actuator error and PPO seed 1 with gain $1.25$, against $64.6$, $60.5$ and $47.3$ percent with $r_a = 0.3$.

\paragraph{Candidate seeds.} A second candidate, trained with another seed, was evaluated on the same 16 closed loops, and both candidates are compared here with one run of the search. The region of the second candidate minus that of the first has a median of 0.00 and a mean of 0.22 percentage points. The two differ by more than one point only on PPO seed 1 on $\SE(3)$: by $+3.53$ points without an actuator error, $+1.97$ with gain $1.25$ and $-1.15$ with gain $1.5$. On the other 13 the largest difference is 0.69 points. With either seed, on the two closed loops of PPO seed 1 on $\SE(3)$ with gain $1.25$ and $1.5$, the quick estimate inside the refinement exceeds the region given by the level rule, as Appendix~\ref{app:training} shows for gain $1.25$.

\paragraph{The region without the actuator error.} The screens of Appendix~\ref{app:screen} leave a gap of at least 5 points on 55 closed loops. For each of them, the region of the candidate trained for the same policy without the actuator error was checked under the closed loop with the actuator error, at the points among $10^{6}$ drawn uniformly in the ball that lie in that region. On 21 of the 32 closed loops with a misalignment, decrease or containment fails at some of these points, on up to 2.2 percent of those outside $N_b$. It fails on none of the 4 closed loops at $5^\circ$, on 1 of the 6 at $10^\circ$, on 4 of the 6 at $15^\circ$, and on all 16 at $30^\circ$ and $40^\circ$. A region computed without a misalignment therefore does not in general satisfy the Lyapunov conditions once the misalignment is present, and the region has to be recomputed. On all 55 closed loops every sampled initial error in such a region still converges in 600 steps, so the region still lies in the region of attraction as far as simulation shows, but the Lyapunov conditions no longer hold on it. With a constant gain a condition fails on 1 of the 14 closed loops, at one point, and with per-axis gains on 5 of the 9. The samples miss failures that the search finds. On 25 of the 29 closed loops of the later screens, the search of the level rule finds points of this region at which a condition fails, against 11 for the samples.

\FloatBarrier
\section{Proof of the regions on \texorpdfstring{$\SO(2)$}{SO(2)}}
\label{app:so2proof}
On $\SO(2)$ the chart is one-dimensional, so the conditions of Section~\ref{sec:learning} can be proved on the whole region by interval arithmetic, with no samples and no simulation. The ball $[-\gamma, \gamma]$ is covered by $2^{14}$ cells. On each cell, interval arithmetic through the policy network, the network of the candidate, the quadratic term and the trigonometric maps bounds the values and derivatives of $F$ and $V$. The constant $\mu(\vecop\Id)$ that the correction subtracts is itself enclosed in an interval, so $F(0) = 0$ holds exactly. Every cell that meets the region must satisfy containment, $F(\text{cell}) \subseteq [-\gamma, \gamma]$. If it also meets the complement of $N_b$, it must satisfy decrease, $\sup \delta V < 0$, and if it meets $N_b$, entry, $\sup V(F(\text{cell})) \le \rho$. A cell that fails is bisected, up to 16 times. If cells still fail, the level is lowered just below the smallest value of $V$ on them and the check is repeated.

Near the goal, $N_b$ is the interval $|z| < \theta_b$ with $\theta_b = \arccos(1 - r_b^{2}/4)$, and every point of it must be shown to converge. On eight closed loops $|F'| < 1$ on $N_b$, so $F$ is a contraction there. On the three SAC closed loops the bound on $|F'|$ exceeds 1 inside $N_b$. There, $|F(z)| < |z|$ is proved instead for $0 < |z| < \theta_b$: on a small interval around $0$ from $|F'| < 1$ and $F(0) = 0$, and elsewhere by bounds on cells. Then $|F^{k}(z)|$ decreases while $F^{k}(z) \neq 0$, and at a limit point $\bar z$ of the orbit $|F(\bar z)| = |\bar z|$ by continuity, so the orbit converges to $0$.

\begin{figure}[!ht]
\centering
\includegraphics[width=\linewidth]{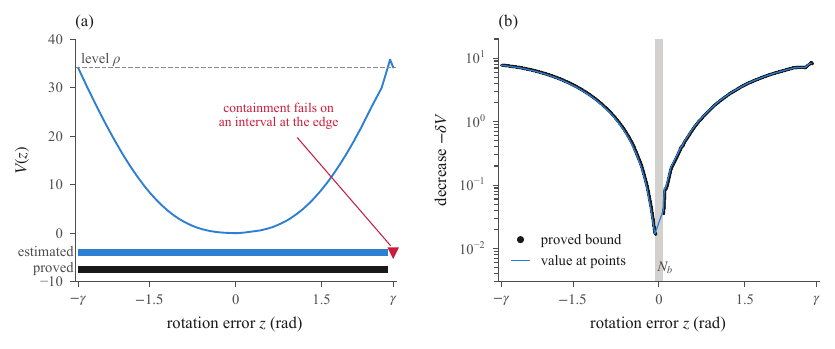}
\vspace{-8pt}
\caption{The closed loop of TRPO seed 2 on $\SO(2)$. (a) The candidate $V$ and the level $\rho$ from the level rule. The bars show the region estimated by the level rule and the region proved by interval arithmetic. They differ only on an interval of width $1.4\times10^{-6}$ at the edge of the ball, where containment fails. The candidate exceeds $\rho$ on $(2.66, 2.75)$, so neither region covers that interval. (b) The proved lower bound on the decrease $-\delta V$ on each cell of the region outside $N_b$ (gray), and its value at points.}
\label{fig:so2proof}
\end{figure}

The conditions hold at the reported level on 10 of the 11 closed loops. On TRPO seed 2, containment fails on an interval of width $1.4\times10^{-6}$ at the edge of the ball, where the held-out check also found its one failing point (Appendix~\ref{app:ablations}). The proved level there is lower than the reported one by a relative $1.2\times10^{-6}$. The proved regions are $87.50$ percent of $\mu_0$ on the ten closed loops, which is the whole ball, and $86.03$ percent on TRPO seed 2. The sizes measured on $2\times10^{5}$ samples of $\mu_0$ are $87.56$ and $86.14$ percent. The smallest proved decrease outside $N_b$ is between $0.0099$ and $0.0170$. Figure~\ref{fig:so2proof} shows both regions on TRPO seed 2. The bounds are computed in double precision with the rounding error of every operation bounded and added to the radii, not with a verified interval library with directed rounding.

\FloatBarrier
\section{The screens for the actuator errors}
\label{app:screen}
Three screens, of 300 closed loops in all, computed for each closed loop the stability test, the part of the ball that converges by simulation, and the regions of $V_d$ and $V_P$ from the first two stages of the level rule. They use no candidate. The first, of 102 closed loops, takes the PPO and TRPO policies, three seeds each, on all four groups, with gain $1.25$ and $1.5$ and with misalignment $15^\circ$, $30^\circ$ and $40^\circ$. The second, of 132, takes the same policies with gain $0.5$ and $0.75$, two draws of per-axis gains, and misalignment $5^\circ$ and $10^\circ$. The third, of 66, takes the SAC and TD3 policies on $\SE(2)$ with all eleven actuator errors. On $\SO(2)$ a misalignment is undefined. A misalignment rotates the increment about the axis $(1, 1, 1)/\sqrt3$ on $\SO(3)$, both its angular and its linear part about this axis on $\SE(3)$, and its translation part in the plane on $\SE(2)$. Every closed loop of the screens passes the stability test, with $\spr(A)$ below $0.80$.

\paragraph{The stability test under an actuator error.} The Jacobian with the actuator error follows from the one without it. An actuator error multiplies the increment by the matrix $E$ of \eqref{eq:closedloopE}, the misalignment followed by the gains, and near the goal both projections onto $\bar B_{\amax}$, before and after the error, are inactive. So $A = \Id + E(A_0 - \Id)$, with $A_0$ the Jacobian of the closed loop without the actuator error. On all 300 closed loops, $\spr$ of this matrix agrees to within $3\times10^{-14}$ with $\spr(A)$ from central differences on the closed loop with the actuator error. On each of the 30 policies of the screens, $\spr(A) < 1$ for every constant gain $k_g$ with $0 < k_g < k_{\max}$, where $k_{\max}$ depends on the policy and lies between $1.67$ and $2.37$. Except on $\SO(2)$, the same holds, in either sense of rotation, for every misalignment below an angle that depends on the policy and exceeds $40^\circ$ on every policy. Every screened constant gain and misalignment therefore lies inside the range on which $\spr(A) < 1$, and the per-axis gains give $\spr(A)$ below $0.61$.

\paragraph{The selected closed loops.} The screens compare the part of the ball that converges with the larger of the two quadratic baseline regions. On $\SO(2)$ and $\SO(3)$ the gap is below 1.7 percentage points on every closed loop. On $\SE(2)$ and $\SE(3)$, 55 closed loops leave a gap of at least 5 points, 26, 25 and 4 in the three screens, and the other 143 a smaller gap. All 55 are evaluated, in Tables~\ref{tab:pert} and \ref{tab:pert-all}. Most of them come from policies on which the baselines already fall short without an actuator error. Of the 55, 44 belong to five policies: PPO seed 1 and TRPO seed 1 on $\SE(2)$ and the three PPO seeds on $\SE(3)$. These are the policies on $\SE(2)$ and $\SE(3)$ on which, without an actuator error, the region of the candidate is larger than both quadratic baseline regions by at least $6.0$ points (Section~\ref{sec:sims}). The other 11 have a misalignment of $10^\circ$ to $40^\circ$. The screens rotate the increment in one sense only. With the samples of the screens and the opposite sense, the larger quadratic baseline region misses at least 5 points of the ball on 22 of the 90 screened closed loops with a misalignment on $\SE(2)$ and $\SE(3)$, against 32 in the sense of the screens. In both senses 18 closed loops reach 5 points, and 4 reach it only in the opposite sense. Which closed loops the screens select therefore depends on the sense. The closed loops with the opposite sense were neither simulated nor evaluated.

\FloatBarrier
\section{Deployment details}
\label{app:deploy}
\begin{figure}[!ht]
\centering
\includegraphics[width=0.82\linewidth]{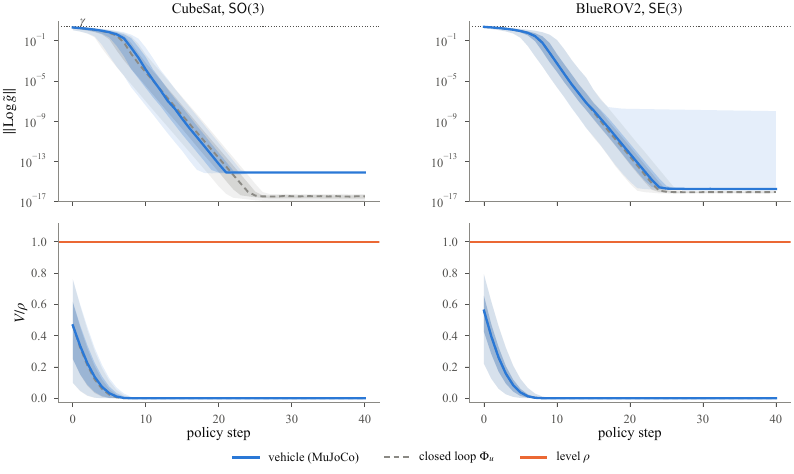}
\vspace{-8pt}
\caption{The two vehicles in MuJoCo (blue) and the closed loop $\Phi_u$ from the same initial errors (dashed gray): medians over 200 episodes, with 25 to 75 and 5 to 95 percentile bands. Top: the error norm $\|\Log\tilde g\|$, with the radius $\gamma$ of the ball (dotted). Bottom: the candidate $V$ divided by its level $\rho$, with the level itself in orange. A few BlueROV2 episodes are still converging at step 40.}
\label{fig:deploy}
\end{figure}

The closed loop $\Phi_u$ describes kinematics only: it has no velocity, inertia or torque. A rigid body fixed at a point is the mechanical system of \citet[Example 5.48]{bullo2005geometric}, and a free rigid body also translates. With inertia matrix $\mathbb{I}$, body angular velocity $\omega_b$ and torque $M_b$, its attitude obeys $\dot R = R\hat\omega_b$ and $\mathbb{I}\dot\omega_b = \mathbb{I}\omega_b\times\omega_b + M_b$. With mass $m_b$, body linear velocity $v_b$ and force $f$, its translation obeys $m_b\dot v_b = m_b\,v_b\times\omega_b + f$. The policy and the vehicle run on two time scales. At each policy step $k$, the policy commands an increment $\tau_k$ from the current pose $g_k$, which defines the target pose $g_d = g_k\Exp(\tau_k)$. A proportional-derivative controller, similar to those of \citet[Section 11.4]{bullo2005geometric}, is tuned to drive the vehicle to $g_d$ and bring it to rest within one sampling period $\Delta t$. If it succeeds, then the error at the next policy step is close to $\Phi_u$ applied to the current error. In the nominal setting, initial errors are drawn uniformly in $\Omega_\rho$ in exponential coordinates, with the vehicle at rest and without disturbances, and the same initial errors are also simulated on $\Phi_u$ for comparison (Figure~\ref{fig:deploy}). The levels are $\rho = 44.40$ for PPO seed 0 on $\SO(3)$ and $\rho = 44.08$ for PPO seed 0 on $\SE(3)$. The largest one-step increase of $V$ in Table~\ref{tab:deploy} is computed from the values of $V$ recorded at every policy step.

\paragraph{A 3U CubeSat.} The first vehicle is a free rigid body of mass 4.0 kg in the $10\times10\times34$ cm envelope of the CubeSat specification \citep{cubesat2022spec}. Its inertia is $\operatorname{diag}(0.049, 0.049, 0.008)$ kg m$^2$, an assumed value 17 to 20 percent above that of a uniform box of this size and mass. Three reaction wheels give torques of at most 1 mN m about three orthogonal axes, as a flight-qualified wheel of this class does \citep{rw001datasheet}. With $R$ the current attitude and $R_d$ the rotation part of the target, the attitude controller is $M_b = -K_P\Log(R_d^{\top}R) - K_D\omega_b + \omega_b\times\mathbb{I}\omega_b$, with gain matrices $K_P$ and $K_D$, saturated at the torque limit. Its last term cancels the gyroscopic torque. The gains are set on each axis, with critical damping, from the rest-to-rest rotation time $2(\amax\mathbb{I}_i/M_{\max})^{1/2}$, where $\mathbb{I}_i$ is the moment of inertia about axis $i$ and $M_{\max}$ the torque limit. This time is 7.85 s on the slowest axis, and the sampling period is twice that, 15.7 s.

\paragraph{A BlueROV2 Heavy.}\looseness -1 The second vehicle is an underwater vehicle with eight thrusters, fully actuated in six degrees of freedom, with mass 13.5 kg and displaced volume 0.0134 m$^3$. Its hydrodynamic constants come from the experimentally validated model of \citet{vonbenzon2022bluerov2}. The model has linear and quadratic damping and a center of buoyancy 1 cm above the center of mass, and each thruster gives at most 30 N. The rotational added mass is added to the inertia, which becomes $0.449$, $0.365$ and $0.592$ kg m$^2$. The translational added mass is $6.36$, $7.12$ and $18.68$ kg, with its Coriolis terms. The translational added mass differs between the axes and is comparable to the body mass, so the scalar mass of a MuJoCo body cannot represent it. At every integration step it is applied as the force that makes the acceleration solve the equation of \citet{vonbenzon2022bluerov2} with the full mass matrix. This force depends only on the current state. Under a constant body force $f$, the simulated acceleration along each axis equals $f/(m_b + m_a)$, with $m_a$ the added mass of that axis, up to the integration error. Thruster dynamics, currents and the tether are not modeled. A thruster was saturated in $0.6$ percent of the integration steps. One unit of translation in exponential coordinates is one meter. The controller tracks the target pose. Its torque is the attitude law above with gains $600$ and $80$ and without the gyroscopic term. With $x$ the current position and $x_d$ the position of the target, its force is $f = K_pR^{\top}(x_d - x) - K_dv_b + R^{\top}f_0$, with gains $400$ and $160$ and $f_0$ the upward force that cancels the weight minus the buoyancy. The force and torque are distributed to the thrusters, each saturated at 30 N. The gains and the sampling period of 2 s were chosen by a sweep so that an increment of length $\amax$, the clip length, is carried out within one period.

\paragraph{Initial rates and a halved period.} Two further settings were run on both vehicles, with 200 episodes each. In the first, the initial body rate is drawn uniformly in a ball whose radius is a quarter or a half of the largest rate reached in the nominal episodes. In the second, the sampling period is halved. Then only 73 percent of an increment is carried out within a period on the CubeSat and 96 percent on the BlueROV2. No episode leaves $\Omega_\rho$ in either setting, and $V$ increases by at most $2\times10^{-6}$ at a policy step.

\paragraph{The BlueROV2 with an actuator error.}\looseness -1 The vehicle was also run on all 28 closed loops of PPO on $\SE(3)$ with an actuator error in Table~\ref{tab:pert-all}. The actuator error distorts each commanded increment as in Section~\ref{sec:sims}, and the vehicle tracks the distorted increment. The settings are those above, with 200 episodes of 40 policy steps on each closed loop and initial errors drawn uniformly in the recomputed region. No episode leaves the region, on the closed loop or on the vehicle, and $V$ does not increase by more than $10^{-12}$ at any policy step. On the 15 closed loops whose level the second run of the search lowered, the initial errors were drawn in the region at the level of the first run, which contains the reported region. Since $V$ does not increase along these episodes, no episode that starts in the reported region leaves it either. Membership is checked at the policy steps only, every 2 s. The carried-out fraction, now the ratio of the achieved length to that of the distorted increment, has medians between $0.9999$ and $1.0018$, and a thruster was saturated in $0.27$ to $0.71$ percent of the integration steps. The same runs were also started in the region computed without the actuator error, which fails a condition under the actuator error on 10 of the 28 closed loops (Appendix~\ref{app:ablations}), at up to $0.26$ percent of its points outside $N_b$. On PPO seed 0 with misalignment $40^\circ$, one episode leaves that region, on the closed loop and on the vehicle, and $V$ increases at five policy steps, by up to $3.8$. On PPO seed 1 at $40^\circ$, $V$ increases at one policy step, by at most $0.62$, and no episode leaves the region. On the other 26 closed loops no episode leaves it and $V$ does not increase by more than $10^{-12}$.

\begin{table}[!ht]
\centering
\caption{Deployment: 200 episodes of 40 policy steps, on the closed loop $\Phi_u$ and on the simulated vehicle. Exits are episodes that leave $\Omega_\rho$. The largest increase of $V$ is over all policy steps on the vehicle. The carried-out fraction is the median ratio of the achieved to the commanded length of an increment, over steps with a commanded length above one percent of $\amax$. The final error is the median of $\|\Log\tilde g\|$ at the last step. The table covers the nominal setting.}
\label{tab:deploy}
\footnotesize
\begin{tabular}{llrrrrrr}
\toprule
vehicle & group & period & exits ($\Phi_u$) & exits (vehicle) & largest increase of $V$ & carried out & final error \\
\midrule
CubeSat & $\SO(3)$ & 15.7 s & 0/200 & 0/200 & $1.4\times10^{-22}$ & 0.936 & $7.6\times10^{-15}$ \\
BlueROV2 & $\SE(3)$ & 2.0 s & 0/200 & 0/200 & $8.2\times10^{-32}$ & 1.000 & $1.8\times10^{-16}$ \\
\bottomrule
\end{tabular}~
\end{table}

\end{document}